\documentclass[final,onefignum,onetabnum]{siamonline171218}

\usepackage{lipsum}
\usepackage{amsfonts}
\usepackage{graphicx}
\usepackage{mathpazo}
\usepackage{epstopdf}
\usepackage{algorithmic}
\usepackage{amsmath}
\usepackage{amssymb}
\ifpdf
  \DeclareGraphicsExtensions{.eps,.pdf,.png,.jpg}
\else
  \DeclareGraphicsExtensions{.eps}
\fi
\newcommand{\dd}{\mathrm{d}}
\newcommand{\ee}{\mathrm{e}}
\newcommand{\ii}{\mathrm{i}}
\renewcommand{\varepsilon}{\epsilon}
\usepackage{enumitem}
\setlist[enumerate]{leftmargin=.5in}
\setlist[itemize]{leftmargin=.5in}

\newsiamremark{remark}{Remark}
\newsiamremark{hypothesis}{Hypothesis}
\crefname{hypothesis}{Hypothesis}{Hypotheses}
\newsiamthm{claim}{Claim}
\newsiamthm{example}{Example}

\headers{Benjamin-Ono Soliton Ensembles}{E. Blackstone, L. Gassot, and P. D. Miller}

\title{On Strong Zero-Dispersion Asymptotics for Benjamin-Ono Soliton Ensembles%\thanks{Submitted to the editors DATE.
}

\author{Elliot Blackstone\thanks{Department of Mathematics, University of Michigan, Ann Arbor, MI 
  (\email{eblackst@umich.edu}).}
\and Louise Gassot\thanks{CNRS and Department of Mathematics, University of Rennes, Rennes, France
	(\email{louise.gassot@cnrs.fr}).}
\and Peter D. Miller\thanks{Department of Mathematics, University of Michigan, Ann Arbor, MI 
  (\email{millerpd@umich.edu}).}
}

\usepackage{amsopn}

\makeatletter
\newcommand*{\addFileDependency}[1]{% argument=file name and extension
  \typeout{(#1)}% latexmk will find this if $recorder=0 (however, in that case, it will ignore #1 if it is a .aux or .pdf file etc and it exists! if it doesn't exist, it will appear in the list of dependents regardless)
  \@addtofilelist{#1}% if you want it to appear in \listfiles, not really necessary and latexmk doesn't use this
  \IfFileExists{#1}{}{\typeout{No file #1.}}% latexmk will find this message if #1 doesn't exist (yet)
}
\makeatother

\newcommand*{\myexternaldocument}[1]{%
    \externaldocument{#1}%
    \addFileDependency{#1.tex}%
    \addFileDependency{#1.aux}%
}
\ifpdf
\hypersetup{
  pdftitle={On Strong Zero-Dispersion Asymptotics for Benjamin-Ono Soliton Ensembles},
  pdfauthor={Elliot Blackstone and Louise Gassot and Peter D. Miller}
}
\fi

\myexternaldocument{ex_supplement}

\usepackage[normalem]{ulem}
\usepackage{comment}
\usepackage{mathrsfs}
\usepackage{subcaption}
\usepackage{caption}
\usepackage{tikz}
\usetikzlibrary{arrows}
\usetikzlibrary{decorations.pathmorphing}
\usetikzlibrary{decorations.pathreplacing}
\usetikzlibrary{decorations.markings}
\usetikzlibrary{patterns}
\usetikzlibrary{automata}
\usetikzlibrary{positioning}
\usepackage{tikz-cd}

\tikzset{->-/.style={decoration={markings,mark=at position #1 with {\arrow[thick]{>}}},postaction={decorate}}}
	
\tikzset{-<-/.style={decoration={markings,mark=at position #1 with{\arrow[thick]{<}}},postaction={decorate}}}

\tikzset{test/.style n args={3}{
    postaction={
    decorate,
    decoration={
    markings,
    mark=between positions 0 and \pgfdecoratedpathlength step 0.5pt with {
    \pgfmathsetmacro\myval{multiply(
        divide(
        \pgfkeysvalueof{/pgf/decoration/mark info/distance from start}, \pgfdecoratedpathlength
        ),
        100
    )};
    \pgfsetfillcolor{#3!\myval!#2};
    \pgfpathcircle{\pgfpointorigin}{#1};
    \pgfusepath{fill};}
}}}}

\newtheorem{conjecture}{Conjecture}

\dedication{To the memory of Igor Krichever, for the enduring legacy of his scientific influence.}

\begin{document}

\maketitle

% REQUIRED
\begin{abstract}
A soliton ensemble is a particular kind of approximation of the solution of an initial-value problem for an integrable equation by a reflectionless potential that is well adapted to singular asymptotics like the small-dispersion limit.  We study soliton ensembles for the Benjamin-Ono equation by using reasonable hypotheses to develop local approximations that capture highly oscillatory features of the solution and hence provide more information than weak convergence results that are easier to obtain.  These local approximations are deduced independently from empirically-observed but unproven distributions of eigenvalues of two related matrices, one Hermitian and another non-Hermitian.  We perform careful numerical experiments to study the asymptotic behavior of the eigenvalues of these matrices in the small-dispersion limit, and formulate conjectures reflecting our observations.  Then we apply the conjectures to construct the local approximations of slowly varying profiles and rapidly oscillating profiles as well.  We show that the latter profiles are consistent with the predictions of Whitham modulation theory as originally developed for the Benjamin-Ono equation by Dobrokhotov and Krichever.
\end{abstract}

% REQUIRED

% REQUIRED

\section{Introduction}

The Benjamin-Ono (BO) equation \cite{Benjamin1967,Ono1975} is given by
\begin{align}\label{BO equation}
    u_{t}+2uu_{x}+\epsilon\mathcal{H}[u_{xx}]=0, ~~~ x\in\mathbb{R}, ~~~ t>0,
\end{align}
where $\epsilon\geq0$ is a parameter and $\mathcal{H}$ is the classical Hilbert transform\footnote{\color{black}This definition follows the ``physicist's'' sign convention, as used for instance by Matsuno in his textbook on the BO equation \cite[Eqn.\@ (3.2)]{MatsunoBook}.  However, more mathematical works (e.g., \cite{Gerard23}) define the Hilbert transform with the opposite sign.}
\begin{align}\label{def: Hilbert transform}
    \mathcal{H}[f](x)=\frac{1}{\pi}\mathrm{P. V.}\int_{-\infty}^{\infty}\frac{f(y)}{y-x}\,\dd y.
\end{align}
The \emph{zero-dispersion limit} refers to the analysis of the solution $u=u(x,t)$ of the Cauchy problem for \eqref{BO equation} where $u(x,0)=u_0(x)$ is an initial condition independent of $\epsilon$.  For convenience, we will assume that $u_0$ is an admissible initial condition in the sense of \cite[Definition 3.1]{MillerX11}.

\begin{figure}
\begin{center}
    \includegraphics[width=0.45\linewidth]{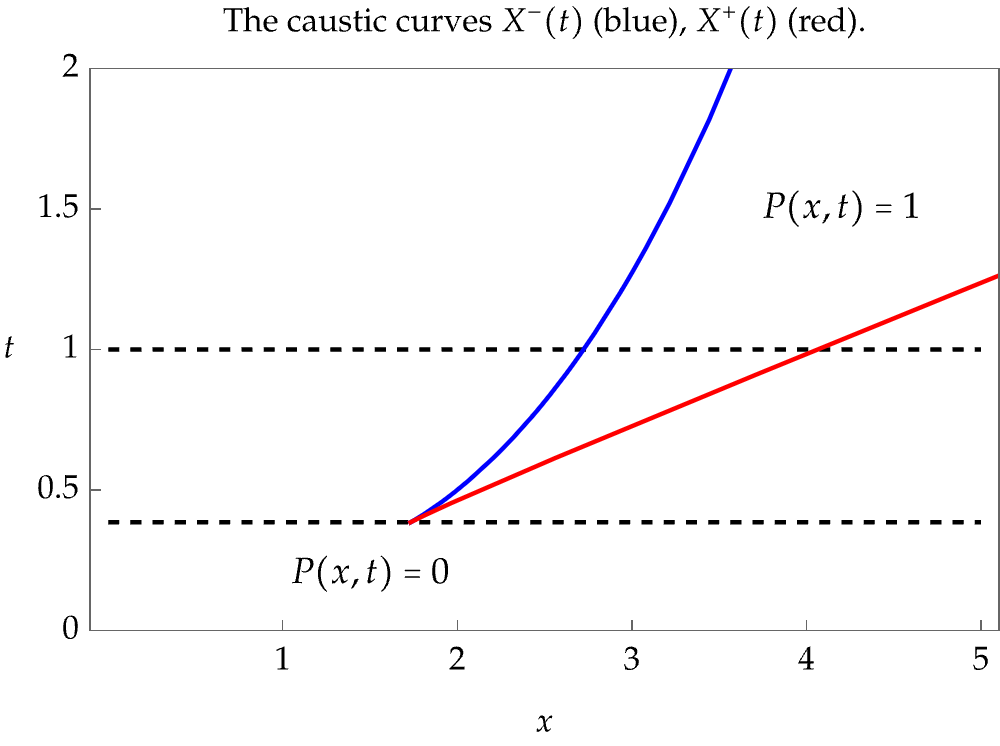}\hfill%
    \includegraphics[width=0.45\linewidth]{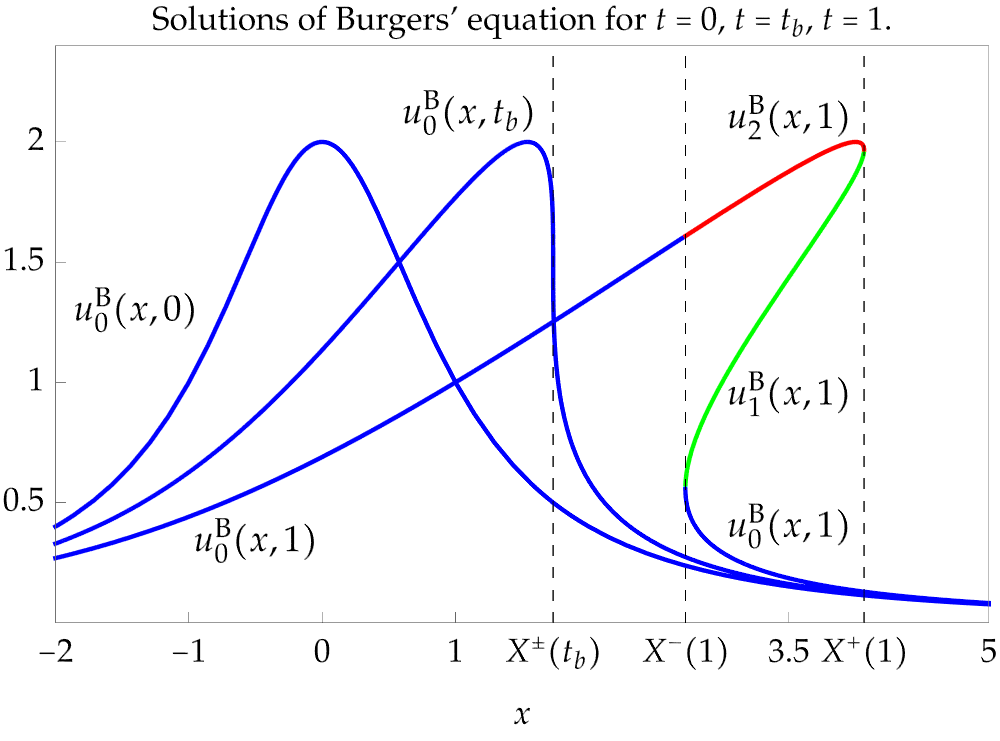}
\end{center}
\caption{The single and triple-valued region for Burgers' equation corresponding to $u_0(x)=2(1+x^2)^{-1}$ are separated by the caustic curves $X^{-}(t)$, $X^{+}(t)$ for $t\geq t_b$, see left panel.  The curves $X^{\pm}(t)$ are the double roots of \eqref{burgers sol}.  In the right panel, the solutions of Burgers' equation are plotted for $t=0,t_b,1$.  The blue, green, red curves are $u_0^{\mathrm{B}}(x,t)$, $u_1^{\mathrm{B}}(x,t)$, $u_2^{\mathrm{B}}(x,t)$, respectively.}
\label{fig: Burgers}
\end{figure}

The BO equation \eqref{BO equation} is an asymptotic model, derived in a small-amplitude and long-wave limit, for internal water waves propagating in one direction.  It applies to gravity-driven motions of the pycnocline separating a lower-density upper fluid layer from a higher-density lower fluid layer in the situation that the lower layer is assumed to be infinitely deep.  The solution $u(x,t)$ is a measure of the vertical displacement of the interface at position $x$ and time $t$.  The parameter $\epsilon>0$ is a measure of the relative strength of dispersion compared to nonlinear effects.  When $\epsilon\ll 1$, the local solution of the Cauchy problem for the inviscid Burgers equation
\begin{align}\label{eq: burgers}
    u^\mathrm{B}_t+2u^\mathrm{B}u^\mathrm{B}_x=0
\end{align}
with the same initial data $u^\mathrm{B}(x,0)=u_0(x)$ is expected to be a good approximation of $u(x,t)$ as long as $u^\mathrm{B}(x,t)$ remains smooth.  The solution of \eqref{eq: burgers} satisfying $u^{\mathrm{B}}(x,0)=u_0(x)$ is given implicitly by
\begin{align}\label{burgers sol}
    u^{\mathrm{B}}(x,t)=u_0(x-2tu^{\mathrm{B}}(x,t)).
\end{align}
However, for typical initial data $u_0$, a gradient catastrophe occurs in $u^\mathrm{B}(x,t)$ at a finite time $t=t_b$, beyond which the dispersion term in \eqref{BO equation} can no longer be neglected.  See Figure \ref{fig: Burgers}, right panel, for solutions of Burgers' equation before and after $t=t_b$.  Numerical experiments show that the effect of small dispersion is to generate a dispersive shock wave, that is, a train of waves with wavelength proportional to $\epsilon$ but with amplitude that is not small.  See Figure~\ref{fig:profile x2} below.  Such waves are described at a formal level by Whitham modulation theory \cite{Whitham}.  That theory starts from a family of exact solutions of (generally multiphase) waves parameterized by arbitrary constant amplitude, wavelength, and wave-shape parameters, and based on multiphase averaging of local conservation laws or a variational principle, posits a system of modulation equations which are partial differential equations governing slowly-varying fields replacing the constant parameters.

For the BO equation, nonlocality makes it challenging both to properly define the multiphase wave solutions and to determine their modulation equations.  Both of these problems were solved by Igor Krichever in a joint work with Dobrokhotov \cite{DobrokhotovK91} (see also the review by Krichever \cite{Krichever91}) that has had a large and lasting influence on the subject.  One of the remarkable observations in that work is that there exist Riemann invariant variables for the modulation equations in which they take the form of a collection of $2P+1$ independent copies of the inviscid Burgers equation \eqref{eq: burgers}, for some $P=0,1,2,\dots$ (the number of phases in the modulating wave), see Figure \ref{fig: Burgers}, left panel.  By contrast, in the Riemann-invariant form of the modulation equations for the Korteweg-de Vries (KdV) equation, one replaces the characteristic speed of a Riemann invariant by a complicated expression involving all of the invariants and represented as a ratio of hyperelliptic integrals, see \cite{FlaschkaFM80}.

One approach to describing the asymptotic behavior of $u(x,t)$ is to adopt a reasonable topology of convergence in which there is a limiting function, denoted $\overline{u}(x,t)$, as $\epsilon\to 0$.  Following the seminal work of Lax and Levermore \cite{LL} on a corresponding small-dispersion limit for the KdV equation, there has been some progress in proving convergence of $u(x,t)$ (or a suitable surrogate, see below) to a limit $\overline{u}(x,t)$ in the weak $L^2(\mathbb{R})$ topology with respect to~$x$, uniformly on compact intervals of $t$.  This has been done both in the setting of $x\in\mathbb{R}$ \cite{MillerX11} and on the torus (periodic boundary conditions) \cite{Gassot23a}.  The results of \cite{MillerX11} have recently been strengthened by G\'erard \cite{Gerard23} using a different approach.  In both the periodic and non-periodic cases, the formula for the weak limit $\overline{u}(x,t)$ is remarkably simple.  Indeed, let $u^\mathrm{B}=u^\mathrm{B}_n(x,t)$, $n=0,\dots,2P$, $u^\mathrm{B}_m(x,t)\le u^\mathrm{B}_n(x,t)$ for $m<n$, denote the generically distinct solutions of the implicit equation \eqref{burgers sol}, which correspond to the ``sheets'' above a given point $(x,t)\in\mathbb{R}^2$ of the multi-valued solution of the inviscid Burgers equation \eqref{eq: burgers} with initial data $u^\mathrm{B}(x,0)=u_0(x)$ (the number $P$ depends on $(x,t)$), see Figure \ref{fig: Burgers}.  Then the weak limit of $u(x,t)$ is given by the alternating sum
\begin{equation}
    \overline{u}(x,t)=\sum_{n=0}^{2P}(-1)^nu^\mathrm{B}_n(x,t).
    \label{eq:ubar}
\end{equation}
This result is far simpler than the corresponding result for KdV obtained in \cite{LL}.  The weak convergence of $u(x,t)$ can be extended \cite{MillerX12} to an infinite number of conserved ``local'' densities (the differential algebra of fields has to be augmented with Hilbert transforms), but even with this additional control, the topology of convergence is insufficient to capture the waveform and phase of the dispersive shock wave that forms in $u(x,t)$ for $t>t_b$.  Indeed, the wild oscillations are simply averaged out upon integration in $x$ against a test function in $L^2(\mathbb{R})$.  The formula \eqref{eq:ubar} suggests that the oscillations occupy the part of the $(x,t)$-plane where $P=P(x,t)\ge 1$, and numerical experiments support this assertion as well.   Moreover, the dispersive shock wave is expected to be described by the one-phase solution (assuming $P=1$).  See~\cite[Section 2.6]{KleinSaut21} for an overview of corresponding results in the context of the KdV equation.

The surrogate for $u(x,t)$ that was mentioned above is an approximation that we will call in this paper a \emph{soliton ensemble}.  It is a family denoted $\tilde{u}(x,t)$ of exact solutions of the BO equation \eqref{BO equation} associated to the specified function $x\mapsto u_0(x)$ via a systematic small-$\epsilon$ approximation of the scattering data for the BO Lax operator with potential $u_0$. In this approximation, one firstly replaces the true discrete eigenvalues $\lambda<0$ with approximate ones obtained from a type of quantization rule (see \eqref{def: tilde lambda n} below), and makes a similar approximation of the auxiliary phase constant $\gamma$ associated with each eigenvalue (see \eqref{def: tilde gamma n} below).  Secondly, one neglects the reflection coefficient defined for $\lambda>0$.  The exact solution of BO corresponding to this modified scattering data for each $\epsilon>0$ is precisely the function $\tilde{u}(x,t)$ (see \eqref{eq: def tilde u} below).  It is a ``nonlinear superposition'' of a large number, proportional to $\epsilon^{-1}$, of solitons, which combine coherently to yield an approximation of the given data $u_0$ in the strong $L^2$ sense when $t=0$, see \cite[Corollary 1.2]{MillerX11}.  The idea of using such a soliton ensemble to analyze the small-dispersion limit originated in the Lax-Levermore theory of the KdV equation \cite{LL}.  In both the BO and KdV cases, the soliton ensemble $\tilde{u}(x,t)$ is expressed explicitly in terms of finite determinants of size proportional to $\epsilon^{-1}$.  While soliton ensembles have proven to be useful to analyze the small-dispersion limit of the BO equation in the sense of weak convergence, in this paper we will study them from the point of view of strong convergence, aiming to capture the oscillatory profile of the dispersive shock wave.

Soliton ensembles may appear to be similar to soliton gases such as studied in \cite{GirottiGJM21}, as both are constructed from a large number of solitons. Soliton gases are limiting solutions in which the number of solitons is actually infinite while the limit process in their construction does not tie the number of solitons to a small parameter in the equation.  Moreover, soliton gases are not specifically constructed to match any particular initial condition $u_0$ at $t=0$.  Hence, one can perhaps think of a soliton ensemble as a type of soliton gas consisting of an arbitrarily large but finite number of solitons that is adapted to the small-dispersion limit and designed to approximate the solution of a specific Cauchy problem.

\subsection{Definition of soliton ensembles}
Define the Cauchy transforms
\begin{align}\label{def: Cauchy transform}
    \mathcal{C}_{\pm}[f](x):=\lim_{\delta\to0^{+}}\frac{1}{2\pi\ii}\int_{-\infty}^{\infty}\frac{f(y)}{y-(x\pm\ii\delta)}\,\dd y,
\end{align}
noting that $\pm\mathcal{C}_\pm$ are self-adjoint orthogonal projections from $L^2(\mathbb{R})$ onto the Hardy space 
\[
H^{\pm}(\mathbb{R})=\Big\{f \text{ holomorphic on } \mathbb{C}^{\pm}:\|f\|_{H^{\pm}}^{2}=\underset{y\in(0,\infty)}{\mathrm{sup}}\int_{\mathbb{R}}|f(x\pm\ii y)|^2\,\dd x<\infty\Big\},
\]
where $\mathbb{C}^{\pm}$ denotes the upper($+$)/lower($-$) half plane.  The inverse-scattering transform solution of the Cauchy problem for the BO equation \eqref{BO equation} is based on the self-adjoint Lax operator
\begin{align}
    \mathcal{L}:=-\ii\epsilon\frac{\partial}{\partial x}-\mathcal{C}_+\circ u\circ\mathcal{C}_+
\end{align}
acting on a domain dense in $H^+(\mathbb{R})$, wherein $u$ denotes the operator of multiplication by a bounded function $u:\mathbb{R}\to\mathbb{R}$.  It has been proved \cite{Wu16} that if $u\in L^1(\mathbb{R})\cap L^\infty(\mathbb{R})$ and $x\mapsto xu(x)$ is in $L^2(\mathbb{R})$, then $\mathcal{L}$ has only finitely many eigenvalues, all negative real numbers $\lambda_1<\lambda_2<\cdots <\lambda_N<0$ with unit geometric multiplicity.  To each eigenvalue $\lambda_j<0$, there corresponds a real number called a phase constant, denoted $\gamma_j$.  Under suitable conditions on $u$ \cite{Wu17} there is also a complex-valued function $\lambda\mapsto \beta(\lambda)$ defined for $\lambda>0$ called the reflection coefficient.  The collection of $N$ pairs $(\lambda_j,\gamma_j)$ and the function $\beta$ are said to constitute the scattering data associated with the initial data $u=u_0$ for \eqref{BO equation} in the inverse-scattering transform first proposed by Fokas and Ablowitz \cite{FokasA83}.  The inversion of the transform has not yet been fully justified, but if $\beta(\lambda)\equiv 0$, the procedure reduces to finite-dimensional linear algebra, resulting in the formula
\begin{align}\label{u Matsuno}
    u(x,t)=2\epsilon\frac{\partial}{\partial x}\mathrm{Im}\left[\log\left(\det\left(\mathbb{I}+\frac{\ii}{\epsilon}\mathbf{A}(x,t)\right)\right)\right],
\end{align}
valid for any choice of branch of the complex logarithm,
where $\mathbb{I}$ denotes the $N\times N$ identity matrix, and $\mathbf{A}(x,t)$ is the $N\times N$ Hermitian matrix with elements
\begin{align}\label{eq: def A}
    A_{jk}(x,t)
    =\begin{cases}
    \displaystyle \frac{2\ii\epsilon\sqrt{\lambda_j\lambda_k}}{\lambda_j-\lambda_k}, &j\neq k, \\
    -2\lambda_j(x+2\lambda_jt+\gamma_j), &j=k.
    \end{cases}
\end{align}
This formula represents a multi-soliton solution of \eqref{BO equation}, and it was first derived by Matsuno \cite{Matsuno1979} using Hirota's bilinear method. The inversion formula has then been proved on the $N$-soliton manifolds in~\cite[Eq.\@ (1.19)]{Sun21}. 

In two papers \cite{Matsuno81,Matsuno82}, Matsuno used formal arguments to investigate the asymptotic behavior of the scattering data described above for fixed $u_0$ in the limit $\epsilon\to 0$.  The main results of \cite{Matsuno81,Matsuno82} were asymptotic formul\ae\ for $|\beta(\lambda)|^2$ for $\lambda>0$ and for the distribution of eigenvalues $\lambda_j<0$, of which there are in general a large number $N\sim\epsilon^{-1}$.  These results were rigorously proven in \cite{MillerW16b} under the additional assumption that $u_0$ is a rational function for which $u_0(x)=-\lambda$ has generically either two or zero real solutions $x$ given $\lambda\in\mathbb{R}$, and for such $u_0$ an asymptotic formula for the phase constant $\gamma_j$ associated with a given eigenvalue $\lambda_j$ was also rigorously established.  We now restate these asymptotic results in the special case that $u_0$ is a smooth positive function with a single critical point (the maximizer) and with sufficient decay as $x\to\pm\infty$.  For such $u_0$, we define positive constants $L,M>0$ by 
\begin{equation}
    L:=\max_{x\in\mathbb{R}} u_0(x),\quad M:=\frac{1}{2\pi}\int_\mathbb{R} u_0(x)\,\dd x,
    \label{eq:LM}
\end{equation}
and for $-L<\lambda<0$ define the turning points $x=x_\pm(\lambda)$, $x_-(\lambda)<x_+(\lambda)$ as the two roots of $u_0(x)=-\lambda$.  Then the following are true:
\begin{itemize}
    \item the reflection coefficient $\beta(\lambda)$ defined for $\lambda>0$ vanishes in the limit $\epsilon\to 0$;
    \item there are $N=M/\epsilon + \mathcal{O}(1)$ eigenvalues in the interval $-L<\lambda<0$, and the number of eigenvalues with $-L<a<\lambda<b<0$ is $N(a,b)$ satisfying
    \begin{equation}\label{eq:number-solitons}
        N(a,b)=\frac{1}{\epsilon}\int_a^b F(\lambda)\,\dd\lambda +\mathcal{O}(1),
    \end{equation}
    where the density of eigenvalues near $\lambda$ is $\epsilon^{-1}F(\lambda)$ with
    \begin{equation}
        F(\lambda):=\frac{1}{2\pi}(x_+(\lambda)-x_-(\lambda)),\quad -L<\lambda<0;
        \label{eq:density of eigenvalues}
    \end{equation}
    \item if $\lambda_j$ is an eigenvalue that converges to a number $\lambda\in (-L,0)$ as $\epsilon\to 0$, then also $\gamma_j\to\gamma(\lambda)$, where
    \begin{equation}
        \gamma(\lambda):=-\frac{1}{2}(x_+(\lambda)+x_-(\lambda)),\quad -L<\lambda<0.
        \label{eq:gamma-formula}
    \end{equation}
\end{itemize}
In the periodic case, analogues of formulas~\eqref{eq:number-solitons}, \eqref{eq:density of eigenvalues} and~\eqref{eq:gamma-formula} should also hold~\cite{Moll20}, and have been proved for bell-shaped initial data in~\cite{Gassot23b}.

\begin{figure}
     \centering
     \begin{subfigure}[b]{0.34\textwidth}
         \centering
         \includegraphics[width=\textwidth]{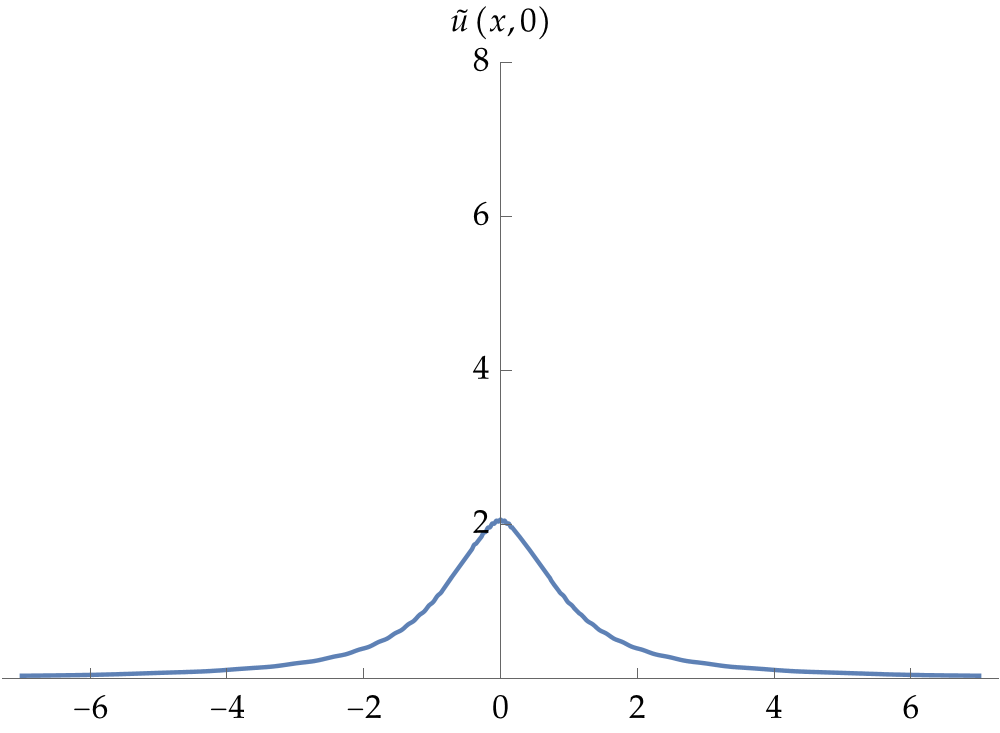}
     \end{subfigure}
     %\hfill
     \hspace{-0.5cm}
     \begin{subfigure}[b]{0.34\textwidth}
         \centering
         \includegraphics[width=\textwidth]{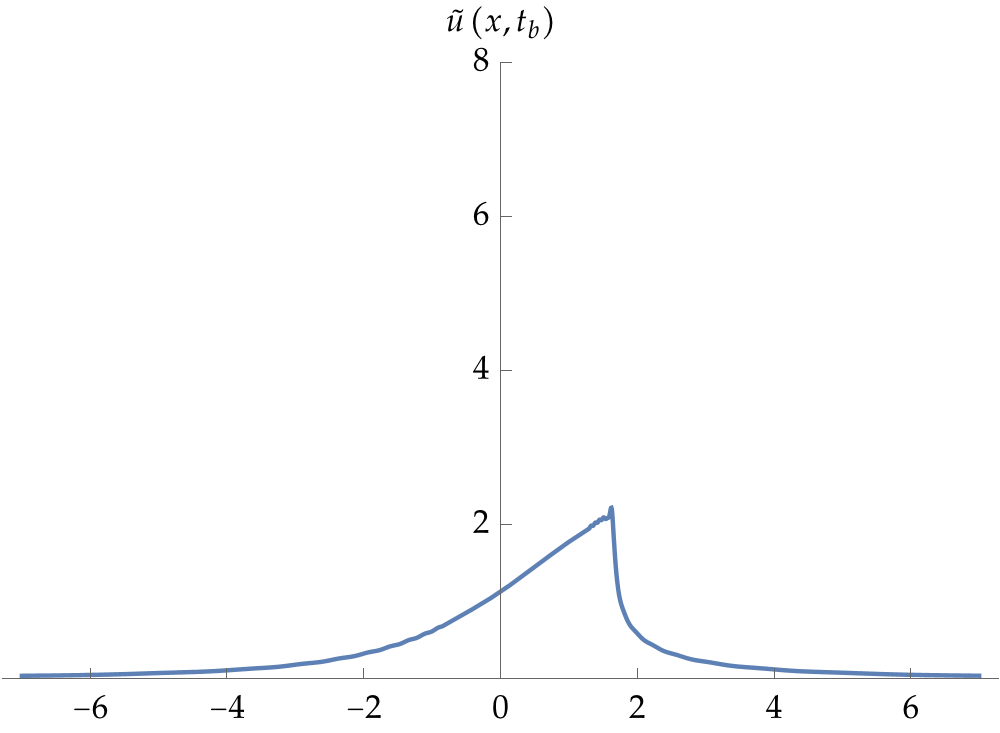}
     \end{subfigure}
     %\hfill
     \hspace{-0.5cm}
     \begin{subfigure}[b]{0.34\textwidth}
         \centering
         \includegraphics[width=\textwidth]{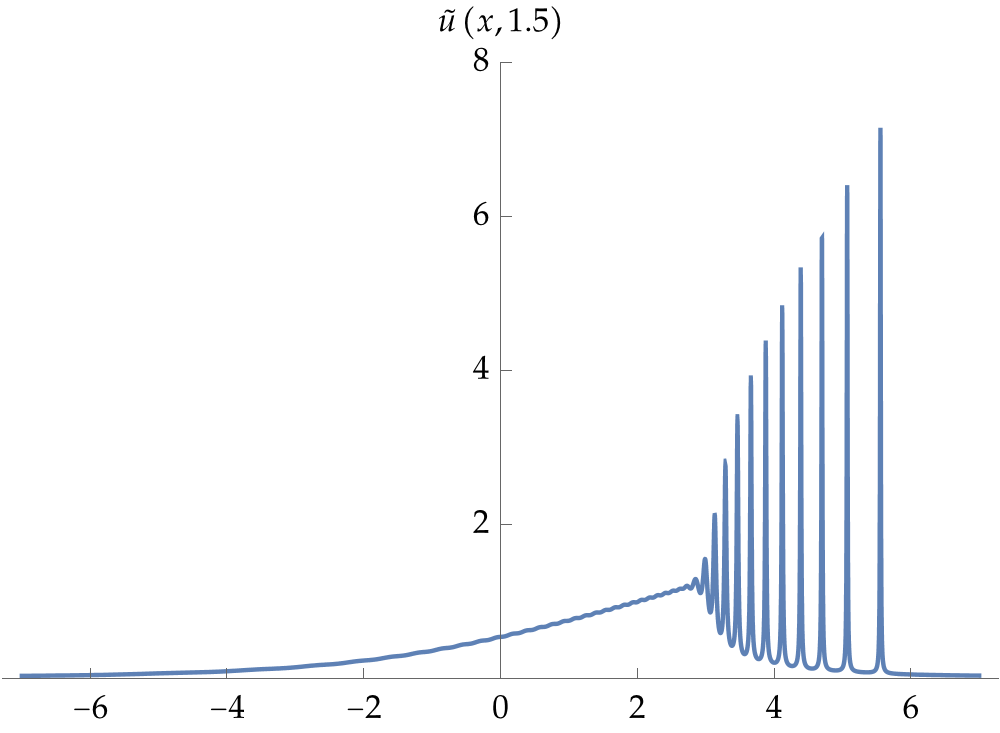}
     \end{subfigure}
        \caption{Time evolution of $\tilde{u}(x,t)$ for $u_0(x)=2(1+x^2)^{-1}$ and $\epsilon=0.02$.}
        \label{fig:profile x2}
\end{figure}

Based on these results, we will now define the soliton ensemble for the BO equation \eqref{BO equation} associated with a positive initial condition $u_0$ of the type described above.  Firstly, we define the exact number of approximate eigenvalues by setting 
\begin{equation}
    N(\epsilon):=\left\lfloor\frac{M}{\epsilon}\right\rfloor.
\end{equation}
Then, we define $N(\epsilon)$ approximate eigenvalues $\tilde{\lambda}_j\in (-L,0)$, $j=1,\dots,N(\epsilon)$ by quantizing the density formula \eqref{eq:density of eigenvalues}:
\begin{equation}\label{def: tilde lambda n}
    \int_{-L}^{\tilde{\lambda}_j}F(\lambda)\,\dd\lambda = \epsilon\left(j-\frac{1}{2}\right),\quad j=1,\dots,N(\epsilon).
\end{equation}
Finally, we define corresponding phase constants by setting
\begin{equation}
    \tilde{\gamma}_j:=\gamma(\tilde{\lambda}_j),\quad j=1,\dots,N(\epsilon),
    \label{def: tilde gamma n}
\end{equation}
where $\gamma(\cdot)$ is defined by \eqref{eq:gamma-formula}.  Then we neglect the reflection coefficient and define an exact multi-soliton solution of the BO equation \eqref{BO equation} using the approximate discrete data.  In detail, first  
define the elements of the $N(\epsilon)\times N(\epsilon)$ Hermitian matrix $\tilde{\mathbf{A}}(x,t)$ as (compare with \eqref{eq: def A}) 
\begin{align}
    \tilde{A}_{jk}(x,t)=\begin{cases} \displaystyle
    \frac{2\ii\epsilon\sqrt{\tilde{\lambda}_j\tilde{\lambda}_k}}{\tilde{\lambda}_j-\tilde{\lambda}_k}, &j\neq k, \\
    -2\tilde{\lambda}_j(x+2\tilde{\lambda}_jt+\tilde{\gamma}_j), &j=k.
    \end{cases}
\end{align}
Then set (compare with \eqref{u Matsuno})
\begin{equation}\label{eq: def tilde u}
\tilde{u}(x,t):=2\epsilon\frac{\partial}{\partial x}\mathrm{Im}\left[\log\left(\det\left(\mathbb{I}+\frac{\ii}{\epsilon}\tilde{\mathbf{A}}(x,t)\right)\right)\right].
\end{equation}
The family of functions $\tilde{u}(x,t)$ parameterized by $\epsilon>0$ constitute the soliton ensemble for the BO equation \eqref{BO equation} for the Cauchy data $u_0$.  See Figure~\ref{fig:profile x2}.

Note that \eqref{def: tilde lambda n} implies that 
\begin{equation}
    \int_{-L}^{\tilde{\lambda}_1}F(\lambda)\,\dd\lambda = \frac{\epsilon}{2}.
\end{equation}
So in the sense measured by differences of the left-hand side of \eqref{def: tilde lambda n}, the first approximate eigenvalue is half the distance to $-L$ as the remaining approximate eigenvalues are to their nearest neighbors.  For an argument appearing below in Section~\ref{sec:A}, we will want to ensure a similar condition for the last approximate eigenvalue, namely that
\begin{equation}
    \int_{\tilde{\lambda}_{N(\epsilon)}}^0 F(\lambda)\,\dd\lambda = \frac{\epsilon}{2}.
\label{eq:result-of-epsilon-quantization}
\end{equation}
It is easy to see that this condition holds precisely when $\epsilon=\epsilon_N:=M/N$ for an integer $N\in\mathbb{Z}_{>0}$, in which case \eqref{eq:number-solitons} returns $N(\epsilon_N)=N$.  For convenience, we will therefore assume below that $\epsilon$ tends to zero within this specific discrete sequence.

Also, since the rest of this paper is concerned only with the $\epsilon$-dependent function $\tilde{u}(x,t)$, henceforth we will drop all tildes for simplicity.

\subsection{Alternate formul\ae\ for the soliton ensemble \texorpdfstring{$u(x,t)$}{u(x,t)}}
Let $\alpha_k(x,t)$, $k=1,\ldots,N(\epsilon)$, denote the real eigenvalues of $\mathbf{A}(x,t)$.  The proof of weak convergence of $u(x,t)$ to $\overline{u}(x,t)$ given in \cite{MillerX11} is based on expressing $u(x,t)$ in the form
\begin{equation}
    u(x,t)=\frac{\partial I}{\partial x}(x,t),\;
    I(x,t):=2\epsilon\mathrm{Im}\left[\log\left(\det\left(\mathbb{I}+\frac{\ii}{\epsilon}\mathbf{A}(x,t)\right)\right)\right] = 2\epsilon\sum_{k=1}^{N(\epsilon)}\arctan(\epsilon^{-1}\alpha_k(x,t)).
    \label{eq:tilde-u-weak}
\end{equation}
Here to define $I(x,t)$ properly, we resolve the ambiguity of the logarithm by taking the principal branch of the arctangent with values in $(-\tfrac{1}{2}\pi,\tfrac{1}{2}\pi)$ in the sum.  The quantity $I(x,t)$ then resembles a Riemann sum for an integral of a discontinuous integrand $\arctan(\epsilon^{-1}\alpha)\to \frac{1}{2}\pi\mathrm{sgn}(\alpha)$ against the distribution of eigenvalues of $\mathbf{A}(x,t)$, and in \cite{MillerX11} this $(x,t)$-dependent distribution is calculated and used to prove locally uniform convergence of $I(x,t)$ to an antiderivative of $\overline{u}(x,t)$.  The weak nature of the convergence of $u(x,t)$ to $\overline{u}(x,t)$ can then be attributed to the $x$-derivative in \eqref{eq:tilde-u-weak}.

Our aim in this article is to build upon the method of \cite{MillerX11} by observing that if it is desired to improve the nature of the convergence of $u(x,t)$ as $\epsilon\to 0$, one should differentiate first in \eqref{eq:tilde-u-weak} and only then analyze $u(x,t)$.  To this end, explicitly differentiating $I(x,t)$ with respect to $x$ in \eqref{eq:tilde-u-weak} (assuming differentiability of the eigenvalues $\alpha_k(x,t)$), we obtain
\begin{equation}
    u(x,t)=\sum_{k=1}^{N(\epsilon)}\frac{2\epsilon^2\alpha_{k,x}(x,t)}{\alpha_k(x,t)^2 + \epsilon^2}, \quad \alpha_{k,x}(x,t):=\frac{\partial\alpha_k}{\partial x}(x,t).
\label{u in terms of alphas}
\end{equation}

Rather than diagonalizing $\mathbf{A}(x,t)=\mathbf{A}(x,t)^\dagger$ first and then differentiating, one could try to take advantage of the fact that the matrix elements of $\mathbf{A}(x,t)$ depend on $x$ in a very simple way.  Indeed, notice that
\begin{equation}
    \mathbb{I}+\ii\epsilon^{-1}\mathbf{A}(x,t)=\ii\epsilon^{-1}\mathbf{D}\left[x\mathbb{I} - \mathbf{B}(t) -\ii\epsilon\mathbf{D}^{-2} \right]\mathbf{D},
    \label{eq:A-to-B}
\end{equation}
where $\mathbf{D}$ is a Hermitian $N(\epsilon)\times N(\epsilon)$ diagonal matrix given by
\begin{equation}\label{def: D matrix}
\mathbf{D}:=\mathrm{diag}\left(\sqrt{-2\lambda_1},\dots,\sqrt{-2\lambda_{N(\epsilon)}}\right)
\end{equation}
and $\mathbf{B}(t)$ is a Hermitian $N(\epsilon)\times N(\epsilon)$ matrix with elements
\begin{equation}\label{def: B matrix}
    B_{jk}(t)=\begin{cases} \displaystyle 
    \frac{-\ii\epsilon}{\lambda_j-\lambda_k}, &j\neq k, \\
    -2\lambda_jt-\gamma_j, &j=k.
    \end{cases}
\end{equation}
Thus, for any branch of the complex logarithm,
\begin{equation}
    \frac{\partial}{\partial x}\log\left(\det\left(\mathbb{I}+\frac{\ii}{\epsilon}\mathbf{A}(x,t)\right)\right) = \frac{\partial}{\partial x}\log\left(\det\left(x\mathbb{I}-\mathbf{B}(t)-\ii\epsilon\mathbf{D}^{-2}\right)\right).
\end{equation}
If we let $\sigma_k(t)=\mu_k(t)+\ii\nu_k(t)$, where $\mu_k(t), \nu_k(t)\in\mathbb{R}$, denote the complex eigenvalues of the non-Hermitian but $x$-independent matrix 
\begin{equation}
    \mathbf{C}(t):=\mathbf{B}(t)+\ii\epsilon\mathbf{D}^{-2}, 
\label{eq:tilde-C-def}
\end{equation}
then \eqref{eq: def tilde u} can be written as
\begin{equation}\label{u in terms of sigmas}
\begin{split}
    u(x,t)&=2\epsilon\mathrm{Im}\left[\frac{\partial}{\partial x}\log(\det(x\mathbb{I}-\mathbf{B}(t)-\ii\epsilon\mathbf{D}^{-2})\right] \\
    &= 2\epsilon\mathrm{Im}\left[\frac{\partial}{\partial x}\sum_{k=1}^{N(\epsilon)}\log(x-\sigma_k(t))\right]\\
    &=\mathrm{Im}\left[\sum_{k=1}^{N(\epsilon)}\frac{2\epsilon}{x-\sigma_k(t)}\right]\\
    &=\sum_{k=1}^{N(\epsilon)}\frac{2\epsilon\nu_k(t)}{(x-\mu_k(t))^2+\nu_k(t)^2}.
    \end{split}
\end{equation}
By substituting the formula \eqref{u in terms of sigmas} into the BO equation \eqref{BO equation} and using the fact that
\begin{align*}
    \mathcal{H}[u(\cdot,t)](x)=\sum_{k=1}^{N(\epsilon)}\frac{-2\epsilon(x-\mu_k(t))}{(x-\mu_k(t))^2+\nu_k(t)^2},
\end{align*}
it can be verified that the $\sigma_k(t)$'s satisfy the system of differential equations
\begin{align}
    \sigma_k'(t)=\sum_{j=1}^{N(\epsilon)}\frac{2\ii\epsilon}{\sigma_k(t)-\sigma_j(t)^*}-\sum_{\substack{j=1 \\ j\neq k}}^{N(\epsilon)}\frac{2\ii\epsilon}{\sigma_k(t)-\sigma_j(t)}.
    \label{eq:pre-CM}
\end{align}
This was first discovered in \cite[Eq.\@ (7)]{CLP1979}, wherein \eqref{eq:pre-CM} was shown to be equivalent to the famous Calogero-Moser $N$-body system.  Let $\mathbf{w}_k(t)=(w_{k,1},\dots,w_{k,N(\epsilon)})^\top$ be the normalized eigenvector of $\mathbf{C}(t)$ corresponding to the eigenvalue $\sigma_k(t)$.  Clearly
\begin{align}\label{sigma in terms of e vec}
    \sigma_k(t)=\mathbf{w}_k(t)^\dagger\mathbf{C}(t)\mathbf{w}_k(t)=\mathbf{w}_k(t)^{\dag}\left(\mathbf{B}(t)+\ii\epsilon\mathbf{D}^{-2}\right)\mathbf{w}_k(t),
\end{align}
where $\dag$ denotes the conjugate transpose.  By comparing real and imaginary parts of \eqref{sigma in terms of e vec}, we have 
\begin{align}
    \mu_k(t)&=-\sum_{j=1}^{N(\epsilon)}(2t\lambda_j+\gamma(\lambda_j))|w_{k,j}(t)|^2-\ii\epsilon\sum_{j=1}^{N(\epsilon)}\sum_{\substack{s=1 \\ s\neq j}}^{N(\epsilon)}\frac{w_{k,j}(t)^*w_{k,s}(t)}{\lambda_s-\lambda_j},\\ 
    \nu_k(t)&=-\frac{\epsilon}{2}\sum_{j=1}^{N(\epsilon)}\frac{|w_{k,j}(t)|^2}{\lambda_j}.\label{eq:sign-nu}
\end{align}
Since $-L<\lambda_1<\cdots<\lambda_j<\cdots<\lambda_{N(\epsilon)}<0$, using $|w_{k,1}(t)|^2+\cdots+|w_{k,N(\epsilon)}(t)|^2=1$ we obtain
\begin{equation}
0<    \frac{1}{2}L\cdot\epsilon<\nu_k(t)<\frac{1}{2}\cdot\frac{\epsilon}{|\lambda_{N(\epsilon)}|},\quad k=1,\dots,N(\epsilon).
\label{eq:lower-upper-bounds-for-nu}
\end{equation}
Likewise, using also $|a||b|\le\frac{1}{2}(a^2+b^2)$,
\begin{equation}
    |\mu_k(t)|\le \left(2tL+\sup_{-L<\lambda<0}|\gamma(\lambda)|\right) +\frac{\epsilon N(\epsilon)}{\min_{s\neq j}|\lambda_s-\lambda_j|},\quad k=1,\dots,N(\epsilon).
    \label{eq:upper-bounds-for-mu}
\end{equation}
Below in \eqref{eq:smallest-lambdas} it is shown that if $u_0(x)$ decays algebraically as $x\to\pm\infty$, i.e., $u_0(x)\sim Cx^{-2p}$ for some $p\ge 1$, then assuming that $\epsilon=\epsilon_N=M/N$ so that \eqref{eq:result-of-epsilon-quantization} holds one finds that $|\lambda_{N(\epsilon)}|\gtrsim \epsilon^{2p/(2p-1)}$ and that $\min_{s\neq j}|\lambda_s-\lambda_j|\gtrsim\epsilon^{2p/(2p-1)}$, in which case the imaginary parts of the eigenvalues have uniform bounds of the form 
\begin{equation}
\frac{1}{2}L\cdot\epsilon<\nu_k(t)\lesssim\epsilon^{-1/(2p-1)}, \quad
k=1,\dots,N(\epsilon),
\label{eq:rational-nu-bound}
\end{equation}
which is equivalent to the rescaled form
\begin{equation}
    \frac{1}{2}L\cdot\frac{1}{\ln(\epsilon^{-1})}<\frac{\nu_k(t)}{\epsilon\ln(\epsilon^{-1})}\lesssim\frac{\epsilon^{-2p/(2p-1)}}{\ln(\epsilon^{-1})},\quad k=1,\dots,N(\epsilon),
    \label{eq:rational-nu-bound-rescaled}
\end{equation}
and similarly assuming also that $\gamma(\cdot)$ is a bounded function (it vanishes identically if $u_0(\cdot)$ is even),
\begin{equation}
    |\mu_k(t)|\lesssim \epsilon^{-2p/(2p-1)},\quad k=1,\dots,N(\epsilon).
    \label{eq:rational-mu-bound}
\end{equation}

\subsection{Outline of the paper}
The purpose of this work is to report the results of several numerical experiments undertaken to study the small-$\epsilon$ asymptotic behavior of the real eigenvalues $\{\alpha_k(x,t)\}$ and the complex eigenvalues $\{\sigma_k(t)\}$ of the matrices $\mathbf{A}(x,t)$ and $\mathbf{C}(t)$ respectively.  We recall that $u(x,0)=u_0(x)$ is an admissible initial condition (see \cite[Definition 3.1]{MillerX11}). For simplicity, we assume further that $u_0$ has exactly one inflection point to the right of the maximizer. Consequently, we always have that for generic $(x,t)$ with $t>0$, $P(x,t)=0$ or $P(x,t)=1$. Moreover, for fixed $t>t_b$ the multi-valued region (in $x$) for Burgers' equation is an interval $(X^{-}(t),X^{+}(t))$, see Figure \ref{fig: Burgers}.
For simplicity, all of our numerical experiments will be for the soliton ensemble associated with the initial condition 
\begin{equation}
    u_0(x)=\frac{2}{1+x^2}.
    \label{eq:Lorentzian}
\end{equation}
For this initial condition, the breaking time for the inviscid Burgers (dispersionless, i.e., $\epsilon=0$) approximation of the BO equation \eqref{BO equation} is exactly $t=t_b=2\sqrt{3}/9\approx 0.3849$.  When $t>t_b$, the dispersive terms are expected to be important and form a highly-oscillatory dispersive shock wave (see Figure \ref{fig:profile x2}, right panel) in the $t$-dependent interval of $x$ on which the method of characteristics predicts a triple-valued solution of the dispersionless approximation.

As will be seen, the numerical experiments suggest that when $\epsilon$ is small, the complex eigenvalues $\{\sigma_k(t)\}$ and the real eigenvalues $\{\alpha_k(x,t)\}$ that are so small as to contribute substantially to the sum in \eqref{u in terms of alphas} are evidently distributed in a regular fashion.  After formulating reasonable conjectures based on the numerical observations, we prove that they imply many of the properties of the dispersive shock wave.  These properties go beyond the weak convergence result that $u(x,t)\rightharpoonup\overline{u}(x,t)$ with limit $\overline{u}(x,t)$ given by \eqref{eq:ubar}.  On the other hand, in general the soliton ensemble only approximates $u_0$ in the strong $L^2(\mathbb{R})$ sense when $t=0$, and translation in $x$ is continuous on $L^2(\mathbb{R})$.  Therefore it is possible that $u(\cdot,0)$ is nearly a translation of $u_0(\cdot)$ by a fraction of the wavelength of the dispersive shock wave that appears for $t>t_b$, leading to $\mathcal{O}(1)$ errors for the solution of the Cauchy problem when measured in $L^\infty(\mathbb{R})$.  In this sense, our main results below, Theorems~\ref{thm:oscillations} and \ref{prop: approx after breaking via A}, apply to the soliton ensemble but would require phase corrections to apply also to the solution of the Cauchy problem for \eqref{BO equation} with initial data $u_0$.

Our study of the complex eigenvalues of $\mathbf{C}(t)$ and how their asymptotic properties imply the most important features of $u(x,t)$ via \eqref{u in terms of sigmas} is presented in Section~\ref{sec:C}.  Then in Section~\ref{sec:A} we give a parallel analysis logically independent of Section~\ref{sec:C} for the eigenvalues of $\mathbf{A}(x,t)$ and the formula \eqref{u in terms of alphas}.  
A completely different approach to strong asymptotics for the BO equation \eqref{BO equation} that is based instead on a remarkable formula of G\'erard \cite{GerardExplicit} is the subject of a forthcoming work \cite{BlackstoneGGM}.

To end this introduction, we emphasize that this paper falls into the category of experimental mathematics.  Our results are not entirely definitive as they are partly based on numerical experiments.  However in our study of this topic we have uncovered numerous interesting connections with traditional areas of mathematical analysis that we want to report, in hopes of stimulating further research.  Indeed, different parts of this paper are related to the theory of singular asymptotics of dispersive nonlinear waves, to the spectra of large structured matrices, and to geometric quantization and semiclassical analysis.

\section{Asymptotic properties of the complex eigenvalues \texorpdfstring{$\sigma_k(t)$}{sigma-k} and their implications}
\label{sec:C}

This section of the paper concerns the nonhermitian matrix $\mathbf{C}(t)$.  We present numerical experiments leading to conjectures about the distribution of its eigenvalues, and then we show how those conjectures lead to strong asymptotics of $u(x,t)$ consistent with Whitham modulation theory.

\subsection{Numerical experiments}
\begin{figure}
     \centering
     \begin{subfigure}[b]{0.32\textwidth}
         \centering
         \includegraphics[width=\textwidth]{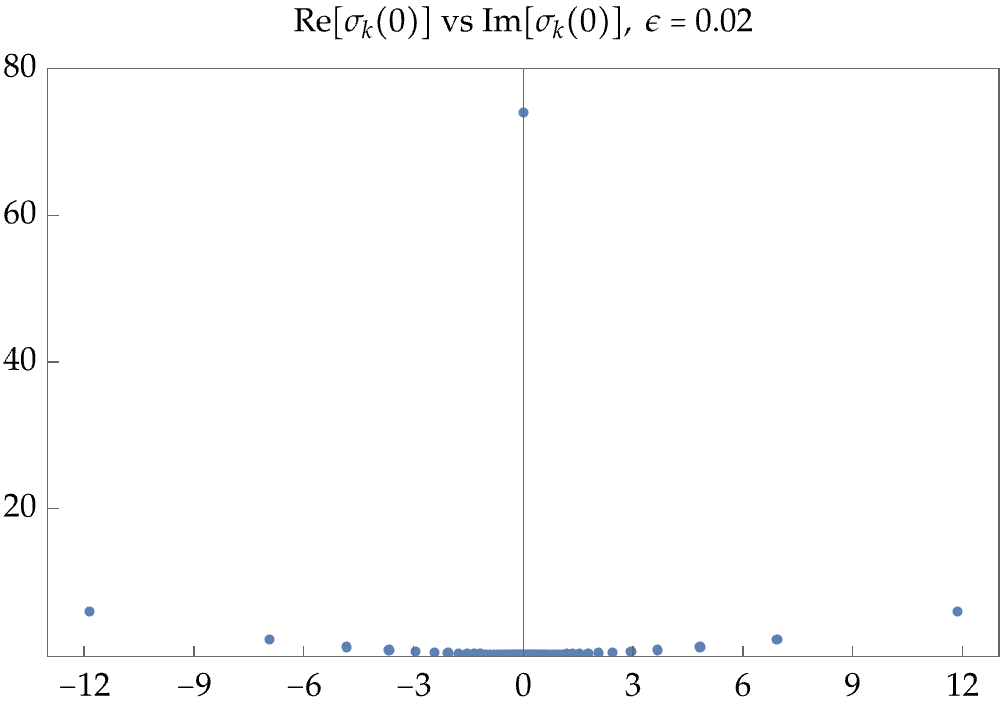}
     \end{subfigure}
     %\hfill
     \begin{subfigure}[b]{0.32\textwidth}
         \centering
         \includegraphics[width=\textwidth]{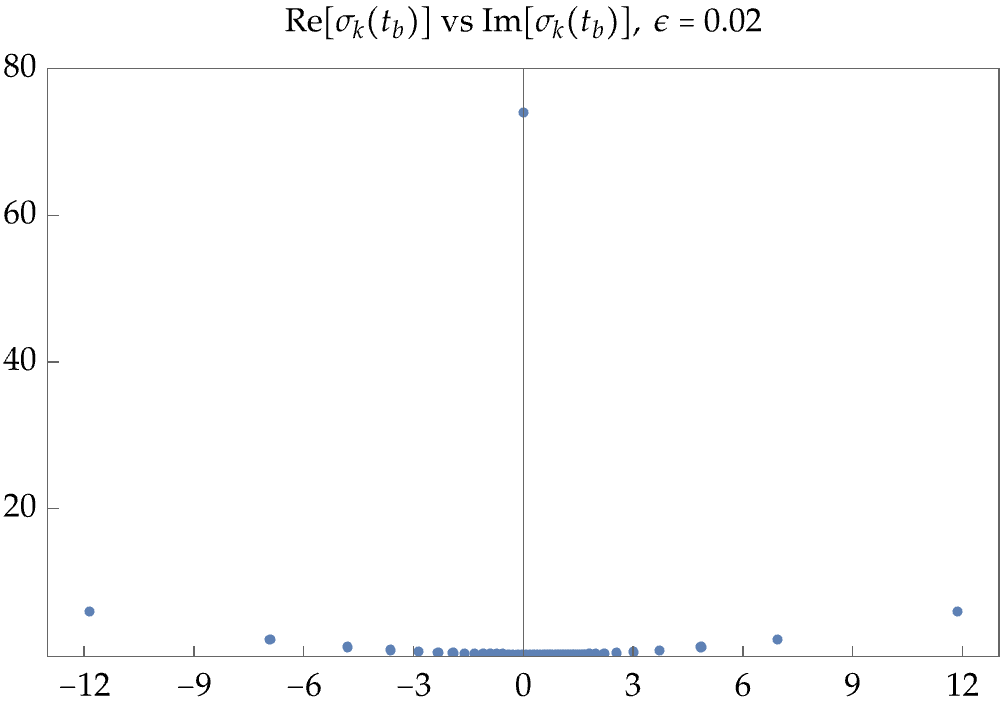}
     \end{subfigure}
     %\hfill
     \begin{subfigure}[b]{0.32\textwidth}
         \centering
         \includegraphics[width=\textwidth]{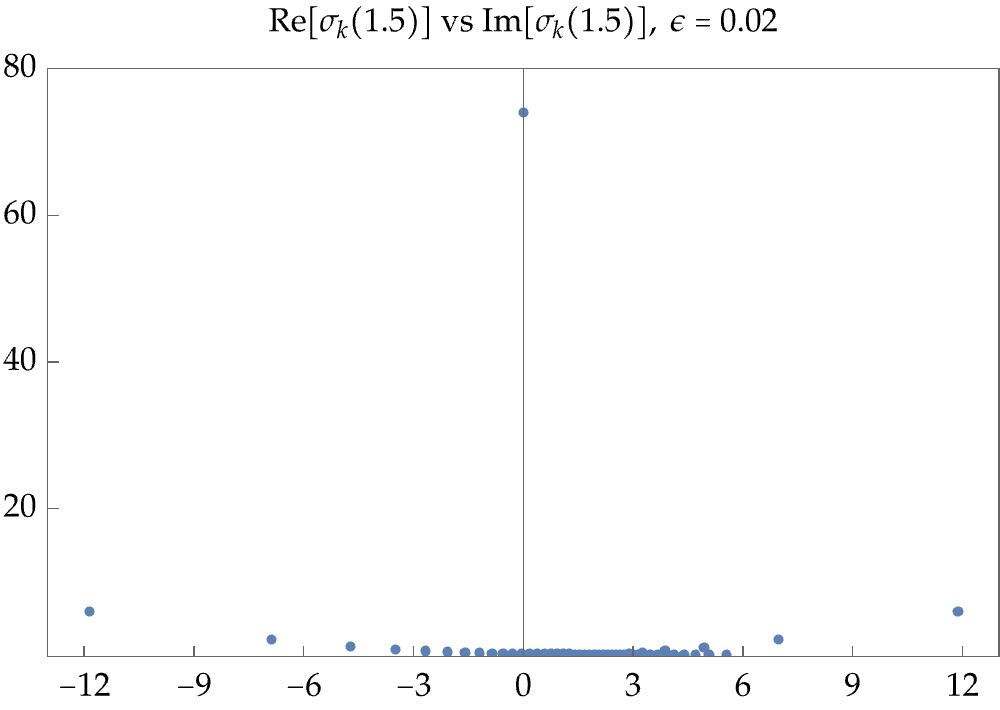}
     \end{subfigure}
        \caption{Outlier eigenvalues for the initial condition $u_0(x)=2(1+x^2)^{-1}$ with $\epsilon=0.02$.}
    \label{fig:outliers x2}
\end{figure}

Given an initial condition $u_0(\cdot)$, a time $t\ge 0$, and a value of $\epsilon$, it is straightforward to construct the nonhermitian matrix $\mathbf{C}(t)$ and to numerically extract its eigenvalues $\{\sigma_k(t)\}$ with high accuracy.  The first observation is that, roughly speaking, most of the complex eigenvalues are close to the real axis.  However depending on the initial condition, there can be some ``outliers'' that do not follow this rule, as shown in Figure~\ref{fig:outliers x2}.

\begin{figure}[h]
     \centering
     \begin{subfigure}[b]{0.32\textwidth}
         \centering
         \includegraphics[width=\textwidth]{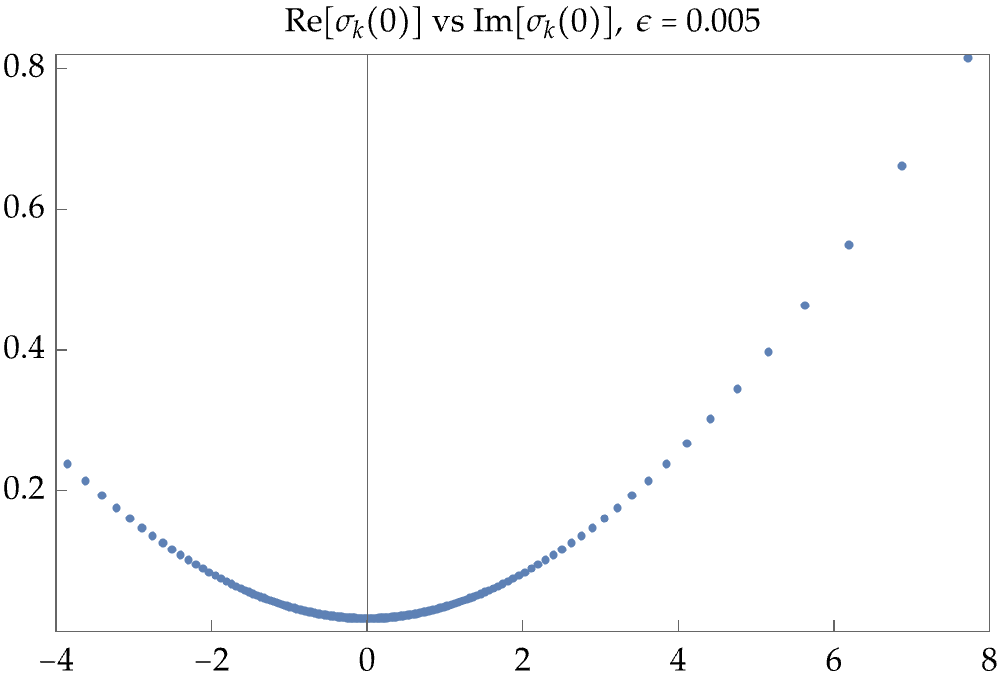}
     \end{subfigure}
     %\hfill
     \begin{subfigure}[b]{0.32\textwidth}
         \centering
         \includegraphics[width=\textwidth]{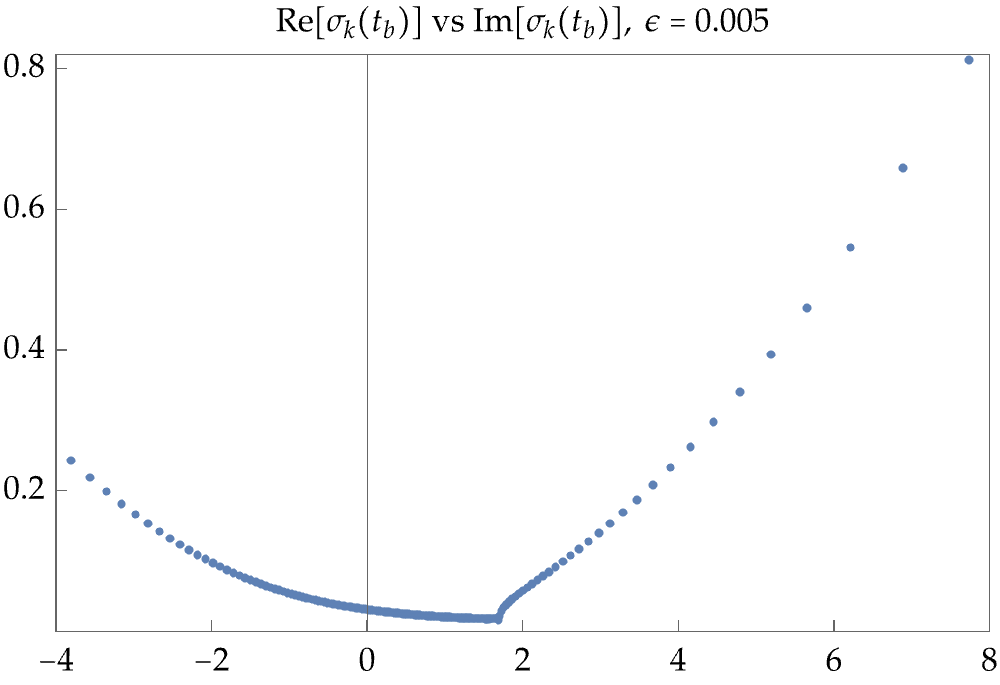}
     \end{subfigure}
     %\hfill
     \begin{subfigure}[b]{0.32\textwidth}
         \centering
         \includegraphics[width=\textwidth]{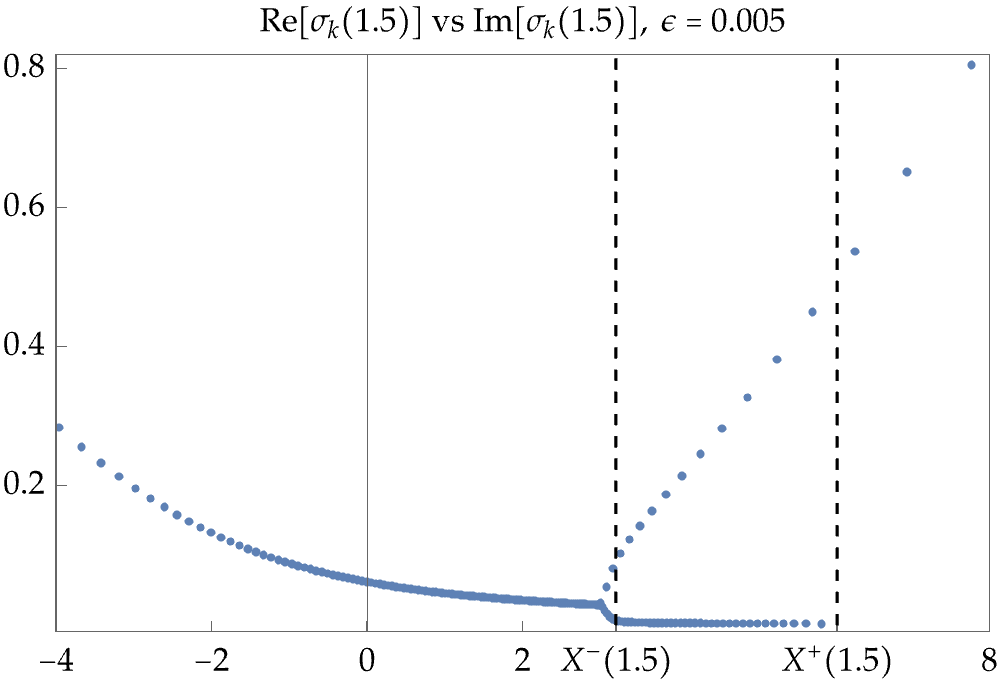}
     \end{subfigure}
        \caption{Time evolution of $\sigma_k(t)$ for $u_0(x)=2(1+x^2)^{-1}$ and $\epsilon=0.005$.}
        \label{fig: sigma time evo}
\end{figure}

What evidently distinguishes the outlier eigenvalues from the rest is that they have small real parts and imaginary parts that are possibly large instead of small.  We will make the following definition:
\begin{definition}[Outliers]
    For fixed $\beta>0$ and $B>0$ independent of $t$ but possibly depending on the initial condition $u_0(\cdot)$, an  eigenvalue $\sigma=\mu+\ii\nu$ will be called an outlier if  $\nu\ge \beta$ and $|\mu|\le B$.
    \label{def:outliers}
\end{definition}
For suitable initial conditions $u_0(\cdot)$, this definition is saying that the outliers include in particular any eigenvalues with bounded real part that saturate the upper bound given in \eqref{eq:rational-nu-bound}.  
Let $S_o\subset\{\sigma_k(t)\}$ denote the subset of outliers, and let the remaining eigenvalues constitute the ``bulk'' $S_b$, so that $S_b\sqcup S_0=\{\sigma_k(t)\}$.  The numerics suggest that there are relatively few outliers compared to the total number $N(\epsilon)\sim\epsilon^{-1}$ of eigenvalues.  Hence we formulate the following conjecture:
\begin{conjecture}[Outliers]\label{conjecture: outliers}
There is an exponent $0\le p<1$ and a constant $K>0$ independent of $t$ such that $|S_o|\le K\epsilon^{-p}$ holds for all $\epsilon>0$ sufficiently small.
\end{conjecture}
In Proposition~\ref{prop: outliers} below we will use this to estimate the contribution of the outliers to the sum \eqref{u in terms of sigmas}.  The bound $|\mu|\le B$ on the real part is part of Definition~\ref{def:outliers} because there are evidently some eigenvalues with both real and imaginary parts that are large when $\epsilon$ is small (see the upward trend with increasing real part of the points near the real line in the plots shown in Figures~\ref{fig:outliers x2} and \ref{fig: sigma time evo}), and as they apparently lie along certain curves and there can be many of them, it is more natural to include them in $S_b$.

Next, we examine the bulk  $S_b$ of the complex eigenvalues, which requires a different scale for the imaginary part coordinate.  As shown in Figure~\ref{fig: sigma time evo}, these bulk eigenvalues appear to be distributed along curves in the complex upper half-plane, and these curves have interesting dynamical behavior as $t\ge 0$ varies.
From these and similar plots it seems clear that when $t\le t_b$, the bulk eigenvalues arrange themselves along a single curve in the upper half-plane; however as soon as $t>t_b$, a second curve bifurcates from the former curve into the part of the upper half-plane below.  We therefore further partition $S_b$ into an ``upper branch'' subset $S_{U}\subseteq S_b$ and a complementary ``lower branch'' subset $S_{L}:=S_b\setminus S_{U}$ (when $t\le t_b$, $S_{L}=\emptyset$).  For $t>t_b$ fixed, both $S_{U}$, $S_{L}$ are non-empty and judging from numerics both have cardinality  $N_U(\epsilon)$, $N_L(\epsilon)$ proportional to $N(\epsilon)$.  Moreover, numerics suggest that the real parts of the eigenvalues belonging to $S_{L}$ lie within the interval $(X^{-}(t),X^{+}(t))$.  For example, when $t=1.5$, $X^{-}(t)\simeq3.20$ and $X^{+}(t)\simeq6.04$, see Figure \ref{fig: sigma time evo}, right panel.

From the plots, it appears that the real parts of the eigenvalues on either branch are locally equally spaced.  We formulate the following conjecture to express this observation in detail.

\begin{conjecture}[Real parts]\label{conjecture: mu}
For all $t\ge 0$ and $\epsilon>0$, the eigenvalues comprising $S_{U}$ have distinct real parts; similarly for $S_{L}$ when $t>t_b$.  Assume that the points of $S_{U}$ and $S_{L}$ are indexed by increasing real part, i.e., $\mu_{U,1}(t)<\mu_{U,2}(t)<\cdots<\mu_{U,N_U(\epsilon)}(t)$ for all $t\ge 0$ and $\mu_{L,1}(t)<\mu_{L,2}(t)<\cdots<\mu_{L,N_L(\epsilon)}(t)$ for all $t>t_b$.  Then, the sets $\{\mu_{U,k}(t)\}$, $\{\mu_{L,k}(t)\}$ are approximate samplings of two respective $\epsilon$-independent functions $\mu_U:(0,1)\times[0,\infty)\to\mathbb{R}$, $\mu_{L}:(0,1)\times(t_b,\infty)\to(X^{-}(t),X^{+}(t))$. More precisely, with $\prime$ denoting differentiation with respect to the first argument,
\begin{equation}
    \mu_{U,k}(t)=\mu_{U}(y_k,t)+\mathcal{O}(\epsilon^2\mu_U'(y_k,t)),\quad y_k:=\frac{k-\frac{1}{2}}{N_U(\epsilon)}
    \label{eq:mu-U-sampling}
\end{equation} 
and, for $t>t_b$, 
\begin{equation}
    \mu_{L,k}(t)=\mu_{L}(y_k,t)+\mathcal{O}(\epsilon^2\mu_L'(y_k,t)),\quad y_k:=\frac{k-\frac{1}{2}}{N_L(\epsilon)}
    \label{eq:mu-L-sampling}
\end{equation}
where the error terms are uniform for $k=1,\dots,N_U(\epsilon)$ and $k=1,\dots,N_L(\epsilon)$ respectively.
The functions $\mu_{U}$, $\mu_{L}$ have the following properties.
    \begin{itemize}
    \item[] \textbf{Invertibility:} $\mu_{L}(\cdot,t)$, $\mu_{U}(\cdot,t)$ are strictly increasing surjective functions for any $t$ on their respective domains, so there exists an inverse function $\mu^{-1}_{L,U}(\cdot,t)$, in the sense that $\mu^{-1}_{L,U}(\mu_{L,U}(y,t),t)=y$.
    
    \item[] \textbf{Smoothness and bounds for $\mu_U$:} For $0\le t<t_b$, $y\mapsto\mu_U(y,t)$ is of class $C^2((0,1))$. For $t\ge t_b$, $y\mapsto\mu_U(y,t)$ of class $C^2((0,1)\setminus\{y^-(t)\})$ where $y^{-}(t):=\mu_{U}^{-1}(X^{-}(t),t)$, and is Lipschitz continuous on $(0,1)$. In both cases there are exponents $q_\pm>0$ and constants $0<c<C$ such that 
    \begin{equation*}
        -Cy^{-q_-}<\mu_U(y,t)<-cy^{-q_-}\quad\text{and}\quad  \mu_U'(y,t)=\mathcal{O}(y^{-(q_-+1)}),\quad y\downarrow 0
    \end{equation*} 
    and 
    \begin{equation*}
        c(1-y)^{-q_+}<\mu_U(y,t)<C(1-y)^{-q_+}\quad\text{and}\quad\mu_U'(y,t)=\mathcal{O}((1-y)^{-(q_++1)}),\quad y\uparrow 1.
    \end{equation*}

    \item[] \textbf{Smoothness and bounds for $\mu_L$:} For $t>t_b$, $y\mapsto\mu_{L}(y,t)$ extends by continuity to $y\in[0,1]$ and is of class $C^2((0,1))$.

    \end{itemize}
    All estimates involving $\mu_U$ are uniform for bounded $t\ge 0$, in which case $\epsilon N_U(\epsilon)$ has a finite nonzero limit as $\epsilon\to 0$.  Likewise all estimates involving $\mu_L$ are uniform for bounded $t>t_b$ with $t-t_b$ bounded below by a positive quantity, in which case also $\epsilon N_L(\epsilon)$ has a finite nonzero limiting value.
\end{conjecture}

Note that $\mu_{U,L}(y_{k+1},t)-\mu_{U,L}(y_k,t)\approx \mu_{U,L}'(y_k,t)/N_{U,L}(\epsilon)$.  Since $N_{U,L}(\epsilon)$ are inversely proportional to $\epsilon$, equations \eqref{eq:mu-U-sampling} and \eqref{eq:mu-L-sampling} assert that the sampling error is $\epsilon$ times a uniform multiple of the local spacing.

In particular, the power-law behavior of $\mu_U(y,t)$ asserted in this conjecture is strongly supported by numerical experiments.  For the initial condition \eqref{eq:Lorentzian}, we estimated the exponents $q_\pm$ from slopes of best-fit lines in plots of $\ln(|\mu_{U,1}(t)|)$ and $\ln(\mu_{U,N_U(\epsilon)}(t))$ versus $\ln(\epsilon)$ over the range from $2^{-6}$ through $2^{-12}$ and obtained
\[
\begin{tabular}{|c||c|c|}
\hline
$t$ & $q_-$ & $q_+$ \\
\hline\hline
$0$ & $1.00021$ & $1.00021$\\
$t_b$ & $1.00023$ & $1.00019$\\
$1.5$ & $1.00028$ & $1.00013$\\
\hline
\end{tabular}
\]
Therefore, for this initial condition, the exponents $q_\pm$ appear to be approximately equal to $1$ regardless of whether $t<t_b$, $t=t_b$, or $t>t_b$.  We also found that the linear fit of the data is extremely accurate over the full range of scales.

\begin{figure}
     \centering
     \begin{subfigure}[b]{0.34\textwidth}
         \centering
         \includegraphics[width=\textwidth]{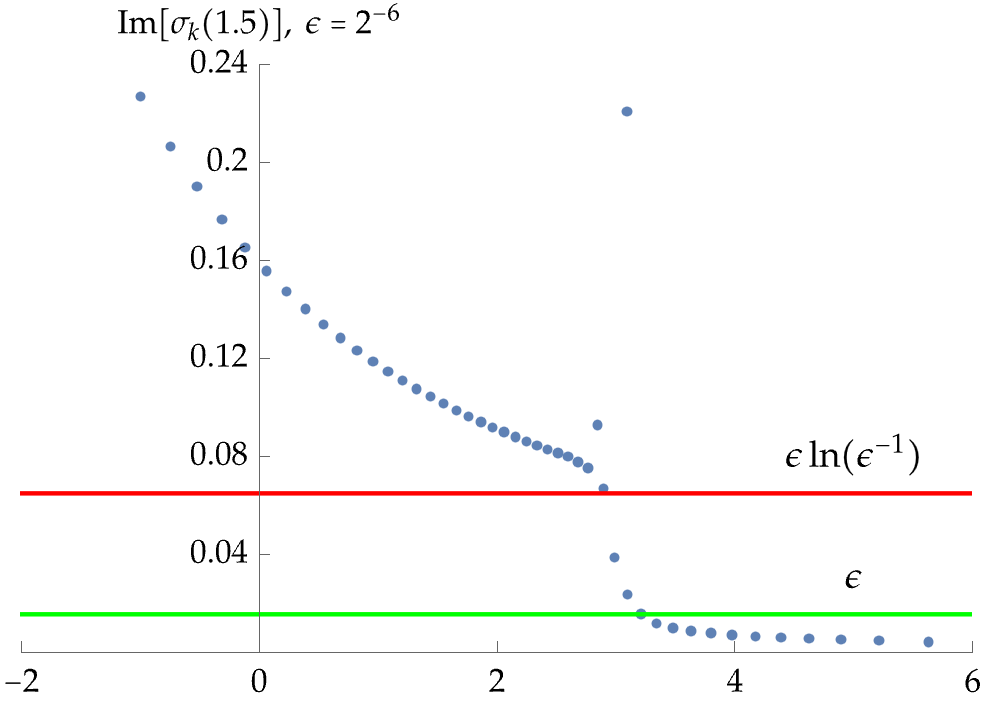}
     \end{subfigure}
     %\hfill
     \hspace{-0.4cm}
     \begin{subfigure}[b]{0.34\textwidth}
         \centering
         \includegraphics[width=\textwidth]{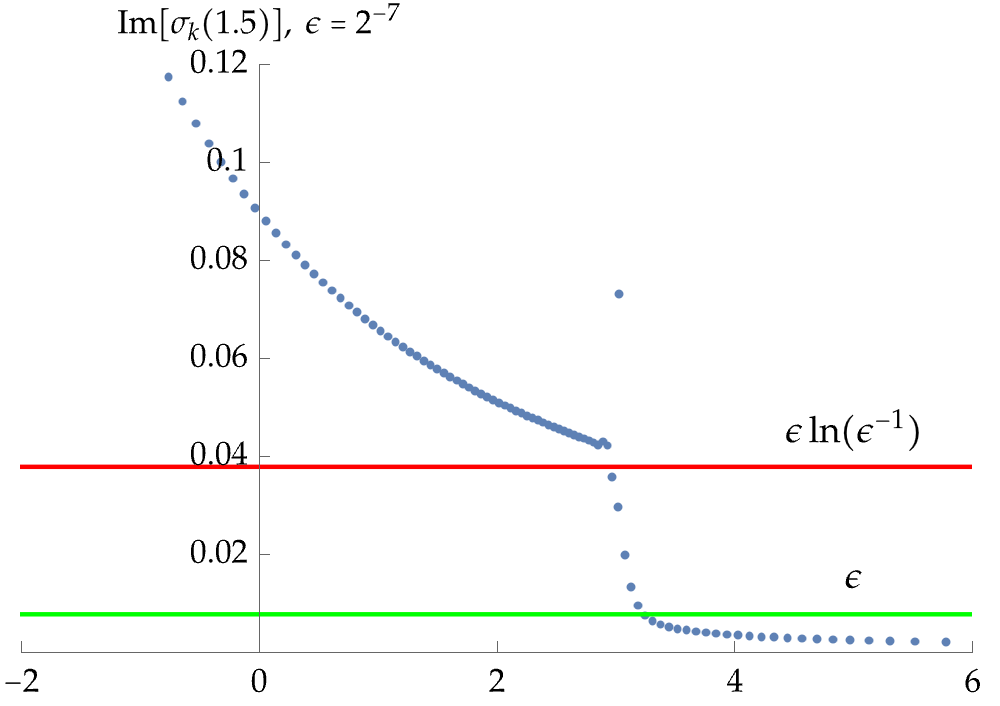}
     \end{subfigure}
     %\hfill
     \hspace{-0.4cm}
     \begin{subfigure}[b]{0.34\textwidth}
         \centering
         \includegraphics[width=\textwidth]{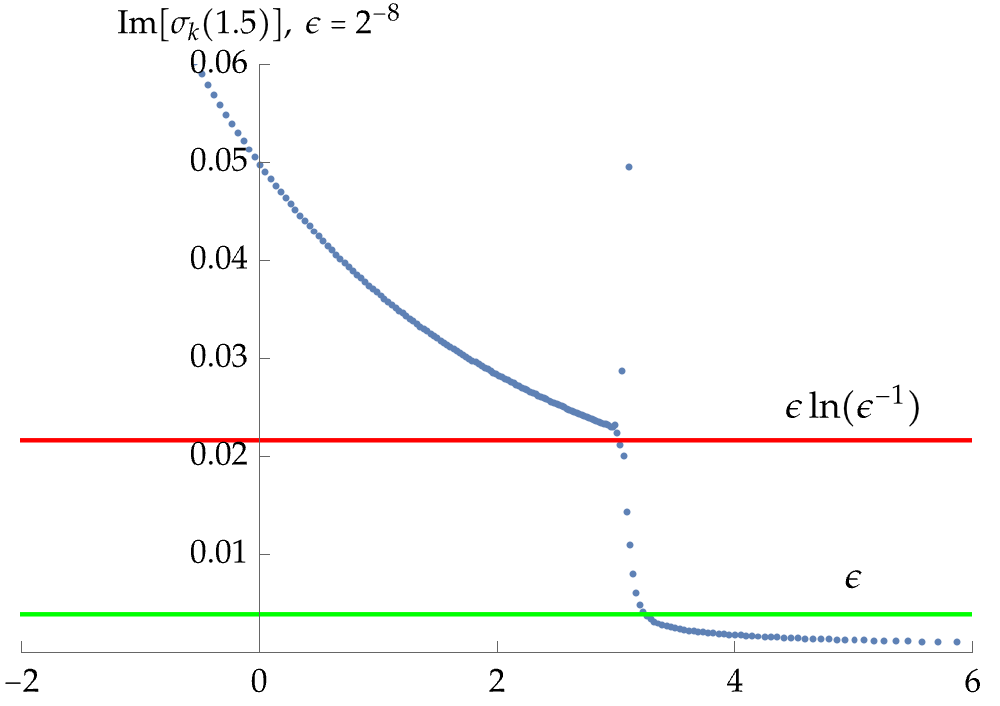}
     \end{subfigure}
        \caption{Upper and lower branch scaling for $\epsilon=2^{-6},2^{-7},2^{-8}$ and $t=1.5$.  The horizontal axis is $\mathrm{Re}[\sigma_k]$ in each plot.}
        \label{fig: branch scaling}
\end{figure}

In order to formulate similar conjectures regarding the imaginary parts of the eigenvalues on the upper and lower branches, it is necessary first to notice that unlike the real parts, the imaginary parts of eigenvalues on both branches are asymptotically small as $\epsilon\to 0$, with different asymptotic scales.  See Figure~\ref{fig: branch scaling}. These plots strongly suggest that the imaginary parts of the eigenvalues in $S_L$ scale proportionally with $\epsilon$, while those on the upper branch are somewhat larger.  In fact, the plots in Figure~\ref{fig: branch scaling} provide good evidence that the scaling of imaginary parts of eigenvalues in $S_U$ is proportional to $\epsilon\ln(\epsilon^{-1})$.  Observing $\mathrm{log}$-type growth/decay via numerics is notoriously difficult.  However, we are fortunate in that the conditional results we will prove in Section~\ref{sec:conditional-asymptotics-sigmas} below are rather insensitive to the precise scale of the imaginary parts of eigenvalues on the upper branch.  It will be enough that they scale as $\delta(\epsilon)$ lying in the asymptotic range $\epsilon\ll\delta(\epsilon)\ll 1$.  

\begin{conjecture}[Upper branch imaginary parts]\label{conjecture: nu upper}
    There is a scale $\delta(\epsilon)$ in the asymptotic range $\epsilon\ll\delta(\epsilon)\ll 1$ such that
    the set $\{\nu_{U,k}(t)\}$ is the approximate sampling of a scaled $\epsilon$-independent function $\nu_U:(0,1)\times[0,\infty)\to(0,\infty)$, i.e., as $\epsilon\to 0$,
    \begin{equation}
        \nu_{U,k}(t)=\delta(\epsilon)\nu_U(y_k,t)(1+\mathrm{o}(1)),\quad y_k:=\frac{k-\frac{1}{2}}{N_{U}(\epsilon)}
        \label{eq:nu-U-sampling}
    \end{equation}  
    holds uniformly over all indices $k=1,\dots,N_U(\epsilon)$.
    The function $\nu_U$ has the following properties.
    \begin{itemize}
        \item[]\textbf{Lower bound:} there exists a constant $c_U>0$ such that $c_U<\nu_{U}(y,t)$.
        
        \item[]\textbf{Edge behavior:}  there are exponents $0<r_\pm<2q_\pm+1$ (see Conjecture~\ref{conjecture: mu}) and constants $0<c<C$ such that 
        \[
        cy^{-r_-}(1-y)^{-r_+}<\nu_U(y,t)<Cy^{-r_-}(1-y)^{-r_+},\quad y\to 0,1.
\]
        \item[] \textbf{Smoothness:} For $0\le t<t_b$, $y\mapsto\nu_{U}(y,t)$ is 
        of class $C^1((0,1))$.  For $t\ge t_b$, $y\mapsto\nu_U(y,t)$ is absolutely continuous on $(0,1)$ and is of class $C^1((0,1)\setminus\{y^-(t)\})$, where $y^-(t)$ is as defined in Conjecture~\ref{conjecture: mu}.
     
    \end{itemize}
    All estimates are uniform for bounded $t\ge 0$.
\end{conjecture}

Again, the power-law behavior of $\nu_U(y,t)$ as $y\to 0,1$ is supported by numerics.  For the same initial condition \eqref{eq:Lorentzian}, by determining the scaling of $\nu_{U,k}(t)/(\epsilon\ln(\epsilon^{-1}))$ for $k=1$ and $k=N_U(\epsilon)$ with $\epsilon$, we obtained the following best-fit values of $r_\pm$, again with very accurate approximation over the whole range of scales:
\[
\begin{tabular}{|c||c|c|}
\hline
$t$ & $r_-$ & $r_+$ \\
\hline\hline
$0$ & $1.83460$ & $1.83460$\\
$t_b$ & $1.83456$ & $1.83464$\\
$1.5$ & $1.83445$ & $1.83475$\\
\hline
\end{tabular}
\]
Comparing with the experimental values of $q_\pm$, we see that indeed the inequalities $0<r_\pm<2q_\pm+1$ are evidently satisfied.

Finally, we offer a conjecture summarizing our numerical observations concerning the imaginary parts of eigenvalues on the lower branch if $t>t_b$.
\begin{conjecture}[Lower branch imaginary parts]\label{conjecture: nu lower}
    The set $\{\nu_{L,k}(t)\}$ is the approximate sampling of a scaled $\epsilon$-independent function $\nu_L:(0,1)\times(t_b,\infty)\to(0,\infty)$, i.e., as $\epsilon\to 0$,
    \begin{equation}
        \nu_{L,k}(t)=\epsilon\nu_{L}(y_k,t)(1+\mathrm{o}(1)),\quad y_k:=\frac{k-\frac{1}{2}}{N_{L}(\epsilon)}
        \label{eq:nu-L-sampling}
    \end{equation}  
    holds uniformly over all indices $k=1,\dots,N_L(\epsilon)$.
    The function $\nu_L$ has the following properties, in which $\delta(\epsilon)$ is the scale from Conjecture~\ref{conjecture: nu upper}, here asserted to have the additional property that
     \begin{equation}
        \int_0^\epsilon\frac{\delta(y)}{y}\,\dd y<\infty
        \label{eq:integrability}
    \end{equation}
    holds for $\epsilon>0$ sufficiently small.  
    \begin{itemize}
        \item[]\textbf{Lower bound:} there exists a constant $c_L>0$ such that $c_L<\nu_L(y,t)$.
        \item[]\textbf{Left edge behavior:}  
    There exist constants $c, C>0$ such that 
        \[
        c\delta(y)y^{-1}<\nu_{L}(y,t)<C\delta(y)y^{-1},\quad y\downarrow 0.
        \]
    
        \item[]\textbf{Smoothness:}  For $t>t_b$, $y\mapsto\nu_{L}(y,t)$ is of class $C^1((0,1))$.
    \end{itemize}
    All estimates are uniform for bounded $t>t_b$ with $t-t_b$ bounded below by a positive quantity.
\end{conjecture}

\begin{remark}
    The statement that $\nu_L(y,t)\sim \delta(y)y^{-1}$ as $y\downarrow 0$ is intended to capture the phenomenon that the lower branch bifurcates from the upper branch at $y=0$ (according to the parametrization of the lower branch; it is instead the point $y^-(t)$ in that of the upper branch) and the imaginary parts of the eigenvalues on the upper branch are asymptotically large compared to those of the lower branch eigenvalues.   The assertion of integrability of $\delta(y)y^{-1}$ at $y=0$ is reasonable given the numerical evidence that $\delta(\epsilon)\sim \epsilon\ln(\epsilon^{-1})$; see Figure~\ref{fig: branch scaling}.
\end{remark}

\begin{figure}[h]
\begin{center}
    \includegraphics[width=0.45\linewidth]{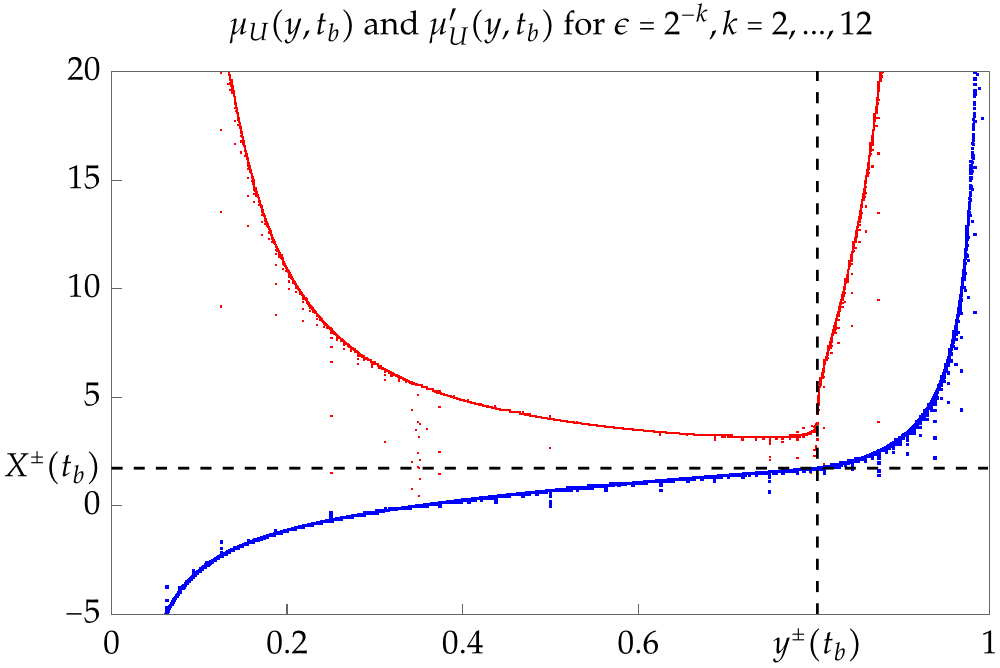}\hfill%
    \includegraphics[width=0.45\linewidth]{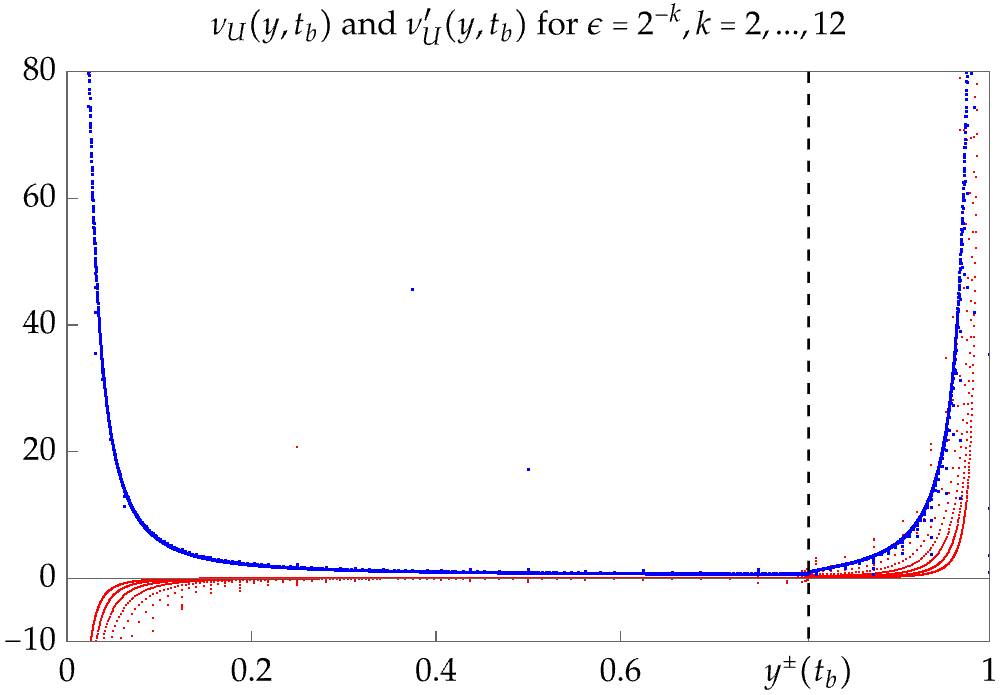}\hfill\\
    \includegraphics[width=0.45\linewidth]{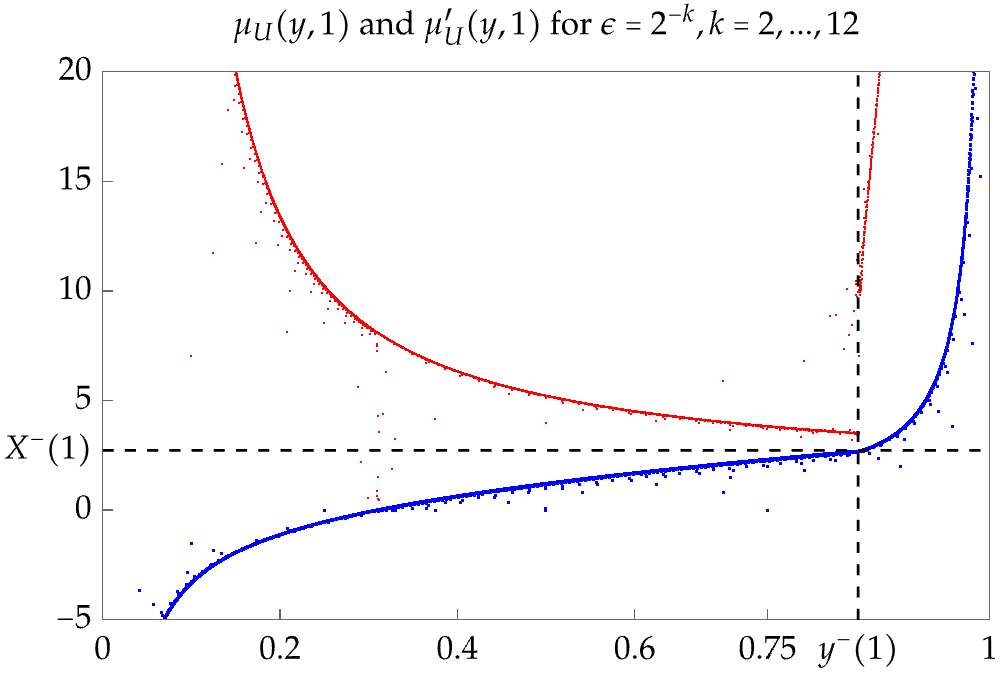}\hfill%
    \includegraphics[width=0.45\linewidth]{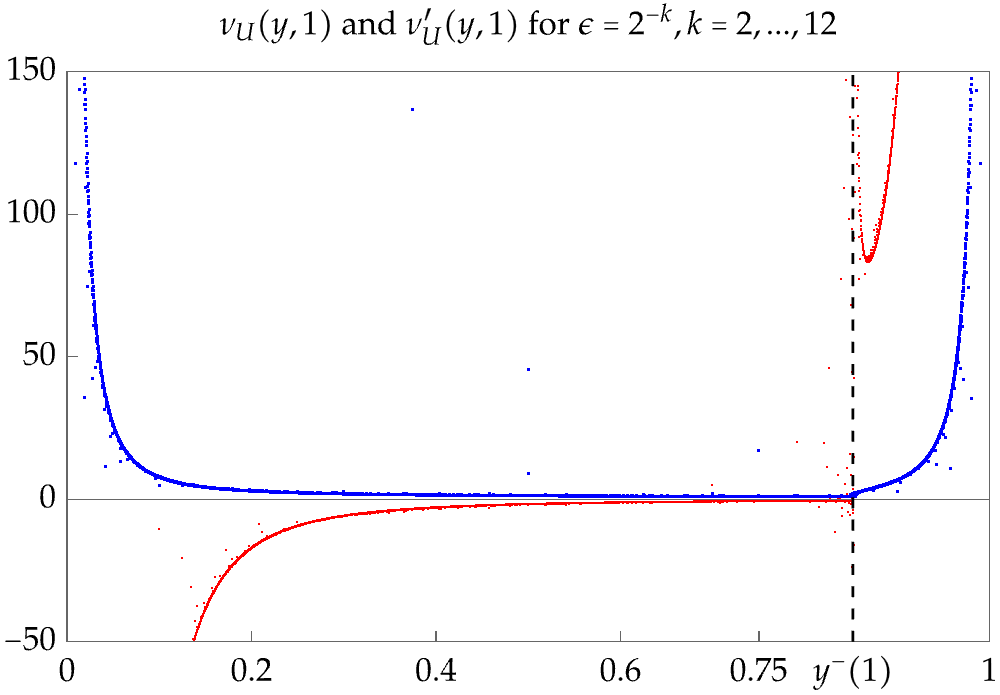}\hfill\\
    \includegraphics[width=0.45\linewidth]{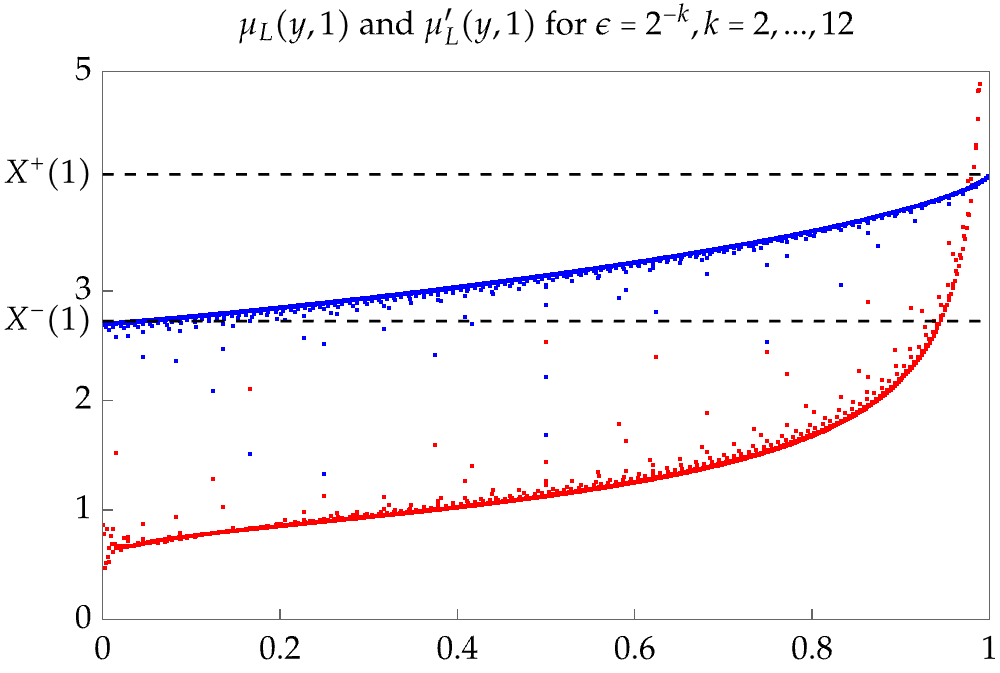}\hfill%
    \includegraphics[width=0.45\linewidth]{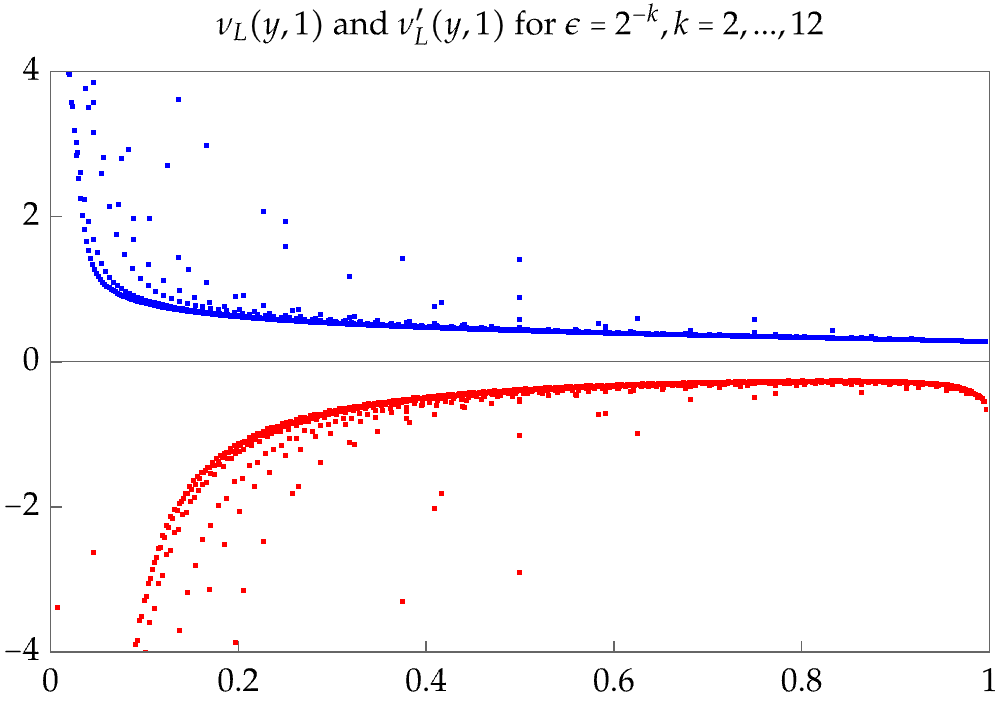}\hfill
\end{center}
\caption{The real and rescaled imaginary parts of the eigenvalues (blue) in the sets $\{\sigma_{U,k}(t)\}=S_{U}$ and $\{\sigma_{L,k}(t)\}=S_{L}$,  obtained from the matrix $\mathbf{C}(t)$ for the initial condition \eqref{eq:Lorentzian}, plotted versus the parameter $y\in (0,1)$ for $t=t_b$ (upper two panels) and $t=1>t_b$ (lower four panels).  Also shown are the numerical derivatives with respect to $y$ (red).}
\label{fig: mu nu prime}
\end{figure}

The approximate sampling properties asserted in Conjectures~\ref{conjecture: mu}, \ref{conjecture: nu upper} and \ref{conjecture: nu lower} are easily illustrated.  After computing the eigenvalues of $\mathbf{C}(t)$ for a given small value of $\epsilon$, one omits the outliers and partitions the bulk into $S_b=S_U\sqcup S_L$ based on the size of the imaginary parts, and these sets of complex numbers are then ordered by increasing real parts.  Setting $N_U(\epsilon):=|S_U|$, each point $\sigma_{U,k}(t)=\mu_{U,k}(t)+\ii\nu_{U,k}(t)$ is then assigned a value of the parameter $y=y_k:=(k-\frac{1}{2})/N_U(\epsilon)$.  Likewise, if $t>t_b$, setting $N_L(\epsilon):=|S_L|$, each point $\sigma_{L,k}(t)=\mu_{L,k}(t)+\ii\nu_{L,k}(t)$ is assigned a value of the parameter $y=y_k:=(k-\frac{1}{2})/N_L(\epsilon)$.  Assuming the hypothetical scale $\delta(\epsilon)=\epsilon\ln(\epsilon^{-1})$, the points $\{(y_k,\mu_{U,k}(t))\}$ and  $\{(y_k,\nu_{U,k}/\delta(\epsilon))\}$, and if $t>t_b$, $\{(y_k,\mu_{L,k}(t))\}$ and $\{(y_k,\nu_{L,k}(t)/\epsilon)\}$ can be plotted on the same axes for a range of scales for $\epsilon$; these are shown in blue in the panels of Figure~\ref{fig: mu nu prime}.
One can clearly observe the blue points condensing onto limiting fixed curves as $\epsilon\to 0$.  These limiting curves should be the graphs of the sampling functions $y\mapsto \mu_U(y,t)$, $y\mapsto \nu_U(y,t)$, and for $t>t_b$, $y\mapsto \mu_L(y,t)$ and $y\mapsto \nu_L(y,t)$. It is also straightforward to compute from the data difference quotient approximations (numerical derivatives) as follows:

\begin{align}
\mu'_{U}\bigg(\frac{k-\frac{1}{2}}{N_U(\epsilon)},t\bigg)&\simeq\frac{\mu_{U,k+1}(t)-\mu_{U,k}(t)}{N_U(\epsilon)^{-1}}, &&\nu_{U}'\bigg(\frac{k-\frac{1}{2}}{N_U(\epsilon)},t\bigg)\simeq\frac{\nu_{U,k+1}(t)-\nu_{U,k}(t)}{N_U(\epsilon)^{-1}\epsilon\ln(\epsilon^{-1})}, \\
\mu'_{L}\bigg(\frac{k-\frac{1}{2}}{N_L(\epsilon)},t\bigg)&\simeq\frac{\mu_{L,k+1}(t)-\mu_{L,k}(t)}{N_{L}(\epsilon)^{-1}}, &&\nu_{L}'\bigg(\frac{k-\frac{1}{2}}{N_L(\epsilon)},t\bigg)\simeq\frac{\nu_{L,k+1}(t)-\nu_{L,k}(t)}{N_L(\epsilon)^{-1}\epsilon}.
\end{align}
These are plotted against $y$ as red points in the same figure.  The red points also condense onto limiting curves, although it is clear that the limiting derivative curves for the data obtained from the upper-branch eigenvalues are discontinuous at $y=y^-(t)$ whenever $t>t_b$.

We have not yet been able to prove these conjectures.  We are aware of formal techniques applicable to nonselfadjoint eigenvalue problems for differential equations with analytic coefficients, and such methods have been applied to deduce curves in the complex plane that attract eigenvalues in a semiclassical limit similar to $N(\epsilon)\to\infty$  (see, e.g., \cite{Miller01}).  Unfortunately, the nonhermitian matrix $\mathbf{C}(t)$ does not have fixed bandwidth, so making an analogy with differential or difference equations is challenging. Even though $\mathbf{C}(t)$ does have an approximate Toeplitz structure, we are not aware of methods in the theory of Toeplitz quantization for nonhermitian matrices that would be sufficiently powerful to prove the above conjectures.  We will however apply elements of Toeplitz quantization to the Hermitian matrix $\mathbf{A}(x,t)$ in Section~\ref{sec:A} below.

\subsection{Conditional small-\texorpdfstring{$\epsilon$}{epsilon} asymptotics of the sum \texorpdfstring{\eqref{u in terms of sigmas}}{}}
\label{sec:conditional-asymptotics-sigmas}
\begin{figure}
     \centering
     \begin{subfigure}[b]{0.32\textwidth}
         \centering
         \includegraphics[width=\textwidth]{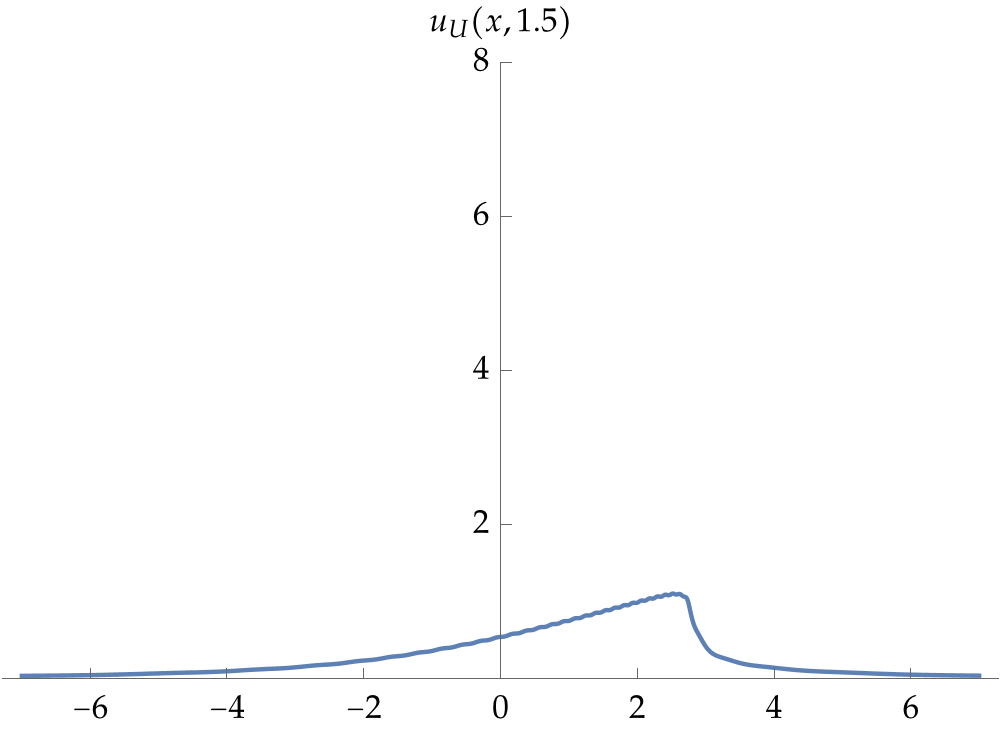}
     \end{subfigure}
     %\hfill
     \hspace{-0.3cm}
     \begin{subfigure}[b]{0.32\textwidth}
         \centering
         \includegraphics[width=\textwidth]{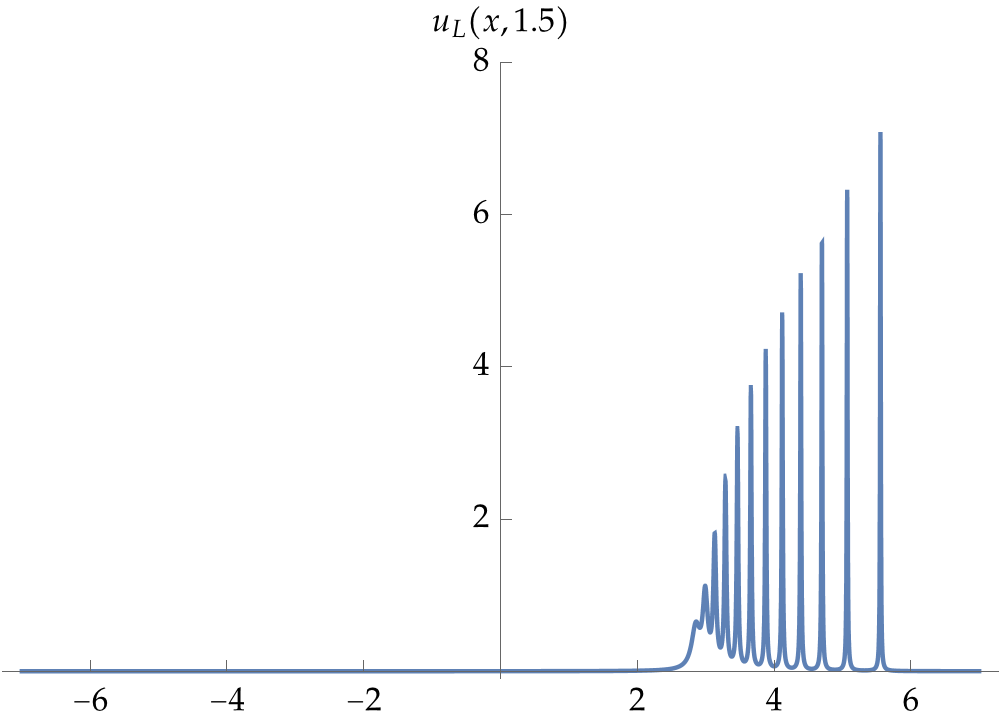}
     \end{subfigure}
     %\hfill
     \hspace{-0.3cm}
     \begin{subfigure}[b]{0.32\textwidth}
         \centering
         \includegraphics[width=\textwidth]{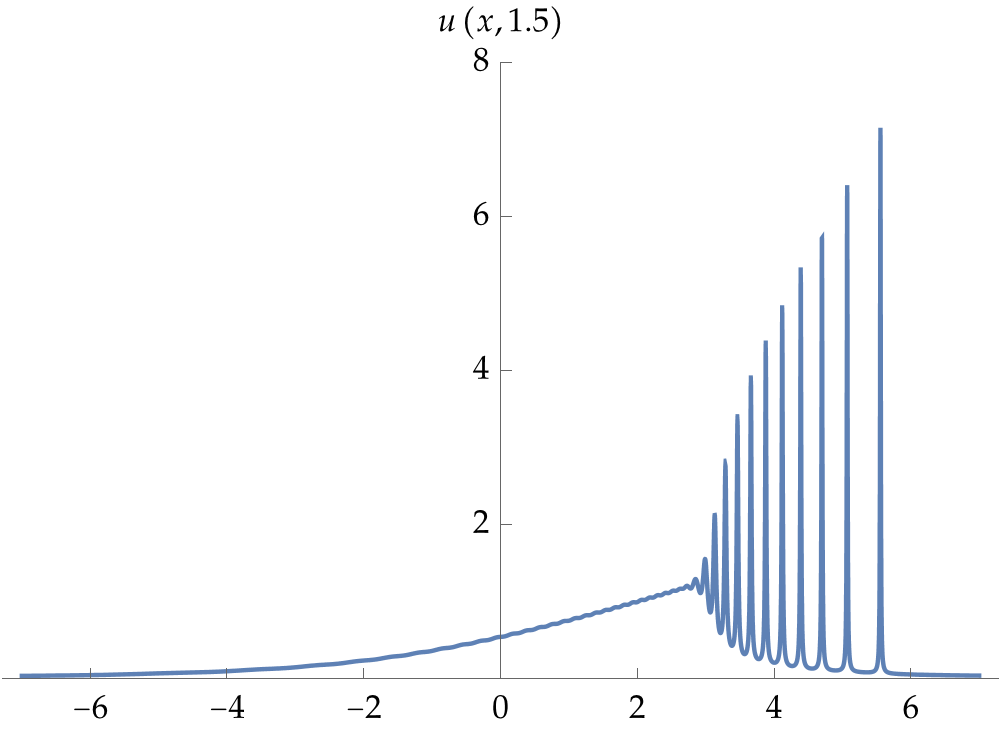}
     \end{subfigure}
        \caption{The two essential components $u_U(x,t)$ and $u_L(x,t)$ (see \eqref{def: lower upper sum}) of $u(x,t)$ compared with $u(x,t)$ itself for $u_0(x)=2(1+x^2)^{-1}$, $t=1.5$, and $\epsilon=0.02$.}
        \label{fig:split sum x2}
\end{figure}

Now we will use the conjectures inspired by numerical experiments to study the soliton ensemble for the BO equation with initial data $u_0(x)$.  
We can split the sum \eqref{u in terms of sigmas} into three parts corresponding to outliers in $S_o$, upper branch eigenvalues in $S_{U}$, and lower branch eigenvalues in $S_{L}$, i.e.,
\begin{align}\label{eq: u split sum}
    u(x,t)=u_{U}(x,t)+u_{L}(x,t)+u_{o}(x,t),
\end{align}
where
\begin{equation}\label{def: lower upper sum}
\begin{split}
   u_{U}(x,t)&:=\sum_{\sigma(t)\in S_{U}}\frac{2\epsilon\nu(t)}{(x-\mu(t))^2+\nu(t)^2}=\sum_{k=1}^{N_U(\epsilon)}\frac{2\epsilon\nu_{U,k}(t)}{(x-\mu_{U,k}(t))^2+\nu_{U,k}(t)^2}\\
    u_{L}(x,t)&:=\sum_{\sigma(t)\in S_{L}}\frac{2\epsilon\nu(t)}{(x-\mu(t))^2+\nu(t)^2}=\sum_{k=1}^{N_L(\epsilon)}\frac{2\epsilon\nu_{L,k}(t)}{(x-\mu_{L,k}(t))^2+\nu_{L,k}(t)^2}
\end{split}
\end{equation}
and $u_{o}(x,t):=u(x,t)-u_U(x,t)-u_L(x,t)$ is the sum over the remaining outlier eigenvalues.  It will turn out that $u_o(x,t)$ is of lower order compared to the other components and can be ignored (see Proposition \ref{prop: outliers} below).  The other two terms evidently make quite different contributions to the overall sum, as can be seen in Figure \ref{fig:split sum x2}, which compares $u_U(x,t)$ and $u_L(x,t)$ with $u(x,t)$ for a value of $t$ that exceeds the breaking time $t_b$.  
These plots suggest that the rapid oscillations (which occur only after the breaking time) come from the term $u_{L}(x,t)$ while the smooth background comes instead from $u_{U}(x,t)$.

First, we formulate a lemma showing that replacing the eigenvalues with their respective sampling functions in either of the two summands \eqref{def: lower upper sum} produces a relatively small error term that is uniform with respect to the sum index.

\begin{lemma}[Sampling error]
Assume that Conjecture~\ref{conjecture: mu} holds and fix constants $K>0$ and $\eta>0$.  Then in the limit $\epsilon\to 0$, the summand in $u_U(x,t)$ can be written in the form
\begin{equation}
    \frac{2\epsilon\nu_{U,k}(t)}{(x-\mu_{U,k}(t))^2+\nu_{U,k}(t)^2}=\frac{2\epsilon\delta(\epsilon)\nu_U(y_k,t)}{(x-\mu_U(y_k,t))^2+\delta(\epsilon)^2\nu_U(y_k,t)^2}(1+o(1)) ,\quad y_k:=\frac{k-\frac{1}{2}}{N_U(\epsilon)}
    \label{eq:U-summand-approx}
\end{equation}
provided also Conjecture~\ref{conjecture: nu upper} holds, where the error term is uniform for $0\le t\le K$, $|x|\le K$, and $k=1,\dots,N_U(\epsilon)$.  Likewise, as $\epsilon\to 0$, 
the summand in $u_L(x,t)$ can be written in the form
\begin{equation}
        \frac{2\epsilon\nu_{L,k}(t)}{(x-\mu_{L,k}(t))^2+\nu_{L,k}(t)^2}=\frac{2\epsilon^2\nu_L(y_k,t)}{(x-\mu_L(y_k,t))^2+\epsilon^2\nu_L(y_k,t)^2}(1+o(1)) ,\quad y_k:=\frac{k-\frac{1}{2}}{N_L(\epsilon)}
        \label{eq:L-summand-approx}
\end{equation}
provided also Conjecture~\ref{conjecture: nu lower} holds,
where the error term is uniform for $t_b+\eta\le t\le K$, $X^-(t)+\eta\le x\le X^+(t)-\eta$ and $k=1,\dots,N_L(\epsilon)$.
\label{lemma: sampling summand}
\end{lemma}

\begin{proof}
For convenience we suppress the dependence on $t$, which should be taken in different intervals depending on which summand is considered.
According to Conjectures~\ref{conjecture: nu upper} and \ref{conjecture: nu lower}, we have $\nu_{U,k}=\delta(\epsilon)\nu_U(y_k)(1+o(1))$ with $y_k=(k-\frac{1}{2})/N_U(\epsilon)$ uniformly for $k=1,\dots,N_U(\epsilon)$, and 
$\nu_{L,k}=\epsilon\nu_L(y_k)(1+o(1))$ with $y_k=(k-\frac{1}{2})/N_L(\epsilon)$ uniformly for $k=1,\dots,N_L(\epsilon)$. Therefore, it remains to get a corresponding uniform $o(1)$ relative error estimate for the denominator of each summand.  

In other words, we want to obtain uniform $o(1)$ estimates for
\begin{equation}
Q_{U,k}:=\left|\frac{[(x-\mu_{U,k})^2+\nu_{U,k}^2]-[(x-\mu_U(y_k))^2+\delta(\epsilon)^2\nu_U(y_k)^2]}{(x-\mu_U(y_k))^2+\delta(\epsilon)^2\nu_U(y_k)^2}\right|,\quad y_k:=\frac{k-\frac{1}{2}}{N_U(\epsilon)}    
\end{equation}
and
\begin{equation}
Q_{L,k}:=\left|\frac{[(x-\mu_{L,k})^2+\nu_{L,k}^2]-[(x-\mu_L(y_k))^2+\epsilon^2\nu_L(y_k)^2]}{(x-\mu_L(y_k))^2+\epsilon^2\nu_L(y_k)^2}\right|,\quad y_k:=\frac{k-\frac{1}{2}}{N_L(\epsilon)}. 
\end{equation}
Let us write $\mu_{U,k}=\mu_U(y_k)+\Delta\mu_{U,k}$, $\mu_{L,k}=\mu_L(y_k)+\Delta\mu_{L,k}$, $\nu_{U,k}=\delta(\epsilon)\nu_U(y_k)(1+\Delta\nu_{U,k})$, and $\nu_{L,k}=\epsilon\nu_L(y_k)(1+\Delta\nu_{L,k})$.  We can assume that $|\Delta\nu_{U,k}|< 1$ and $|\Delta\nu_{L,k}|<1$, so
\begin{equation}
    Q_{U,k}\le \frac{2|x-\mu_U(y_k)||\Delta\mu_{U,k}|+|\Delta\mu_{U,k}|^2 +3\delta(\epsilon)^2\nu_U(y_k)^2|\Delta\nu_{U,k}|}{(x-\mu_U(y_k))^2+\delta(\epsilon)^2\nu_U(y_k)^2}\le R_{U,k}  + 3|\Delta\nu_{U,k}|
\end{equation}
where
\begin{equation}
    R_{U,k}:=\frac{2|x-\mu_U(y_k)||\Delta\mu_{U,k}|+|\Delta\mu_{U,k}|^2 }{(x-\mu_U(y_k))^2+\delta(\epsilon)^2\nu_U(y_k)^2},\quad y_k:=\frac{k-\frac{1}{2}}{N_U(\epsilon)}
\end{equation}
and
\begin{equation}
        Q_{L,k}\le \frac{2|x-\mu_L(y_k)||\Delta\mu_{L,k}|+|\Delta\mu_{L,k}|^2 +3\epsilon^2\nu_L(y_k)^2|\Delta\nu_{L,k}|}{(x-\mu_L(y_k))^2+\epsilon^2\nu_L(y_k)^2}\le R_{L,k} + 3|\Delta\nu_{L,k}|
\end{equation}
where
\begin{equation}
    R_{L,k}:=\frac{2|x-\mu_L(y_k)||\Delta\mu_{L,k}|+|\Delta\mu_{L,k}|^2 }{(x-\mu_L(y_k))^2+\epsilon^2\nu_L(y_k)^2},\quad y_k:=\frac{k-\frac{1}{2}}{N_L(\epsilon)}.
\end{equation}
Since $\Delta\nu_{U,k}=o(1)$ and $\Delta\nu_{L,k}=o(1)$ both hold uniformly for $k=1,\dots,N_U(\epsilon)$ and $k=1,\dots,N_L(\epsilon)$ respectively, it is enough to estimate $R_{U,k}$ and $R_{L,k}$.

First we consider $R_{U,k}$.  If $k$ is such that $|x-\mu_U(y_k)|\ge 1$, then by neglecting $\delta(\epsilon)^2\nu_U(y_k)^2$ in the denominator,
\begin{equation}
    R_{U,k}\le 2\frac{|\Delta\mu_{U,k}|}{|x-\mu_U(y_k)|} + \left(\frac{|\Delta\mu_{U,k}|}{|x-\mu_U(y_k)|}\right)^2.
\end{equation}
But then we have the estimate by Conjecture~\ref{conjecture: mu} 
\begin{equation}
\frac{|\Delta\mu_{U,k}|}{|x-\mu_U(y_k)|}=\mathcal{O}\left(\frac{\epsilon^2\mu_U'(y_k)}{|x-\mu_U(y_k)|}\right)
\end{equation}
which is uniform for $k$ in the range $1\le k\le N_U(\epsilon)$.  This is $\mathcal{O}(\epsilon^2)$ as long as $y_k$ is bounded away from the endpoints $0,1$.  However, according to Conjecture~\ref{conjecture: mu}, $\mu'_U(y)/\mu_U(y)=\mathcal{O}(y^{-1})$ as $y\downarrow 0$ while $\mu'_U(y)/\mu_U(y)=\mathcal{O}((1-y)^{-1})$ as $y\uparrow 1$.  Using these estimates and the fact that $y_k$ and $1-y_k$ are both greater than or equal to $\frac{1}{2}N_U(\epsilon)^{-1}$ shows that 
\begin{equation}
    |x-\mu_U(y_k)|\ge 1\implies \frac{|\Delta\mu_{U,k}|}{|x-\mu_U(y_k)|}=\mathcal{O}(\epsilon^2N_U(\epsilon))=\mathcal{O}(\epsilon) \implies R_{U,k}=\mathcal{O}(\epsilon)
\end{equation}
holds uniformly over the indicated subset of indices $k$.  On the other hand, if $|x-\mu_U(y_k)|<1$, then by neglecting $(x-\mu_U(y_k))^2$ in the denominator
\begin{equation}
    R_{U,k}<\frac{2|\Delta\mu_{U,k}| + |\Delta\mu_{U,k}|^2}{\delta(\epsilon)^2\nu_U(y_k)^2}\le \frac{2|\Delta\mu_{U,k}| + |\Delta\mu_{U,k}|^2}{\delta(\epsilon)^2c_U^2},
\end{equation}
where $\nu_U(y)\ge c_U$ holds for $0<y<1$ by Conjecture~\ref{conjecture: nu upper}.  Since $\mu'(y)$ is uniformly bounded if $\mu(y)$ itself is, $|x-\mu_U(y_k)|<1$ controls $\mu'(y_k)$ uniformly for bounded $x$, and hence $\Delta\mu_{U,k}=\mathcal{O}(\epsilon^2\mu'(y_k))=\mathcal{O}(\epsilon^2)$.  Therefore
\begin{equation}
    |x-\mu_U(y_k)|<1\implies R_{U,k}=\mathcal{O}\left(\frac{\epsilon^2}{\delta(\epsilon)^2}\right)=o(1)
\end{equation}
holds uniformly for bounded $x$ and in the indicated range of indices $k$ because $\epsilon\ll\delta(\epsilon)\ll 1$ as asserted in Conjecture~\ref{conjecture: nu upper}.  This shows that $R_{U,k}=o(1)$ holds uniformly for $k=1,\dots,N_U(\epsilon)$ if $x$ is uniformly bounded.  Hence also $Q_{U,k}=o(1)$ in the same sense.

Next, we consider $R_{L,k}$.  If $k$ is such that $|x-\mu_L(y_k)|\ge \epsilon$, then neglecting $\epsilon^2\nu_L(y_k)^2$ from the denominator as before,
\begin{equation}
 |x-\mu_L(y_k)|\ge\epsilon\implies   R_{L,k}\le 2\frac{|\Delta\mu_{L,k}|}{|x-\mu_L(y_k)|}+\left(\frac{|\Delta\mu_{L,k}|}{|x-\mu_L(y_k)|}\right)^2\le 2\frac{|\Delta\mu_{U,k}|}{\epsilon} + \left(\frac{|\Delta\mu_{U,k}|}{\epsilon}\right)^2.
\end{equation}
Since $\mu_L'(y)$ is uniformly bounded on $0<y<1$, we use the uniform estimate $\Delta\mu_{L,k}=\mathcal{O}(\epsilon^2\mu'_L(y_k))=\mathcal{O}(\epsilon^2)$ from Conjecture~\ref{conjecture: mu} to obtain $R_{L,k}=\mathcal{O}(\epsilon)$ as a uniform bound for indices satisfying the indicated condition.  Conversely, if $|x-\mu_L(y_k)|<\epsilon$, we omit $(x-\mu_L(y_k))^2$ from the denominator instead and obtain
\begin{equation}
    |x-\mu_L(y_k)|<\epsilon\implies R_{L,k}< \frac{2\epsilon|\Delta\mu_{L,k}| + |\Delta\mu_{L,k}|^2}{\epsilon^2\nu_L(y_k)^2}\le \frac{2\epsilon|\Delta\mu_{L,k}| + |\Delta\mu_{L,k}|^2}{\epsilon^2c_L^2},
\end{equation}
using $\nu_L(y)\ge c_L$ for $0<y<1$ as follows from Conjecture~\ref{conjecture: nu lower}.  Again using $\Delta\mu_{L,k}=\mathcal{O}(\epsilon^2)$ and combining with the result for $|x-\mu_L(y_k)|\ge \epsilon$ shows that the bound $R_{L,k}=\mathcal{O}(\epsilon)$ holds uniformly for all indices $k=1,\dots,N_L(\epsilon)$. Hence also $Q_{L,k}=o(1)$ holds in the same sense.
\end{proof}

We now investigate the small-$\epsilon$ asymptotics of each sum, starting with the outlier sum $u_{o}(x,t)$.

\begin{proposition}[Outlier sum]\label{prop: outliers}
Let $\beta>0$ be as in Definition~\ref{def:outliers},  and assume that Conjecture~\ref{conjecture: outliers} holds with some constant $K>0$ and exponent $0\le p<1$.  Then, $1-p>0$, and
\begin{align}\label{eq: outlier sum}
    u_{o}(x,t)=\sum_{\sigma_k(t)\in S_o}\frac{2\epsilon\nu_k(t)}{(x-\mu_k(t))^2+\nu_k(t)^2}=\mathcal{O}(\epsilon^{1-p})
\end{align}
for sufficiently small $\epsilon>0$ where the error term is uniform with respect to $x, t$.
\end{proposition}

\begin{proof}
The sum \eqref{eq: outlier sum} is positive because all its terms are positive, so
\begin{align}
    \sum_{\sigma_k(t)\in S_o}\frac{2\epsilon\nu_k(t)}{(x-\mu_k(t))^2+\nu_k(t)^2}\leq\sum_{\sigma_k(t)\in S_o}\frac{2\epsilon}{\nu_k(t)}\leq \tilde{K}\epsilon^{1-p},
\end{align}
where $\tilde{K}=2K/\beta>0$ is independent of $x,t,\epsilon$, as desired.
\end{proof}

Next we consider the upper sum $u_U(x,t)$.  The main idea here is that the sum resembles a Riemann sum for an integral, which in turn has an integrand involving a highly-peaked function that approximates a Dirac delta.  The first observation we make is that the conjectured asymptotic behavior of $\mu_U(y,t)$ and $\nu_U(y,t)$ near the endpoints $y=0,1$ is enough to neglect the contributions of the extreme eigenvalues on the upper branch.

\begin{lemma}[Real parts close to $x$]\label{lemma:omit-extremes}
Assume that Conjectures~\ref{conjecture: mu} and \ref{conjecture: nu upper} both hold.  Then, uniformly for bounded $t\ge 0$ and bounded $x\in\mathbb{R}$, and arbitrary fixed $\eta>0$,
\begin{equation}
u_U(x,t)=\sum_{\substack{\sigma_k(t)\in S_U\\|\mu_k(t)-x|<\eta}}\frac{2\epsilon\nu_k(t)}{(x-\mu_k(t))^2+\nu_k(t)^2}+\mathcal{O}(\delta(\epsilon)),\quad\epsilon\to 0.
\end{equation}
\end{lemma}
\begin{proof}
    We estimate the complementary part of $u_U(x,t)$ by first invoking Lemma~\ref{lemma: sampling summand} to obtain (denoting $y_k:=(k-\frac{1}{2})/N_U(\epsilon)$)
    \begin{equation}
        \label{eq:first-line-bound} \sum_{\substack{\sigma_k(t)\in S_U\\|\mu_k(t)-x|\ge \eta}}\frac{2\epsilon\nu_k(t)}{(x-\mu_k(t))^2+\nu_k(t)^2}=(1+o(1))\sum_{\substack{\sigma_k(t)\in S_U\\|\mu_k(t)-x|\ge \eta}}\frac{2\epsilon\delta(\epsilon)\nu_U(y_k,t)}{(x-\mu_U(y_k,t))^2+\delta(\epsilon)^2\nu_U(y_k,t)^2}
        \end{equation}
        where we used the fact that the terms are all positive. But neglecting $\delta(\epsilon)^2\nu_U(y_k,t)^2$ from the denominator and 
    using Conjectures~\ref{conjecture: mu} and \ref{conjecture: nu upper} we get 
    \begin{equation}
    \begin{split}
\sum_{\substack{\sigma_k(t)\in S_U\\\mu_k(t)\le x-\eta}}\frac{2\epsilon\delta(\epsilon)\nu_U(y_k,t)}{(x-\mu_U(y_k,t))^2+\delta(\epsilon)^2\nu_U(y_k,t)^2}&\le
        \sum_{\substack{\sigma_k(t)\in S_U\\\mu_k(t)\le x-\eta}}\frac{2\epsilon\delta(\epsilon)\nu_U(y_k,t)}{(x-\mu_U(y_k,t))^2}\\ &\lesssim\sum_{\substack{\sigma_k(t)\in S_U\\\mu_k(t)\le x-\eta}}\frac{2\epsilon\delta(\epsilon)y_k^{-r_-}}{y_k^{-2q_-}}.
    \end{split}
    \end{equation}
Letting $N_-(\epsilon)=\mathcal{O}(\epsilon^{-1})$ denote the greatest index $k$ for which $\mu_k(t)\le x-\eta$, this latter sum is just
\begin{equation}
\begin{split}
    \sum_{k=1}^{N_-(\epsilon)}\frac{2\epsilon \delta(\epsilon)y_k^{-r_-}}{y_k^{-2q_-}} &= 
    2\epsilon\delta(\epsilon)N_U(\epsilon)^{-(2q_--r_-)}\sum_{k=1}^{N_-(\epsilon)}(k-\tfrac{1}{2})^{-(2q_--r_-)}\\ &\lesssim 2\epsilon\delta(\epsilon)N_U(\epsilon)^{-(2q_--r_-)}N_-(\epsilon)^{1-(2q_--r_-)},
\end{split}
\end{equation}
because $2q_--r_->-1$, which is $\mathcal{O}(\delta(\epsilon))$ as desired.  The estimate of the part of the upper bound in \eqref{eq:first-line-bound} with $\mu_k(t)\ge x+\eta$ is similar, using $2q_+-r_+>-1$.
\end{proof}

Next, we have a Lemma that will help us compute the remaining terms in the upper sum $u_{U}(x,t)$.

\begin{lemma}\label{lemma: delta}
Suppose that $x\in\mathbb{R}$ and $\eta>0$.  Let $f:[x-\eta,x+\eta]\to\mathbb{R}$ be a continuous function with positive lower bound $f(z)\ge c_1>0$ that is differentiable with uniformly bounded derivative: $|f'(z)|\le c_2$. The function
\begin{equation}
    D_\delta(z;x):=\frac{1}{\pi}\cdot\frac{\delta f(z)}{(z-x)^2+\delta^2f(z)^2}
\end{equation}
with parameter $\delta>0$ is an approximate delta function in the sense that for any bounded Lipschitz continuous function $g:[x-\eta,x+\eta]\to\mathbb{R}$,
\begin{equation}
    \int_{x-\eta}^{x+\eta}D_\delta(z;x)g(z)\,\dd z = g(x)+\mathcal{O}(\delta\ln(\delta^{-1})),\quad\delta\to 0.
\end{equation}
\end{lemma}

\begin{proof}
Setting $s(z):=(z-x)/(\delta f(z))$, we have
\begin{equation}
\begin{split}
    \int_{x-\eta}^{x+\eta}D_\delta(z;x)\,\dd z &= \int_{x-\eta}^{x+\eta}\frac{1}{s(z)^2+1}\left(s'(z)+\frac{(z-x)f'(z)}{\delta f(z)^2}\right)\frac{\dd z}{\pi}\\
    &=\int_{-\eta/(\delta f(x-\eta))}^{\eta/(\delta f(x+\eta))}\frac{1}{s^2+1}\frac{\dd s}{\pi} + \int_{x-\eta}^{x+\eta}\frac{\delta (z-x)f'(z)}{(z-x)^2+\delta^2f(z)^2}\frac{\dd z}{\pi}\\
    &=1+\mathcal{O}(\delta) + \int_{x-\eta}^{x+\eta}\frac{\delta (z-x)f'(z)}{(z-x)^2+\delta^2f(z)^2}\frac{\dd z}{\pi},\quad\delta\to 0,
\end{split}
\end{equation}
where on the last line we used the positive lower bound for $f(z)$.  But using also the upper bound for $|f'(z)|$,
\begin{equation}
\begin{split}
    \left|\int_{x-\eta}^{x+\eta}\frac{\delta (z-x)f'(z)}{(z-x)^2+\delta^2f(z)^2}\frac{\dd z}{\pi}\right|\le \frac{c_2\delta}{\pi}\int_{x-\eta}^{x+\eta}\frac{|z-x|\,\dd z}{(z-x)^2+\delta^2c_1^2}&=\frac{c_2\delta}{\pi}\int_{-(\delta c_1)^{-1}}^{(\delta c_1)^{-1}}\frac{|w|\,\dd w}{w^2+1} \\ &= \mathcal{O}(\delta\ln(\delta^{-1})),\quad\delta\to 0.
\end{split}
\end{equation}
Therefore,
\begin{equation}
    \int_{x-\eta}^{x+\eta}D_\delta(z;x)g(z)\,\dd z = g(x)+\mathcal{O}(\delta\ln(\delta^{-1})) + \int_{x-\eta}^{x+\eta}D_\delta(z;x)(g(z)-g(x))\,\dd z.
\end{equation}
Let $c_g>0$ denote the Lipschitz constant of $g$, so that $|g(z)-g(x)|\le c_g|z-x|$. Then
\begin{equation}
\begin{split}
    \left|\int_{x-\eta}^{x+\eta}D_\delta(z;x)(g(z)-g(x))\,\dd z\right|&\le c_g\int_{x-\eta}^{x+\eta}D_\delta(z;x)|z-x|\,\dd z\\ &\le
    \frac{Cc_g\delta}{\pi}\int_{x-\eta}^{x+\eta}\frac{|z-x|\,\dd z}{(z-x)^2+\delta^2c_1^2}.
\end{split}
\end{equation}
where $C$ denotes the finite positive maximum value of $f(z)$ on $[x-\eta,x+\eta]$.  This upper bound is again $\mathcal{O}(\delta\ln(\delta^{-1}))$, so the proof is finished.
\end{proof}

Now we may compute the upper branch sum.
\begin{proposition}[Upper branch sum]\label{prop: upper branch leading order}
Assume that Conjectures~\ref{conjecture: mu} and \ref{conjecture: nu upper} hold.  Then,
\begin{align}\label{eq: upper bulk sum}
    u_{U}(x,t):=\sum_{\sigma_k(t)\in S_{U}}\frac{2\epsilon\nu_k(t)}{(x-\mu_k(t))^2+\nu_k(t)^2}=\psi_U(x,t)+o(1), \quad\epsilon\to 0,
\end{align}
where 
\begin{align}
    \psi_U(x,t):=\frac{2\pi\epsilon N_{U}(\epsilon)}{\mu'_{U}(\mu_{U}^{-1}(x,t),t)},\quad \mu_U'(y,t):=\partial_y\mu_U(y,t),
\label{eq:psiU}
\end{align}
and the $o(1)$ error term is uniform on bounded subsets of $(x,t)\in \mathbb{R}\times [0,\infty)$ for which $(x,t)$ is bounded away from the curve $(X^-(\tau),\tau)$ parametrized by $\tau\ge t_b$.
\end{proposition}
\begin{remark}
    The product  $\epsilon N_U(\epsilon)$ converges as $\epsilon\to 0$ to a finite nonzero value depending on $t$ only.  In fact, for $0\le t\le t_b$, $\epsilon N_U(\epsilon)\to M$ defined in \eqref{eq:LM} as $\epsilon\to 0$.
\end{remark}

\begin{proof}
Thanks to Lemma~\ref{lemma:omit-extremes} and the fact that $\delta(\epsilon)=o(1)$ as $\epsilon\to 0$ it suffices to study the sum over eigenvalues $\sigma_k(t)\in S_U$ for which $|\mu_k(t)-x|<\eta$ for any suitable $\eta>0$. A suitable value for $\eta$ will be specified at the end of the proof.  Letting $k=k_-,\dots,k_+$ denote the index range for which $|\mu_k(t)-x|<\eta$ holds, applying Lemma~\ref{lemma: sampling summand} and using the positivity of the summand, the contributing terms can be written as
\begin{equation}    \sum_{\substack{\sigma_k(t)\in S_U\\|\mu_k(t)-x|<\eta}}\frac{2\epsilon\nu_k(t)}{(x-\mu_k(t))^2+\nu_k(t)^2} =(1+o(1))\sum_{k=k_-}^{k_+}\frac{2\epsilon\delta(\epsilon)\nu(y_k)}{(x-\mu(y_k))^2+\delta(\epsilon)^2\nu(y_k)^2},
\label{eq:sampling-error}
\end{equation}
where $y_k:=(k-\frac{1}{2})/N_U(\epsilon)$, and for simplicity we are abbreviating $\mu(y):=\mu_U(y,t)$ and $\nu(y):=\nu_U(y,t)$.
Since we have a sum of a sampling of a continuous function over a large number of grid points because $k_+-k_-\gtrsim\epsilon^{-1}$, we can apply the Euler-Maclaurin formula in the form
\begin{equation}
    \sum_{k=k_-}^{k_+}h(k)=\int_{k_-}^{k_+}h(k)\,\dd k + \frac{h(k_+)+h(k_-)}{2}+\int_{k_-}^{k_+}h'(k)(k-\lfloor k\rfloor-\tfrac{1}{2})\,\dd k
    \label{eq:EL}
\end{equation}
with $h(k)$ a function defined on the real interval $[k_-,k_+]$ by
\begin{equation}    h(k):=\frac{2\epsilon\delta(\epsilon)\nu(y_k)}{(x-\mu(y_k))^2+\delta(\epsilon)^2\nu(y_k)^2},\quad y_k:=\frac{k-\frac{1}{2}}{N_U(\epsilon)},\quad k\in [k_+,k_-].
\label{eq:hfunc-def}
\end{equation}
Now, when $k=k_\pm$, $\mu(y_k)\approx x\pm \eta$, so certainly $|x-\mu(y_{k_\pm})|\ge\frac{1}{2}\eta$.  Therefore $|h(k_\pm)|\le8\epsilon\delta(\epsilon)\nu(y_{k_\pm})\eta^{-1}$.  Since $\nu(y)$ is a continuous function of $y$ on any closed subinterval of $(0,1)$, and is independent of $\epsilon$, we easily obtain 
\begin{equation}
    \frac{h(k_+)+h(k_-)}{2}=\mathcal{O}(\epsilon\delta(\epsilon)),\quad\epsilon\to 0.
\end{equation}
The last integral on the right-hand side of \eqref{eq:EL} can be estimated using that $|k-\lfloor k\rfloor-\frac{1}{2}|\le\frac{1}{2}$ holds for all $k\in\mathbb{R}$, and hence with a change of the integration variable we get
\begin{equation}
    \left|\int_{k_-}^{k_+}h'(k)(k-\lfloor k\rfloor-\tfrac{1}{2})\,\dd k\right|\le\frac{1}{2}\int_{k_-}^{k_+}|h'(k)|\,\dd k = 
    \frac{1}{2}\int_{y_{k_-}}^{y_{k_+}}\left|\frac{\dd}{\dd y}h(N_U(\epsilon)y+\tfrac{1}{2})\right|\,\dd y.
\end{equation} 
Thus, using \eqref{eq:hfunc-def}
we have
\begin{equation}
     \left|\int_{k_-}^{k_+}h'(k)(k-\lfloor k\rfloor-\tfrac{1}{2})\,\dd k\right|\le\epsilon\delta(\epsilon)\int_{y_{k_-}}^{y_{k_+}}\left|\frac{\dd}{\dd y}\frac{\nu(y)}{(x-\mu(y))^2+\delta(\epsilon)^2\nu(y)^2}\right|\,\dd y.
\end{equation}
Doing the differentiation, we get a larger upper bound as a sum of three terms:
\begin{equation}
    \int_{y_{k_-}}^{y_{k_+}}\left|\frac{\dd}{\dd y}\frac{\nu(y)}{(x-\mu(y))^2+\delta(\epsilon)^2\nu(y)^2}\right|\,\dd y\le I_1+I_2+I_3,
\end{equation}
where
\begin{equation}
    \begin{split}
        I_1&:=\int_{y_{k_-}}^{y_{k_+}}\frac{|\nu'(y)|\,\dd y}{(x-\mu(y))^2+\delta(\epsilon)^2\nu(y)^2},\\
        I_2&:=\int_{y_{k_-}}^{y_{k_+}}\frac{2\delta(\epsilon)^2\nu(y)^2|\nu'(y)|\,\dd y}{[(x-\mu(y))^2+\delta(\epsilon)^2\nu(y)^2]^2},\\
        I_3&:=\int_{y_{k_-}}^{y_{k_+}}\frac{2|\mu(y)-x|\nu(y)\mu'(y)\,\dd y}{[(\mu(y)-x)^2+\delta(\epsilon)^2\nu(y)^2]^2},
    \end{split}
\end{equation}
where to write $I_3$ we used the fact that $\nu(y)>0$ and $\mu'(y)>0$.  Now, according to Conjecture~\ref{conjecture: nu upper}, we have a lower bound of the form $\nu(y)\ge c_U>0$, and as noted above $\nu(y)$ attains a finite maximum on $[y_{k_-},y_{k_+}]$ that we denote by $C$.  Likewise, Conjecture~\ref{conjecture: nu upper} asserts the absolute continuity of the derivative $\nu'(y)$ on $(0,1)$, and hence 
\begin{equation}
    I_1\le \frac{1}{c_U^2\delta(\epsilon)^2}\int_{y_{k_-}}^{y_{k_+}}|\nu'(y)|\,\dd y\quad\text{and}\quad
    I_2\le \frac{2C^2}{c_U^4\delta(\epsilon)^2}\int_{y_{k_-}}^{y_{k_+}}|\nu'(y)|\,\dd y.
\end{equation}
Since $y_{k_\pm}\approx \mu^{-1}(x\pm \eta)$, and the latter values do not depend on $\epsilon$, it is clear that $I_1+I_2=\mathcal{O}(\delta(\epsilon)^{-2})$ as $\epsilon \to 0$.  To estimate $I_3$, we first use the upper and lower bounds on $\nu(y)$ to obtain
\begin{equation}
    I_3\le C\int_{y_{k_-}}^{y_{k_+}}\frac{2|\mu(y)-x|\mu'(y)\,\dd y}{[(\mu(y)-x)^2+\delta(\epsilon)^2c_U^2]^2}.
    \label{eq:I3-bound}
\end{equation}
We next split the integral at the point $y=\mu^{-1}(x)$; then $\mu(y)>x$ holds for $\mu^{-1}(x)<y<y_{k_+}$ while $\mu(y)<x$ holds for $y_{k_-}<y<\mu^{-1}(x)$.  Therefore,
\begin{equation}
\begin{split}
    \int_{y_{k_-}}^{y_{k_+}}\frac{2|\mu(y)-x|\mu'(y)\,\dd y}{[(\mu(y)-x)^2+\delta(\epsilon)^2c_U^2]^2}&=\left(\int_{\mu^{-1}(x)}^{y_{k_+}}-\int_{y_{k_-}}^{\mu^{-1}(x)}\right)\frac{2(\mu(y)-x)\mu'(y)\,\dd y}{[(\mu(y)-x)^2+\delta(\epsilon)^2c_U^2]^2}\\    &=\left(\int_x^{\mu(y_{k_+})}-\int_{\mu(y_{k_-})}^x\right)\frac{2(\mu-x)\,\dd\mu}{[(\mu-x)^2+\delta(\epsilon)^2c_U^2]^2}.
    \end{split}
\end{equation}
Carrying out the explicit integration and using the result in \eqref{eq:I3-bound} gives
\begin{equation}
    I_3\le \frac{2C}{\delta(\epsilon)^2c_U^2}-\frac{C}{(\mu(y_{k_+})-x)^2+\delta(\epsilon)^2c_U^2}-\frac{C}{(\mu(y_{k_-})-x)^2+\delta(\epsilon)^2c_U^2}\le \frac{2C}{\delta(\epsilon)^2c_U^2}
\end{equation}
so also $I_3=\mathcal{O}(\delta(\epsilon)^{-2})$ as $\epsilon\to 0$.  Combining the estimates and using $\epsilon\ll\delta(\epsilon)\ll 1$ shows that
\begin{equation}
    \sum_{k=k_-}^{k=k_+}h(k)-\int_{k_-}^{k_+} h(k)\,\dd k = \mathcal{O}(\epsilon\delta(\epsilon)) + \mathcal{O}(\epsilon\delta(\epsilon)^{-1})
    = \mathcal{O}(\epsilon\delta(\epsilon)^{-1})=o(1),\quad\epsilon\to 0.
    \label{eq:sum-int-diff-estimate}
\end{equation}

It only remains to consider the integral
\begin{equation}
\begin{split}
    \int_{k_-}^{k_+}h(k)\,\dd k &= N_U(\epsilon)\int_{y_{k_-}}^{y_{k_+}}h(N_U(\epsilon)y+\tfrac{1}{2})\,\dd y\\ &=2\epsilon N_U(\epsilon)\int_{y_{k_-}}^{y_{k_+}}\frac{\delta(\epsilon)\nu(y)\,\dd y}{(x-\mu(y))^2+\delta(\epsilon)^2\nu(y)^2}\\
    &=2\pi \epsilon N_U(\epsilon)\cdot\frac{1}{\pi}
    \int_{\mu(y_{k_-})}^{\mu(y_{k_+})}\frac{\delta(\epsilon)\nu(\mu^{-1}(z))}{(x-z)^2+\delta(\epsilon)^2\nu(\mu^{-1}(z))^2}\frac{1}{\mu'(\mu^{-1}(z))}\,\dd z,
\end{split}
\end{equation}
where on the last line we changed variables by the increasing map $z=\mu(y)$.  Since $\mu(y_{k_\pm})=x\pm \eta+\mathcal{O}(\epsilon)$ as $\epsilon \to 0$ and the integrand is $\mathcal{O}(\delta(\epsilon))$ near the limits of integration, 
\begin{equation}
    \int_{k_-}^{k_+}h(k)\,\dd k = 2\pi\epsilon N_U(\epsilon)\cdot\frac{1}{\pi}\int_{x-\eta}^{x+\eta}\frac{\delta(\epsilon)\nu(\mu^{-1}(z))}{(x-z)^2+\delta(\epsilon)^2\nu(\mu^{-1}(z))^2}\frac{1}{\mu'(\mu^{-1}(z))}\,\dd z + \mathcal{O}(\epsilon\delta(\epsilon))
\end{equation}
because $\epsilon N_U(\epsilon)=\mathcal{O}(1)$.  Finally, we appeal to Lemma~\ref{lemma: delta} with $f(z):=\nu(\mu^{-1}(z))$ and $g(z):=\mu'(\mu^{-1}(z))^{-1}$, which will satisfy the required hypotheses provided we now choose $\eta>0$ appropriately.  If $0\le t<t_b$, we will take $\eta=1$, but if $t\ge t_b$, since $x\neq X^-(t)$, we will take $\eta<\frac{1}{2}|x-X^-(t)|$.  This choice guarantees that, according to Conjectures~\ref{conjecture: mu} and \ref{conjecture: nu upper}, $f(z)$ is continuous with a positive lower bound and is continuously differentiable (hence having bounded derivative) for $z\in [x-\eta,x+\eta]$; likewise $g(z)$ is continuously differentiable (hence Lipschitz) on the same interval.  We conclude that
\begin{equation}
\begin{split}
    \int_{k_-}^{k_+}h(k)\,\dd k&=\frac{2\pi\epsilon N_U(\epsilon)}{\mu'(\mu^{-1}(x))} + \mathcal{O}(\delta(\epsilon)\ln(\delta(\epsilon)^{-1}))\\
    &=\frac{2\pi\epsilon N_U(\epsilon)}{\mu'(\mu^{-1}(x))} + o(1),\quad\epsilon\to 0,
\end{split}
\end{equation}
because $\epsilon\ll \delta(\epsilon)\ll 1$.  Combining with \eqref{eq:sampling-error}, \eqref{eq:hfunc-def}, and \eqref{eq:sum-int-diff-estimate} completes the proof.
\end{proof}

Now we can turn our attention to the sum over the lower branch of eigenvalues, assuming that $(x,t)$ is inside the triple-valued region for Burgers' equation.

\begin{proposition}[Lower branch sum]\label{prop: lower branch leading order}
Fix $t>t_b$, and let $x_0$ be fixed in the interior of the triple-valued region for Burgers' equation, i.e., $X^-(t)<x_0<X^+(t)$.
Assume that Conjectures \ref{conjecture: mu} and \ref{conjecture: nu lower} hold.  Then,
\begin{align}\label{eq: lower bulk sum}
    u_{L}(x,t)=\frac{\psi_{L}(x_0,t)\mathrm{sinh}\left(\phi_{L}(x_0,t)\right)}{\mathrm{cosh}\left(\phi_{L}(x_0,t)\right)-\cos\left(\psi_{L}(x_0,t)\epsilon^{-1}(x-x_0)+2\pi p(x_0)\right)}+o(1),\quad\epsilon\to 0,
\end{align}
holds uniformly for $x-x_0=\mathcal{O}(\epsilon^{3/4})$
where $p(x_0)$ is a uniformly bounded (bound independent of $x_0$) phase shift, see \eqref{def:px0} below, and
\begin{align}
    \psi_{L}(x,t):=\frac{2\pi\epsilon N_{L}(\epsilon)}{\mu'_{L}(\mu_{L}^{-1}(x,t),t)}, ~~~ \phi_{L}(x,t):=\psi_{L}(x,t)\nu_{L}(\mu_{L}^{-1}(x,t),t).
\label{eq:psiL-phiL}
\end{align}
\end{proposition}

\begin{proof}
For brevity, we write $\nu(y):=\nu_L(y,t)$ and $\mu(y):=\mu_L(y,t)$ for the sampling functions defined in Conjectures~\ref{conjecture: mu} and \ref{conjecture: nu lower} as $t$ is fixed and we are only concerned with the lower branch.  For $u_L(x,t)$ defined by \eqref{def: lower upper sum}, we first apply Lemma~\ref{lemma: sampling summand} to get
\begin{align}\label{uL in proof}
    u_{L}(x,t)=(1+o(1))\sum_{k=1}^{N_L(\epsilon)}\frac{2\epsilon^2\nu(y_k)}{(x-\mu(y_k))^2+\epsilon^2\nu(y_k)^2}, \quad y_k:=\frac{k-\frac{1}{2}}{N_L(\epsilon)}
\end{align}
because the terms are all positive.  We will work with the explicit sum and then deal later with the multiplicative factor $1+o(1)$.  The leading-order contribution to the explicit sum will come from the terms where $\mu(y_k)$ is close to $x_0$ which we specify by the condition $|\mu(y_k)-x_0|\le 2\epsilon^r$ where $r$ is an exponent with $0<r<\frac{1}{2}$. Indeed, under the complementary condition $|\mu(y_k)-x_0|>2\epsilon^r$, we have $|x-\mu_k(y)|\ge ||\mu(y_k)-x_0|-|x-x_0||=|\mu(y_k)-x_0|-|x-x_0|\ge \epsilon^r$ for $\epsilon>0$ sufficiently small, because $x-x_0=\mathcal{O}(\epsilon^{3/4})=o(\epsilon^r)$ holds for $r<\frac{1}{2}$.  Therefore, neglecting the term $\epsilon^2\nu(y_k)^2$ in the denominator and extending the sum over the full range of $k$, the sum of complementary terms is estimated as follows:
\begin{equation}
    \sum_{\substack{k=1\\|\mu(y_k)-x_0|>2\epsilon^r}}^{N_L(\epsilon)}\frac{2\epsilon^2\nu(y_k)}{(x-\mu(y_k))^2+\epsilon^2\nu(y_k)^2}\le 2\epsilon N_L(\epsilon)\cdot\epsilon^{1-2r}\sum_{k=1}^{N_L(\epsilon)}\frac{\nu(y_k)}{N_L(\epsilon)}.
    \label{eq:complementary}
\end{equation}
The last sum is a Riemann sum for the integral $\int_0^1\nu(y)\,\dd y$ which is finite according to Conjecture~\ref{conjecture: nu lower}.  Since $\epsilon N_L(\epsilon)=\mathcal{O}(1)$ as $\epsilon\to 0$ and $r<\frac{1}{2}$, the sum \eqref{eq:complementary} of complementary terms tends to zero with $\epsilon$.

The remaining terms in the explicit sum on the right-hand side of \eqref{uL in proof} have indices $k$ in the set $\mathscr{I}^r:=\{j=1,2,\dots,N_L(\epsilon):|\mu(y_j)-x_0|\le 2\epsilon^r\}$.  Next we show that for the terms in the explicit sum with indices in $\mathscr{I}^r$, we can replace $\nu(y_k)$ in the summand with $\nu(\mu^{-1}(x_0))$ because $r>0$.  Indeed, consider the difference
\begin{multline}
    D:=\sum_{k\in \mathscr{I}^r}\left[\frac{2\epsilon^2\nu(y_k)}{(x-\mu(y_k))^2+\epsilon^2\nu(y_k)^2}-\frac{2\epsilon^2\nu(\mu^{-1}(x_0))}{(x-\mu(y_k))^2+\epsilon^2\nu(\mu^{-1}(x_0))^2}\right] \\
{}    = 2\epsilon^2\sum_{k\in \mathscr{I}^r}\frac{[(x-\mu(y_k))^2-\epsilon^2\nu(y_k)\nu(\mu^{-1}(x_0))](\nu(y_k)-\nu(\mu^{-1}(x_0)))}{[(x-\mu(y_k))^2+\epsilon^2\nu(y_k)^2][(x-\mu(y_k))^2+\epsilon^2\nu(\mu^{-1}(x_0))^2]}.
\label{eq:D-def}
\end{multline}
Since the function $x\mapsto \nu(\mu^{-1}(x))$ is continuously differentiable near $x=x_0$, the condition $|\mu(y_k)-x_0|\le 2\epsilon^r$ implies that $|\nu(y_k)-\nu(\mu^{-1}(x_0))|\lesssim \epsilon^r$, and therefore
\begin{equation}
    |D|\lesssim \epsilon^{2+r}\sum_{k\in \mathscr{I}^r}\frac{(x-\mu(y_k))^2 + \epsilon^2\nu(y_k)\nu(\mu^{-1}(x_0))}{[(x-\mu(y_k))^2+\epsilon^2\nu(y_k)^2][(x-\mu(y_k))^2+\epsilon^2\nu(\mu^{-1}(x_0))^2]}.
\end{equation}
Furthermore, according to Conjecture~\ref{conjecture: nu lower} we have the lower bound $\nu(y)\ge c_L>0$, and as $\nu$ is a $C^1$ function near $\mu^{-1}(x_0)>0$, we have also a local upper bound:  $\nu(y)\le K$, so 
\begin{equation}
    |D|\lesssim\epsilon^{2+r}\sum_{k\in\mathscr{I}^r}\frac{(x-\mu(y_k))^2}{[(x-\mu(y_k))^2+\epsilon^2c_L^2]^2} + \epsilon^{4+r}\sum_{k\in \mathscr{I}^r}\frac{1}{[(x-\mu(y_k))^2+\epsilon^2c_L^2]^2}.
\label{eq:absD2}
\end{equation}
Let $\mathscr{I}^r_0:=\{k\in \mathscr{I}^r: |x-\mu(y_k)|\le \epsilon\}$.  Then the cardinality of $\mathscr{I}^r_0$ is $|\mathscr{I}^r_0|=\mathcal{O}(1)$ as $\epsilon\downarrow 0$, so 
\begin{equation}
\begin{split}
    \sum_{k\in \mathscr{I}^r_0}\frac{(x-\mu(y_k))^2}{[(x-\mu(y_k))^2+\epsilon^2c_L^2]^2}&\le \sum_{k\in \mathscr{I}^r_0}\frac{\epsilon^2}{\epsilon^4c_L^4} \lesssim \epsilon^{-2},\\
    \sum_{k\in \mathscr{I}_0^r}\frac{1}{[(x-\mu(y_k))^2+\epsilon^2c_L^2]^2}&\le \sum_{k\in \mathscr{I}_0^r}\frac{1}{\epsilon^4c_L^4}\lesssim\epsilon^{-4},
\end{split}
\end{equation}
and therefore the terms indexed by $\mathscr{I}^r_0$ contribute $\mathcal{O}(\epsilon^r)$ to the right-hand side of the estimate \eqref{eq:absD2}.  Then for the remaining indices in $\mathscr{I}^r$ we can write
\begin{equation}
\begin{split}
    \sum_{k\in \mathscr{I}^r\setminus \mathscr{I}^r_0}\frac{(x-\mu(y_k))^2}{[(x-\mu(y_k))^2+\epsilon^2c_L^2]^2}&\le\sum_{k\in \mathscr{I}^r\setminus \mathscr{I}^r_0}\frac{1}{(x-\mu(y_k))^2},\\
    \sum_{k\in \mathscr{I}^r\setminus \mathscr{I}^r_0}\frac{1}{[(x-\mu(y_k))^2+\epsilon^2c_L^2]^2}&\le\sum_{k\in \mathscr{I}^r\setminus \mathscr{I}^r_0}\frac{1}{(x-\mu(y_k))^4}.
\end{split}
\label{eq:remaining-terms-in-Sr}
\end{equation}
Then, for either exponent $p=2,4$, we can use the positive lower bound on $\mu'(y)=\mu_L'(y,t)$ implied by Conjecture~\ref{conjecture: mu} to get
\begin{equation}
    \sum_{k\in \mathscr{I}^r\setminus \mathscr{I}^r_0}\frac{1}{(x-\mu(y_k))^p}\lesssim N_L(\epsilon)\sum_{k\in \mathscr{I}^r\setminus \mathscr{I}^r_0}\frac{\mu'(y_k)}{(x-\mu(y_k))^p}\cdot\frac{1}{N_L(\epsilon)}.
\end{equation}
Since $1/N_L(\epsilon)$ is exactly the spacing of the points $y_k$, the latter sum is a Riemann sum for an integral, and hence
\begin{equation}
    \sum_{k\in \mathscr{I}^r\setminus \mathscr{I}^r_0}\frac{1}{(x-\mu(y_k))^p}\lesssim N_L(\epsilon)\int_{|\mu-x|\ge\epsilon}\frac{\dd\mu}{(x-\mu)^p},
\end{equation}
wherein the integration is extended to $\mu=\pm\infty$ as a finite upper bound since $p=2,4$.  Performing the integration shows that
\begin{equation}
    \sum_{k\in \mathscr{I}^r\setminus \mathscr{I}^r_0}\frac{1}{(x-\mu(y_k))^p}\lesssim N_L(\epsilon)\epsilon^{1-p}\lesssim\epsilon^{-p},\quad p=2,4.
\end{equation}
Combining with \eqref{eq:remaining-terms-in-Sr} shows that also the terms with indices in $\mathscr{I}^r\setminus \mathscr{I}^r_0$ contribute $\mathcal{O}(\epsilon^r)$ to the right-hand side of the estimate \eqref{eq:absD2} so $D\to 0$ as $\epsilon\to0$.

Therefore, it remains to analyze the sum
\begin{equation}
    \Sigma:=\sum_{k\in \mathscr{I}^r}\frac{2\epsilon^2\nu(\mu^{-1}(x_0))}{(x-\mu(y_k))^2+\epsilon^2\nu(\mu^{-1}(x_0))^2},
\label{eq:Sigma-def}
\end{equation}
which can be written in the form
\begin{equation}
    \Sigma=\frac{2}{\nu(\mu^{-1}(x_0))}\sum_{k\in \mathscr{I}^r}L\left(\frac{x-\mu(y_k)}{\epsilon\nu(\mu^{-1}(x_0))}\right),\quad L(z):=\frac{1}{z^2+1}.
\label{eq:Sigma-rewrite}
\end{equation}
Let $k_0\in \mathscr{I}^r$ be the index such that $|x_0-\mu(y_k)|$ is minimized.  The spacing of the points $\mu(y_k)$ near $x_0$ is approximately $N_{L}(\epsilon)^{-1}\mu'(y_{k_0})$.  More precisely, since $\mu(y)$ is twice continuously differentiable according to Conjecture~\ref{conjecture: mu} with derivative $\mu'(y_{k_0})>0$,
\begin{equation}
\mu(y_{k_0})-\mu(y_k)=\mu'(y_{k_0})\frac{k_0-k}{N_L(\epsilon)}+\mathcal{O}\left(\left(\frac{k_0-k}{N_L(\epsilon)}\right)^2\right)
\label{eq:mu-Taylor}
\end{equation}
We next express $x_0$ in terms of the phase shift $p(x_0)$, which is defined by the relation
\begin{align}\label{def:px0}
    x_0=\mu(y_{k_0})+\mu'(y_{k_0})\frac{p(x_0)}{N_L(\epsilon)}.
\end{align}
The phase shift $p(x_0)$ must be bounded with an upper bound on  $|p(x_0)|$ close to $\frac{1}{2}$ because $k_0$ minimizes the distance between $x_0$ and $\mu(y_{k})$. Combining \eqref{eq:mu-Taylor} and \eqref{def:px0}, we write the argument of $L$ in the summand of \eqref{eq:Sigma-rewrite} as 
\begin{equation}
z=\frac{x-\mu(y_k)}{\epsilon\nu(\mu^{-1}(x_0))}=\tau + \frac{k_0-k}{ \Omega(x_0;\epsilon)}+\mathcal{O}\left(\frac{(k_0-k)^2}{N_L(\epsilon)}\right),
\end{equation}
where 
\begin{equation}
\tau:=\frac{x-x_0}{\epsilon\nu(\mu^{-1}(x_0))}+\frac{p(x_0)}{\Omega(x_0;\epsilon)},    
\end{equation}
and $\Omega(x_0;\epsilon)$ is the quantity
\begin{equation}
    \Omega(x_0;\epsilon):=\frac{\epsilon N_L(\epsilon)\nu(\mu^{-1}(x_0))}{\mu'(\mu^{-1}(x_0))},
\end{equation}
which has a finite nonzero limit $\Omega(x_0;0)$ as $\epsilon\to 0$.
We notice that the error term in the argument of $L$ will be small of order $\mathcal{O}(\epsilon^{1-2q})$ for indices $k\in \mathscr{I}^r$ for which $|k-k_0|\le\epsilon^{-q}$, where $q$ is any exponent with $0<q<\frac{1}{2}$.  Using also $r<\frac{1}{2}$, these terms will be a small fraction of the total cardinality $|\mathscr{I}^r|\sim\epsilon^{r-1}$ because $q+r<1$.  Since $z\mapsto L(z)$ obviously has a uniformly bounded derivative,
\begin{equation}
\begin{split}
    \sum_{\substack{k\in \mathscr{I}^r\\|k-k_0|\le\epsilon^{-q}}}L\left(\frac{x-\mu(y_k)}{\epsilon\nu(\mu^{-1}(x_0))}\right)&=\sum_{\substack{k\in \mathscr{I}^r\\|k-k_0|\le\epsilon^{-q}}}\left[L\left(\tau+\frac{k_0-k}{\Omega(x_0;\epsilon)}\right) + \mathcal{O}(\epsilon^{1-2q})\right] \\
    &= \sum_{\substack{k\in \mathscr{I}^r\\|k-k_0|\le\epsilon^{-q}}}L\left(\tau+\frac{k_0-k}{\Omega(x_0;\epsilon)}\right)+\mathcal{O}(\epsilon^{1-3q}),\quad\epsilon\to 0.
\end{split}
\label{eq:qsum}
\end{equation}
To ensure that the error term in \eqref{eq:qsum} is small we will now further constrain $q$ by assuming $q<\frac{1}{3}$.
Now, for $k$ in the complementary part of $\mathscr{I}^r$ where $|k-k_0|>\epsilon^{-q}$, we see that if we can guarantee the condition $x-x_0=o(\epsilon^{1-q})$ as $\epsilon\to 0$, we have $\tau=o((k_0-k)/\Omega(x_0;\epsilon))$ because also $p(x_0)/\Omega(x_0;\epsilon)$ is bounded.  Therefore,
\begin{equation}
    z=\frac{x-\mu(y_k)}{\epsilon\nu(\mu^{-1}(x_0))}=\frac{k_0-k}{\Omega(x_0;\epsilon)}(1+o(1)),
\end{equation}
so for these terms the argument $z$ of $L$ satisfies $|z|\gtrsim\epsilon^{-q}$.  Since $L(z)=\mathcal{O}(z^{-2})$ as $z\to\infty$, 
\begin{equation}
    \sum_{\substack{k\in \mathscr{I}^r\\|k-k_0|>\epsilon^{-q}}}L\left(\frac{x-\mu(y_k)}{\epsilon\nu(\mu^{-1}(x_0))}\right) = \sum_{\substack{k\in \mathscr{I}^r\\|k-k_0|>\epsilon^{-q}}}\mathcal{O}(\epsilon^{2q}) = \mathcal{O}(\epsilon^{2q+r-1}),
\label{eq:qbarsum}
\end{equation}
because these terms constitute the dominant fraction of those in $\mathscr{I}^r$, and $|\mathscr{I}^r|\sim\epsilon^{r-1}$.  For the same reasons, we also have
\begin{equation}
    \sum_{\substack{k\in \mathscr{I}^r\\|k-k_0|>\epsilon^{-q}}}L\left(\tau+\frac{k_0-k}{\Omega(x_0;\epsilon)}\right) = \mathcal{O}(\epsilon^{2q+r-1}),
\label{eq:qbartarget}
\end{equation}
so combining \eqref{eq:qsum}, \eqref{eq:qbarsum}, and \eqref{eq:qbartarget} we obtain
\begin{equation}
    \sum_{k\in \mathscr{I}^r}L\left(\frac{x-\mu(y_k)}{\epsilon\nu(\mu^{-1}(x_0))}\right)=\sum_{k\in \mathscr{I}^r}L\left(\tau+\frac{k_0-k}{\Omega(x_0;\epsilon)}\right) + \mathcal{O}(\epsilon^{1-3q})+\mathcal{O}(\epsilon^{2q+r-1}),\quad\epsilon\to 0.
\end{equation}
To guarantee that the term $\mathcal{O}(\epsilon^{2q+r-1})$ is negligible as $\epsilon\to 0$, we must put a lower bound on $q$, namely $q>\frac{1}{2}(1-r)$ which also implies $q>\frac{1}{4}$ because $r<\frac{1}{2}$.  With this inequality on $q$ we can verify that because $x-x_0=\mathcal{O}(\epsilon^{3/4})$ it also holds that $x-x_0=o(\epsilon^{1-q})$, exactly as presumed above.

Now that the argument of $L(z):=(z^2+1)^{-1}$ is linear in the index $k$, the corresponding infinite series summing over $k\in\mathbb{Z}$ is explicitly convergent \cite[1.445.9]{GRtable}:
\begin{equation}
\begin{split}
    \sum_{k\in\mathbb{Z}}L\left(\tau+\frac{k_0-k}{\Omega(x_0;\epsilon)}\right) &= \sum_{j\in\mathbb{Z}}L\left(\tau+\frac{j}{\Omega(x_0;\epsilon)}\right)\\ &=\frac{\pi\Omega(x_0;\epsilon)\sinh(2\pi\Omega(x_0;\epsilon))}{\cosh(2\pi\Omega(x_0;\epsilon))-\cos(2\pi\Omega(x_0;\epsilon)\tau)}.
\end{split}
\label{eq:series-identity}
\end{equation}
Hence if $0<r<\frac{1}{2}$ and $\frac{1}{4}<\frac{1}{2}(1-r)<q<\frac{1}{3}$, the bound $x-x_0=\mathcal{O}(\epsilon^{3/4})$ yields
\begin{equation}
    \sum_{k\in \mathscr{I}^r}L\left(\frac{x-\mu(y_k)}{\epsilon\nu(\mu^{-1}(x_0))}\right)=\frac{\pi\Omega(x_0;\epsilon)\sinh(2\pi\Omega(x_0;\epsilon))}{\cosh(2\pi\Omega(x_0;\epsilon))-\cos(2\pi\Omega(x_0;\epsilon)\tau)}+o(1),\quad\epsilon\to 0,
\end{equation}
and referring back to \eqref{eq:Sigma-def}--\eqref{eq:Sigma-rewrite}, this implies that $\Sigma$ is given by
\begin{multline}
   \Sigma= \sum_{k\in \mathscr{I}^r}\frac{2\epsilon^2\nu(\mu^{-1}(x_0))}{(x-\mu(y_k))^2+\epsilon^2\nu(\mu^{-1}(x_0))^2} \\{}= \frac{2\pi\Omega(x_0;\epsilon)}{\nu(\mu^{-1}(x_0))}\frac{\sinh(2\pi\Omega(x_0;\epsilon))}{\cosh(2\pi\Omega(x_0;\epsilon))-\cos(2\pi\Omega(x_0;\epsilon)\tau)} + o(1),\quad\epsilon\to 0.
\end{multline}
Under the same conditions on $r$ we have already seen that $D=o(1)$ (see \eqref{eq:D-def}) and that the right-hand side of \eqref{eq:complementary} is $o(1)$.  Combining these results with \eqref{uL in proof} and noting that $2\pi\Omega(x_0;\epsilon)=\phi_L(x_0,t)$ and $2\pi\Omega(x_0;\epsilon)\tau=\psi_L(x_0,t)\epsilon^{-1}(x-x_0)+2\pi p(x_0)$ (see \eqref{eq:psiL-phiL}) completes the proof.
\end{proof}

\begin{remark}
 The identity in \eqref{eq:series-identity} is a BO analogue of a corresponding identity for the KdV equation in which a periodic traveling wave solution given by the square of a Jacobi elliptic function can be expressed as a linear sum of translates of solitons proportional to the square of the hyperbolic secant function.  See \cite[22.11.14]{DLMF}.
\end{remark}

Combining Propositions \ref{prop: outliers}, \ref{prop: upper branch leading order}, and \ref{prop: lower branch leading order} gives us the complete result.

\begin{theorem}[Asymptotic expansion of the soliton ensemble]\label{thm:oscillations}
Let $u_0$ be an admissible initial condition with one inflection point to the right of the maximizer, and let $u(x,t)$ denote the corresponding soliton ensemble solving \eqref{BO equation}.
Suppose either that $0\le t<t_b$, or that $t\ge t_b$ and $x\in\mathbb{R}\setminus [X^-(t),X^+(t)]$, and that Conjectures~\ref{conjecture: outliers}, \ref{conjecture: mu}, and \ref{conjecture: nu upper} hold.  Then
\begin{equation}
    u(x,t)=\psi_U(x,t)+o(1),\quad\epsilon\to 0
\end{equation}
with the error term being uniform on compact subsets of the indicated domain.  On the other hand, if $t>t_b$ is fixed and $x_0\in (X^-(t),X^+(t))$, and Conjectures~\ref{conjecture: outliers}, \ref{conjecture: mu}, \ref{conjecture: nu upper}, and \ref{conjecture: nu lower} hold, then as $\epsilon\to 0$,
\begin{equation}
    u(x,t)=\psi_U(x,t)+\frac{\psi_L(x_0,t)\sinh(\phi_L(x_0,t))}{\cosh(\phi_L(x_0,t))-\cos(\psi_L(x_0,t)\epsilon^{-1}(x-x_0)+2\pi p(x_0))}+o(1)
    \label{eq: approx after breaking}
\end{equation}
with the error term being uniform for $x-x_0=\mathcal{O}(\epsilon^{3/4})$.
\end{theorem}

\subsection{Comparison with expectations of Whitham modulation theory}
The explicit terms on the right-hand side of \eqref{eq: approx after breaking} give, for each fixed $(x_0,t)$ a periodic function of $x$ that is the profile $f$ of an exact traveling wave solution $u(x,t)=f(x-ct)$ of the BO equation \eqref{BO equation}; see \cite{Benjamin1967,Ono1975}.  The periodic wave has wavelength proportional to $\epsilon$, and the approximation asserted in Theorem~\ref{thm:oscillations} is valid for a range of values of $x-x_0$ that includes a large number ($\mathcal{O}(\epsilon^{-1/4})$) of wavelengths.  Rapidly oscillatory time dependence of the wave enters implicitly via the phase $p(x_0)$.  Resolving the time dependence of $p(x_0)$ would require further conjectures regarding that of the eigenvalues $\sigma_k(t)$ going beyond the scope of our paper.  However, as the parameters $(x_0,t)$ vary within the domain $t>t_b$ with $X^-(t)<x_0<X^+(t)$, the parameters of the periodic wave vary as well, which means that the soliton ensemble is actually a relatively slowly modulated periodic wavetrain.

\begin{figure}[h]
\begin{center}
    \includegraphics[width=0.45\linewidth]{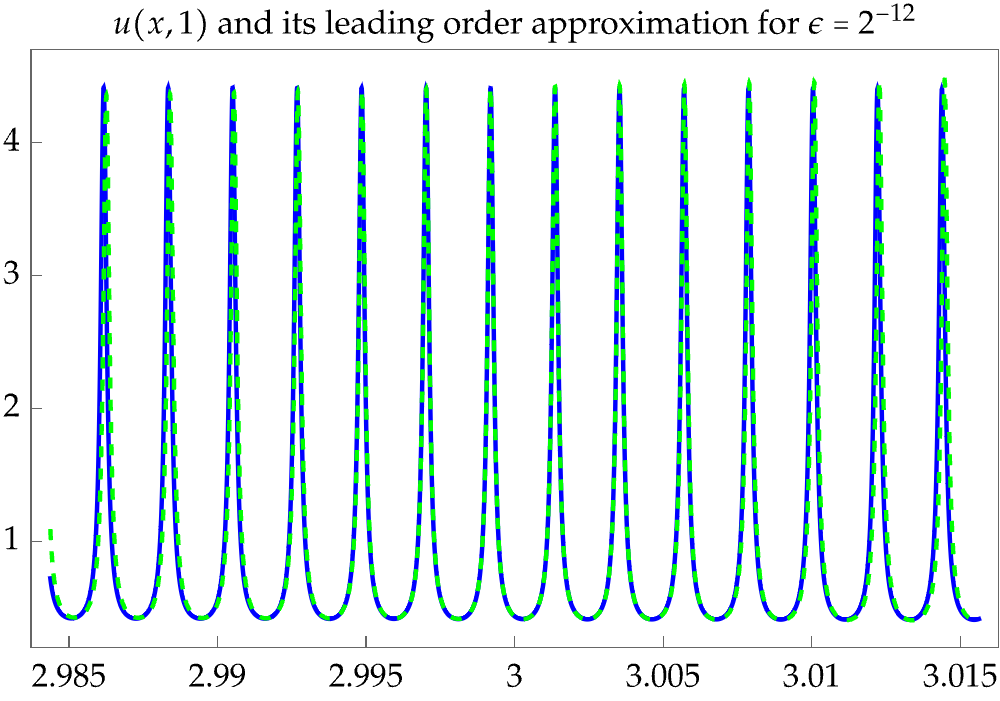}\hfill%
    \includegraphics[width=0.45\linewidth]{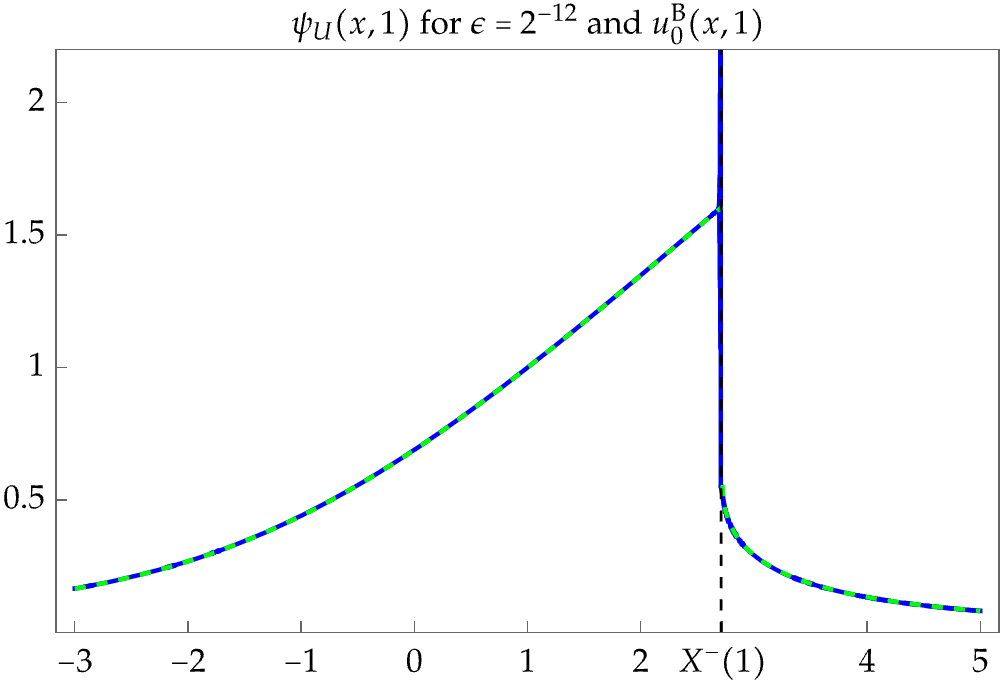}\hfill\\
    \includegraphics[width=0.45\linewidth]{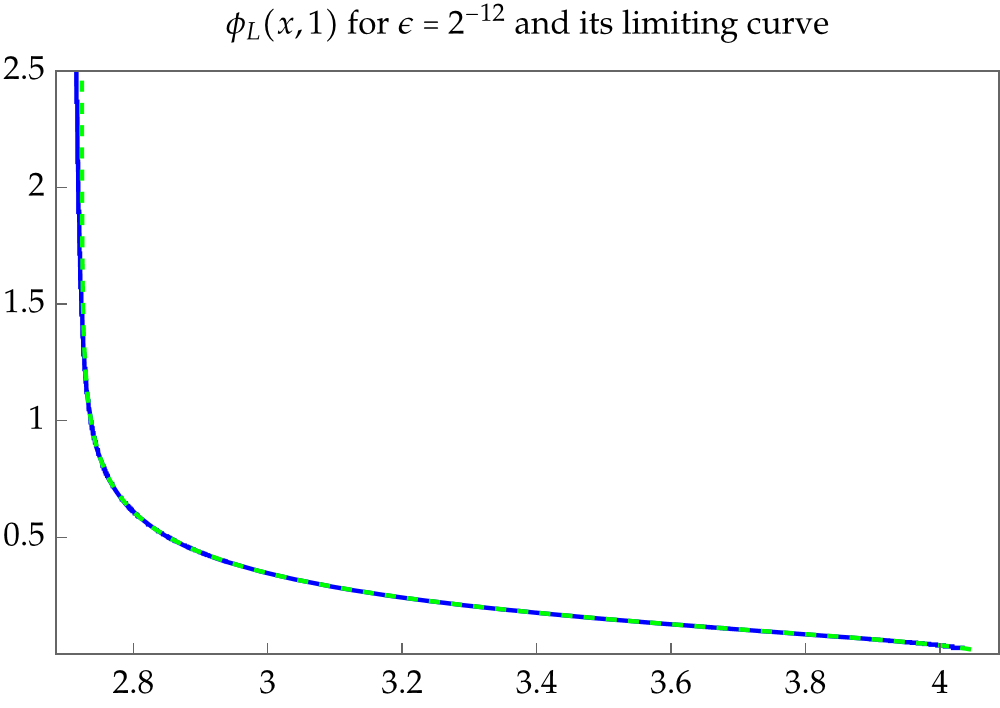}\hfill%
    \includegraphics[width=0.45\linewidth]{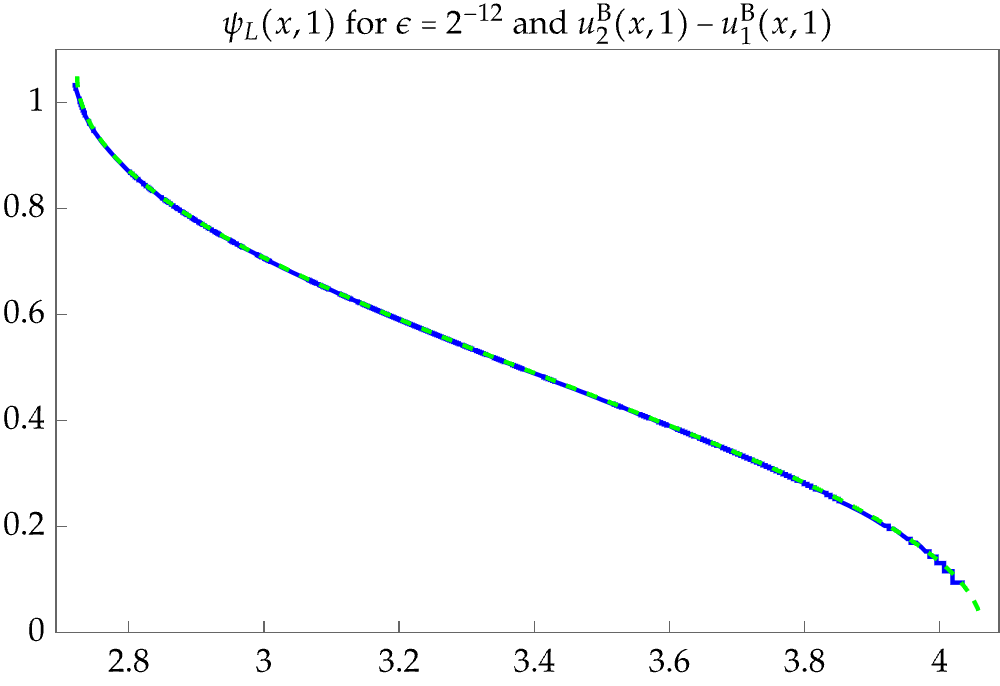}\hfill
\end{center}
\caption{The theoretical asymptotic expressions  for the soliton ensemble $u(x,t)$ (\eqref{eq: approx after breaking} with $x_0=3$) and modulation fields $\psi_{U,L}(x,t)$, and $\phi_{L}(x,t)$ (see \eqref{eq:psiU-B}--\eqref{eq:psiphiL-B}) shown with dashed green curves compared with the soliton ensemble for fixed nonzero $\epsilon=2^{-12}$ (see \eqref{u in terms of sigmas}) and corresponding numerical approximations of the modulation fields obtained from eigenvalues of $\mathbf{C}(t)$ and the coincident discretizations of the functions $\mu_{U,L}(y)$, $\mu_{U,L}'(y)$, and $\nu_{L}(y)$, shown with blue curves.}
\label{fig: limiting curves}
\end{figure}

According to the formal Whitham modulation theory for the BO equation developed by Dobrokhotov and Krichever in \cite{DobrokhotovK91}, modulated periodic waves should have the form given in \eqref{eq: approx after breaking}, but the quantities $\psi_U$, $\psi_L$, and $\phi_L$ should be given as functions of $(x_0,t)$ in terms of three Riemann invariants solving a system of uncoupled Burgers equations. We identify these Riemann invariants with the branches of the multi-valued solution of Burgers' equation with initial data $u_0$.  Matching the formula \eqref{eq: approx after breaking} with the Whitham theory 
requires the following explicit identifications:

\begin{align}
    \psi_{U}(x,t)=u_0^\mathrm{B}(x,t),
    \label{eq:psiU-B}
\end{align}
and, for $(x,t)$ in the multi-valued Burgers region,
\begin{align}
    \psi_{L}(x,t)=u_2^\mathrm{B}(x,t)-u_1^\mathrm{B}(x,t), ~~~ \phi_{L}(x,t)=\frac{1}{2}\ln\left(\frac{u_2^\mathrm{B}(x,t)-u_0^\mathrm{B}(x,t)}{u_1^\mathrm{B}(x,t)-u_0^\mathrm{B}(x,t)}\right).
    \label{eq:psiphiL-B}
\end{align}
We do not have direct proof of either \eqref{eq:psiU-B} or \eqref{eq:psiphiL-B}, although we do have solid numerical evidence that both are true.  Indeed, the numerical approximation of the multiscale soliton ensemble $u(x,t)$, and the slowly varying modulation parameters $\psi_{L,U}(x_0,t)$, and $\phi_{L}(x_0,t)$ obtained from the eigenvalues $\sigma_k(t)$ are plotted along with the expressions above in Figure \ref{fig: limiting curves}.  
This shows that the distribution of complex eigenvalues $\sigma_k(t)$ of $\mathbf{C}(t)$ indeed produces an approximate formula \eqref{eq: approx after breaking} for the BO soliton ensemble associated with the initial data $u_0$ that is fully consistent with Whitham modulation theory.  We formalize these expectations in a final conjecture.
\begin{conjecture}[Connection with Whitham modulation theory]
The function $\psi_U(x,t)$ defined in terms of the scaling functions of the upper-branch eigenvalues by \eqref{eq:psiU} and the functions $\psi_L(x,t)$ and $\phi_L(x,t)$ defined in terms of those of the lower-branch eigenvalues by \eqref{eq:psiL-phiL} are related to the generally multi-valued solution of Burgers' equation with initial data $u_0$ by \eqref{eq:psiU-B}--\eqref{eq:psiphiL-B}.
\end{conjecture}

\section{Asymptotic properties of the real eigenvalues \texorpdfstring{$\alpha_k(x,t)$}{alpha-k} and their implications}
\label{sec:A}
In this section, we study the implications of a  different formula for $u(x,t)$ in terms of the eigenvalues of the matrix $\mathbf{A}(x,t)$.  The approach in this section is completely independent of that presented in Section~\ref{sec:C}.

\subsection{Estimates on derivatives of the eigenvalues and the importance of small eigenvalues}

Let $\alpha_k(x,t)$, $k=1,\ldots,N$, denote the eigenvalues of $\mathbf{A}(x,t)$ and let $\mathbf{u}_k(x,t)$ denote the corresponding orthonormalized eigenvectors.  
Recall the formula \eqref{u in terms of alphas} for $u(x,t)$ in terms of the eigenvalues of $\mathbf{A}(x,t)$ and their derivatives with respect to $x$. By differentiation of the relation $\mathbf{A}(x,t)\mathbf{u}_k(x,t)=\alpha_k(x,t)\mathbf{u}_k(x,t)$ and noting that the coordinates $(x,t)$ only occur in the diagonal entries of $\mathbf{A}(x,t)$, we find that 
the $x$ and $t$ derivatives of $\alpha_k(x,t)$ can be expressed in terms of the eigenvector $\mathbf{u}_k(x,t)$ as
\begin{equation}\label{eq:alpha-x}
\frac{\partial\alpha_k}{\partial x}(x,t)=\sum_{j=1}^{N(\epsilon)}(-2\lambda_j)|u_{k,j}(x,t)|^2,
\quad \frac{\partial\alpha_k}{\partial t}(x,t)=\sum_{j=1}^{N(\epsilon)}(-4\lambda_j^2)|u_{k,j}(x,t)|^2.
\end{equation}
Since $-L<\lambda_j<0$ for $j=1,\ldots,N(\epsilon)$, it is clear that
\begin{align}\label{eq:bound_alpha_x}
    0<\frac{\partial\alpha_k}{\partial x}(x,t)<2L, ~~~ -4L^2<\frac{\partial\alpha_k}{\partial t}(x,t)<0.
\end{align}
We next note that only small values of $\alpha$ actually contribute significantly to $u(x,t)$.  Indeed, we know from~\eqref{eq:bound_alpha_x} that $0\leq \partial_x\alpha_k\leq 2L$, hence
\begin{equation}\label{eq:large-A-eigenvalues}
\sum_{|\alpha_k|\geq \varepsilon^{r}}\frac{2\partial_x\alpha_{k}(x,t)}{\left(\epsilon^{-1}\alpha_k(x,t)\right)^2 + 1}
	\leq 4L\sum_{|\alpha_k|\geq \varepsilon^{r}}\frac{1}{\varepsilon^{2r-2}}\le 4LN(\epsilon)\epsilon^{2-2r}
	\lesssim \varepsilon^{1-2r}.
\end{equation}
This goes to zero as $\varepsilon\to 0$ when $r<\frac{1}{2}$.

\subsection{Numerical experiments}
\begin{figure}[h]
\begin{center}
    \includegraphics[width=0.3\linewidth]{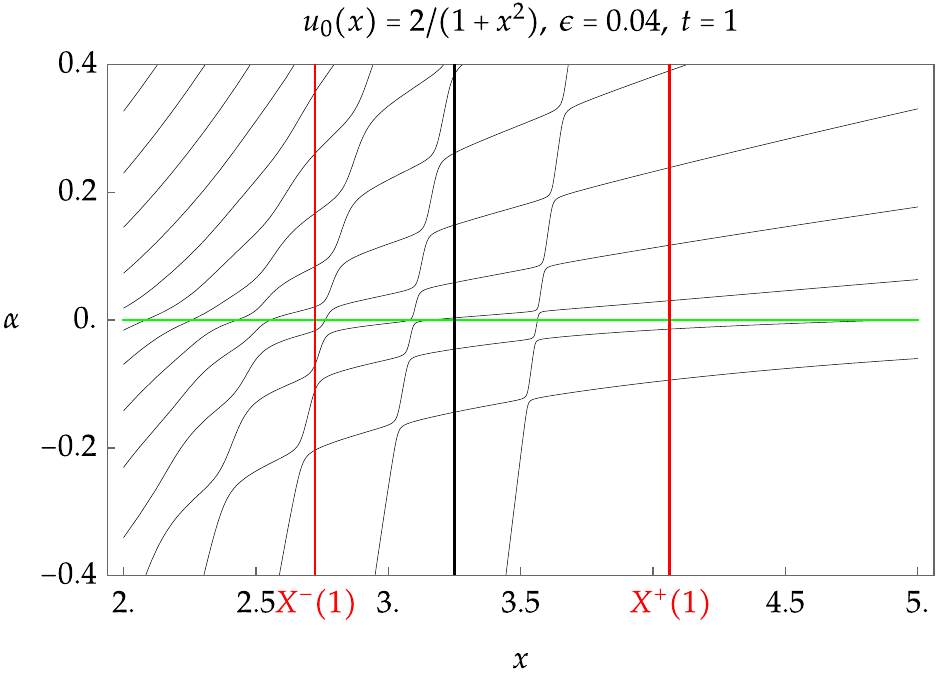}\hfill%
    \includegraphics[width=0.3\linewidth]{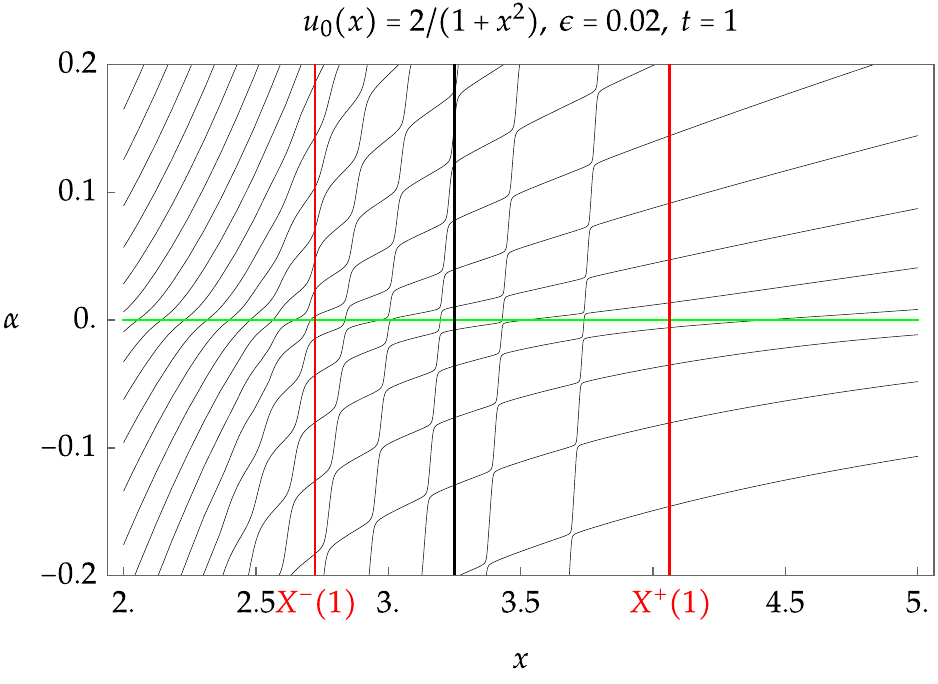}\hfill%
    \includegraphics[width=0.3\linewidth]{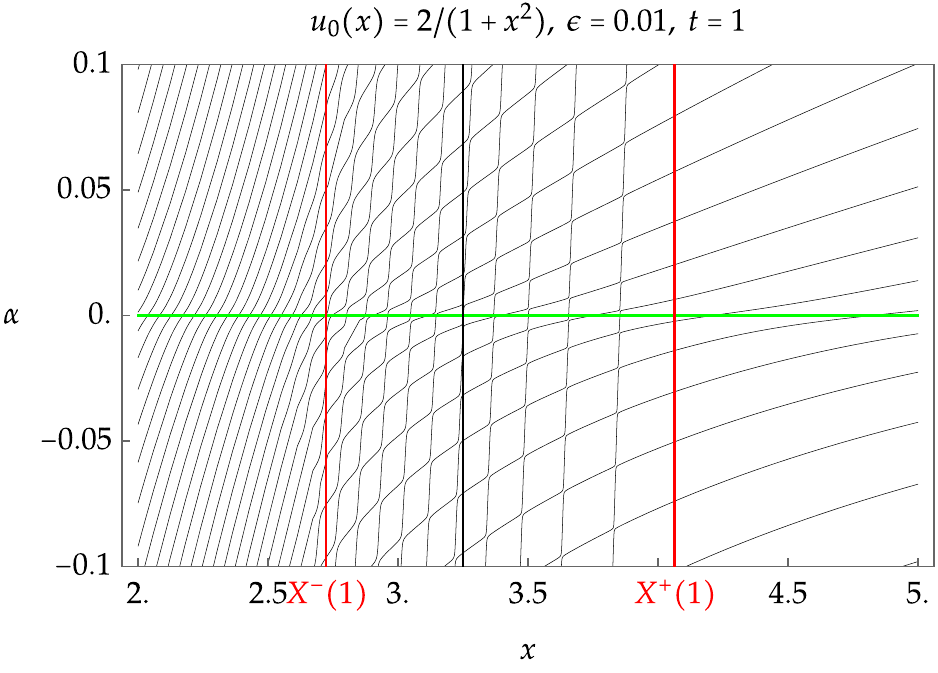}
\end{center}
\caption{Numerical plots of small eigenvalues $\alpha_k(x,1)$ versus $x$, for three different values of $\epsilon$.  The value of $t=1>t_b$ is selected so that for the range of $x$ in the plots, there is an interval $X^-(t)<x<X^+(t)$ delineated with vertical red lines on which the solution of Burgers' equation is multi-valued.  It is clear that for $x\in (X^-(t),X^+(t))$, there are both slow-moving and fast-moving eigenvalues whose trajectories actually form a system of non-intersecting paths. (The vertical line at $x=3.25$ is for reference only; see Figure~\ref{fig:alpha-x} below.)}
\label{fig:alphas}
\end{figure}

We first observe the distribution of small eigenvalues of $\mathbf{A}(x,t)$ and their corresponding (normalized) eigenvectors  via numerics.  See Figure~\ref{fig:alphas}.
These computations suggest that on subintervals of $x\in(X^-(t),X^+(t))$ where the method-of-characteristics solution of the inviscid Burgers equation with initial data $u_0(x)$ is multi-valued at some fixed time $t>t_b$, most small eigenvalues $\alpha_k(x,t)$ have either a specific large $x$-velocity $\alpha_{k,x}(x,t)$ or a small $x$-velocity.  The small eigenvalues $\alpha_k(x,t)$ also appear to never coincide for any $x$; hence there are also numerous near-collisions between fast-moving and slow-moving eigenvalues in which what actually occurs is that a slow eigenvalue rapidly accelerates while a fast eigenvalue decelerates.  On the other hand, if there is only one branch of the Burgers solution above $(x,t)$, then all of the small eigenvalues appear to be of the slow-moving type.  

The overall distribution of the eigenvalues $\alpha_k(x,t)$ in the limit $\epsilon\to 0$ is known according to \cite{MillerX11}.  Indeed, the following limit holds in the weak-$*$ sense:
\begin{equation}
\begin{split}
    \lim_{\epsilon\to 0}\frac{M}{N(\epsilon)}\sum_{k=1}^{N(\epsilon)}\delta_{\alpha_k(x,t)}(\alpha) &= G(\alpha;x,t)\,\dd\alpha,\\
   G(\alpha;x,t)&:=-\frac{1}{4\pi}\int_{-L}^0\chi_{[-2\lambda(x+2\lambda t-x_+(\lambda)),-2\lambda(x+2\lambda t -x_-(\lambda))]}(\alpha)\frac{\dd\lambda}{\lambda},
\end{split}
\label{eq:distribution-of-alphas}
\end{equation}
where $\chi_{I}$ denotes the characteristic function of the interval $I$, and $x_-(\lambda)<x_+(\lambda)$ are the two turning points satisfying $u_0(x_\pm(\lambda))=-\lambda$.  So, given $\alpha\in\mathbb{R}$, the limiting density of eigenvalues near $\alpha$ is computed as $-1/(4\pi)$ times the integral of $1/\lambda$ over those subintervals of $\lambda\in (-L,0)$ where the inequalities $-2\lambda(x+2\lambda t-x_+(\lambda))<\alpha<-2\lambda(x+2\lambda t-x_-(\lambda))$ both hold.  In the case that the inviscid Burgers solution with initial data $u_0(x)$ has three branches, $u_0^\mathrm{B}(x,t)<u_1^\mathrm{B}(x,t)<u_2^\mathrm{B}(x,t)$, and that $\alpha$ is small, there are two such subintervals:  $[-u_2^\mathrm{B}(x,t)+o(1),-u_1^\mathrm{B}(x,t)+o(1)]$ and $[-u_0^\mathrm{B}(x,t)+o(1),\eta(\alpha)]$, where $\eta(\alpha)\le 0$ is a small quantity that satisfies the implicit equation 

\begin{equation}
    u_0\left(\frac{\alpha}{2\eta(\alpha)}+x+2\eta(\alpha) t\right)=-\eta(\alpha).
\label{eq:implicit-lambda-alpha}
\end{equation}
For positive rational $u_0(x)$ there is an integer $p>0$ and a constant $C>0$ such that $u_0(x)=Cx^{-2p}(1+\mathcal{O}(x^{-1}))$ as $x\to\pm\infty$.  Solutions $\eta(\alpha)$ of \eqref{eq:implicit-lambda-alpha} that are small as $\alpha\to 0$ necessarily satisfy $\eta(\alpha)\ll\alpha$, in which case the above large-$x$ approximation of $u_0$ yields that 
\begin{equation}
    \eta(\alpha)=-K|\alpha|^{2p/(2p-1)}(1+o(1)),\quad \alpha\to 0,\quad K:=\left(\frac{1}{2^{2p}C}\right)^{1/(2p-1)}.
\end{equation}
Therefore, for small $\alpha$, the density $G(\alpha;x,t)$ is approximated by
\begin{equation}
\begin{split}
    G(\alpha;x,t)&=-\frac{1}{4\pi}\left[\ln(u^\mathrm{B}_1(x,t)+o(1))-\ln(u^\mathrm{B}_2(x,t)+o(1))\right.\\
    &\qquad\qquad\qquad\left.{}+\ln(-\eta(\alpha))-\ln(u^\mathrm{B}_0(x,t)+o(1))\right]\\
    &=\frac{1}{4\pi}\frac{2p}{2p-1}\ln(|\alpha|^{-1}) + \mathcal{O}(1),\quad \alpha\to 0.
\end{split}
\label{eq:G-small-alpha}
\end{equation}
Thus, the overall asymptotic density of eigenvalues $\alpha_k(x,t)$ diverges logarithmically as $\alpha\to 0$.  There are therefore many small eigenvalues $\alpha_k(x,t)$, and when $x\in(X^{-}(t),X^{+}(t))$, the plots in Figure~\ref{fig:alphas} suggest that the majority of these are slow-moving eigenvalues\footnote{The majority of the intersections with any given vertical reference line such as $x=3.25$ in Figure~\ref{fig:alphas} are evidently with curves having the smaller of the two slopes.  On the other hand, the majority of the intersections with the horizontal line $\alpha=0$ as $x$ varies in the multi-valued interval are with curves having the larger of the two slopes, i.e. most of the eigenvalues crossing the origin $\alpha=0$ with varying $x$ are of the fast-moving variety.  Note that according to \eqref{eq:A-to-B}, we have $\mathbf{A}(x,t)=\mathbf{D}(x\mathbb{I}-\mathbf{B}(t))\mathbf{D}$, and hence $\mathbf{A}(x,t)$ has $\alpha=0$ as an eigenvalue precisely when $x$ is an eigenvalue of $\mathbf{B}(t)$.}.  One may think of the origin $\alpha=0$ as locating a kind of ``traffic jam'' of eigenvalues with small positive $x$-velocities through which a small number of fast-moving eigenvalues repeatedly pass with increasing $x$.

\begin{figure}
\begin{center}
    \includegraphics[width=0.3\linewidth]{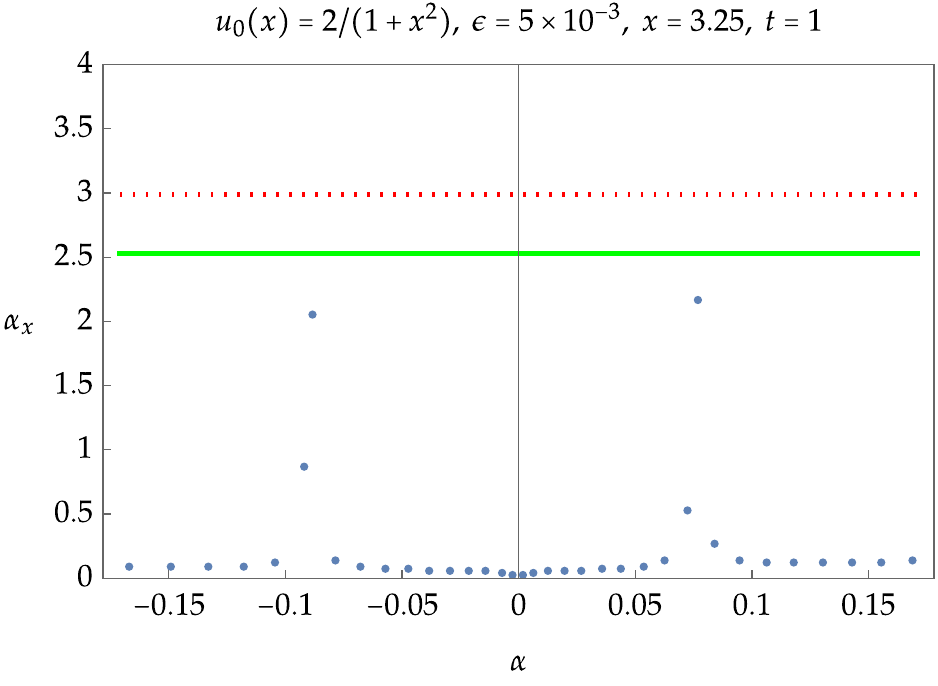}\hfill%
    \includegraphics[width=0.3\linewidth]{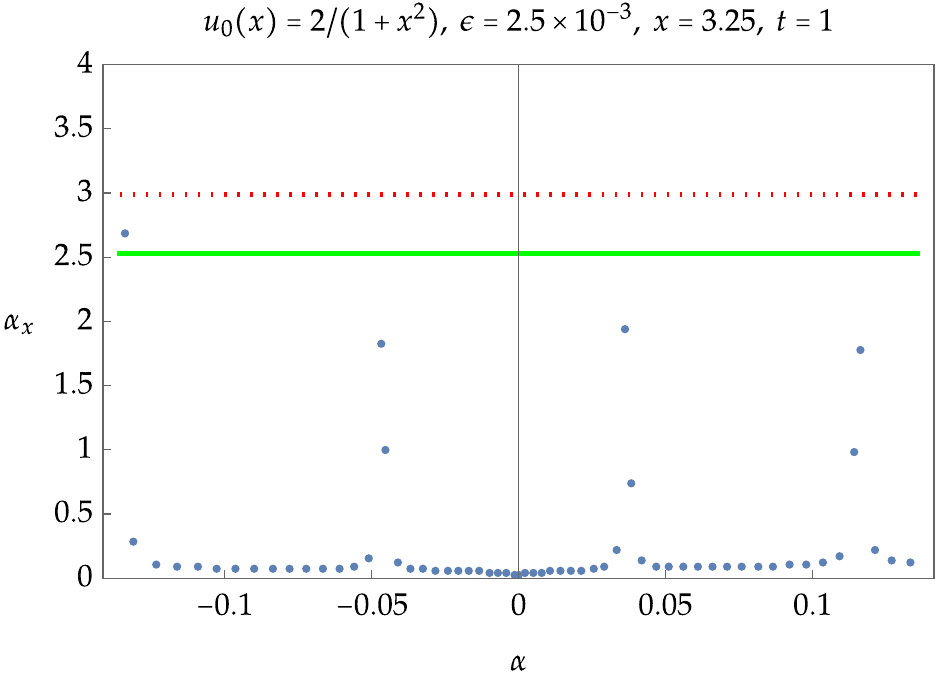}\hfill%
    \includegraphics[width=0.3\linewidth]{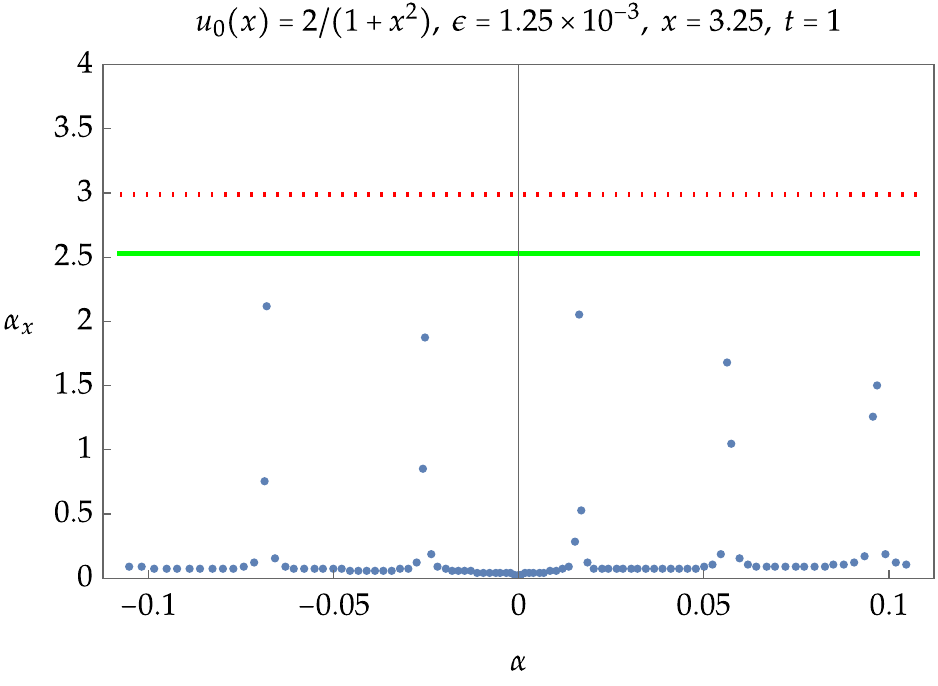}\\
    \includegraphics[width=0.3\linewidth]{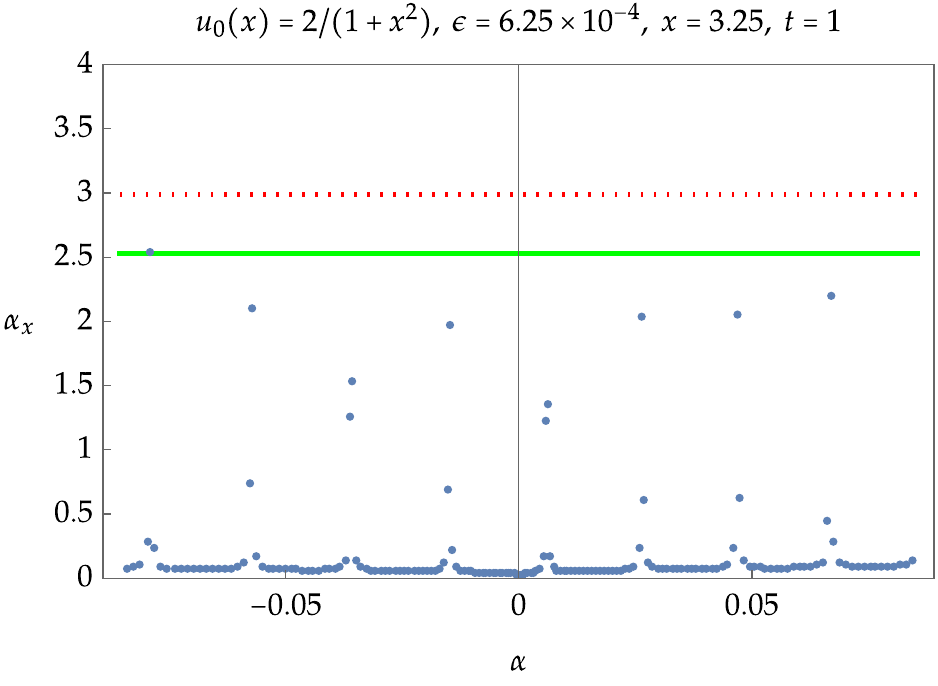}\hfill%
    \includegraphics[width=0.3\linewidth]{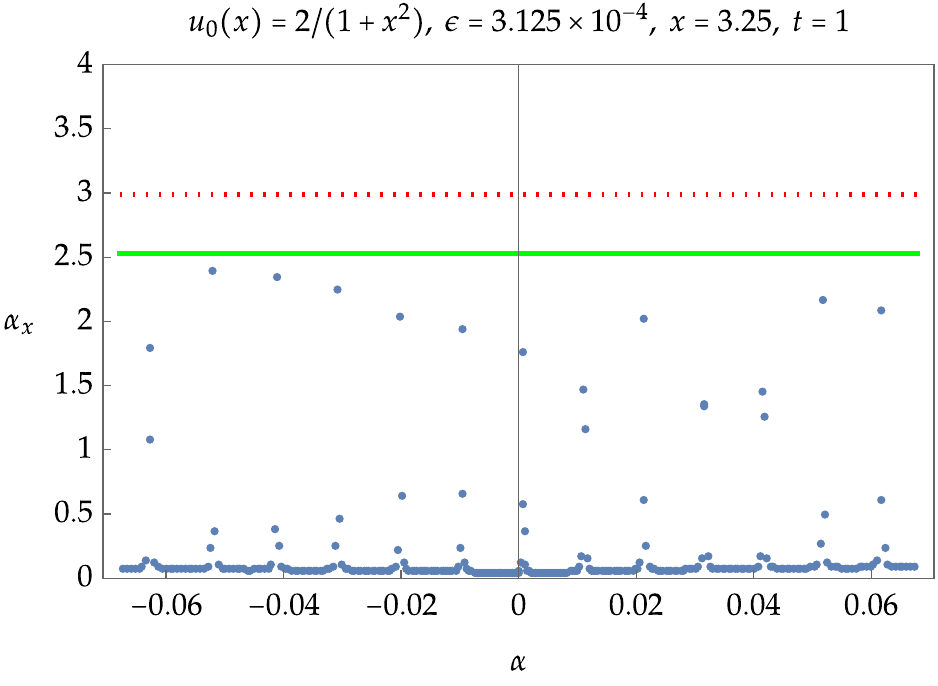}\hfill%
    \includegraphics[width=0.3\linewidth]{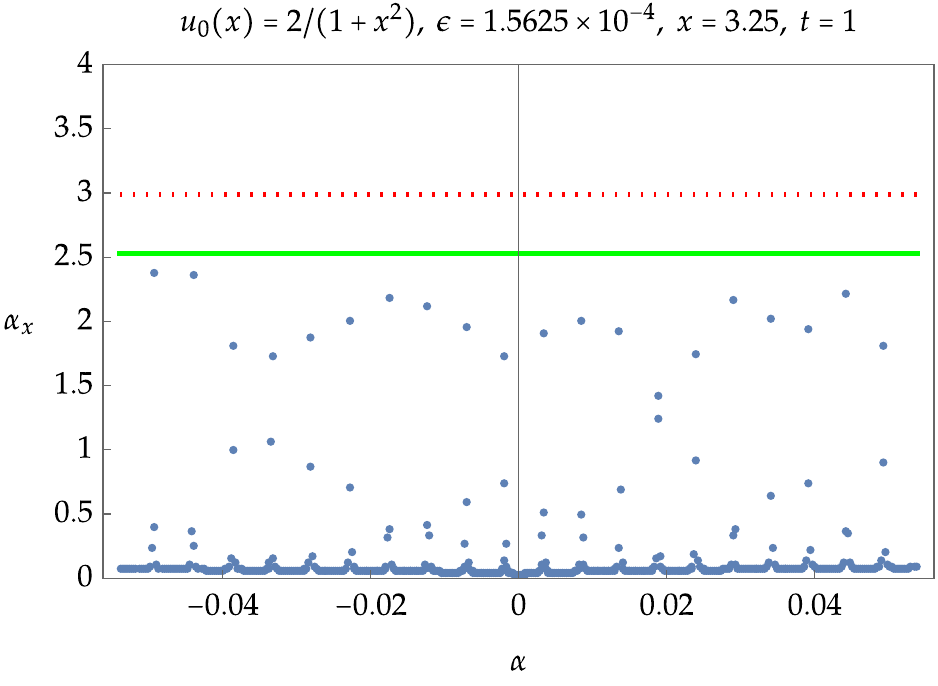}
\end{center}
\caption{The eigenvalue velocity $\alpha_{k,x}(x,t)$ determined from the eigenvector according to \eqref{eq:alpha-x} plotted against the eigenvalues $\alpha_k(x,t)$ in the range $|\alpha|<\epsilon^{1/3}$, for $(x,t)=(3.25,1)$ and various indicated values of $\epsilon$.  For the indicated initial condition, $L=2$, so all derivatives $\alpha_x$ lie in the range $0<\alpha_x<2L=4$.  The dotted red line is $\alpha_x=2(u_2^\mathrm{B}(x,t)-u_1^\mathrm{B}(x,t))/(\ln(u_2^\mathrm{B}(x,t))-\ln(u_1^\mathrm{B}(x,t)))$, and the green line is $\alpha_x=2(u_2^\mathrm{B}(x,t)-u_1^\mathrm{B}(x,t))/(\ln(u_2^\mathrm{B}(x,t)-u_0^\mathrm{B}(x,t))-\ln(u_1^\mathrm{B}(x,t)-u_0^\mathrm{B}(x,t)))$.}
\label{fig:alpha-x}
\end{figure}

The $x$-velocity of an eigenvalue $\alpha=\alpha_k(x,t)$ can be expressed explicitly in terms of its corresponding normalized eigenvector $\mathbf{u}_k(x,t)$ by \eqref{eq:alpha-x}.  One can then calculate the $x$-velocities of all the (small) eigenvalues directly from numerically computed eigenvectors.  Selecting the value of $x$ at $t=1$ corresponding to the vertical line at $x=3.25$ in the plots in Figure~\ref{fig:alphas}, the velocities of the eigenvalues in the range $|\alpha|<\epsilon^{1/3}$ (larger eigenvalues make a negligible contribution to the sum in \eqref{u in terms of alphas} by \eqref{eq:large-A-eigenvalues}) are plotted against the eigenvalues $\alpha$ for a series of decreasing values of $\epsilon$ in Figure~\ref{fig:alpha-x}.

These plots suggest that as $\epsilon\to 0$, the ``fast'' small eigenvalues are in the minority, they have a regular spacing proportional to $\epsilon$, and (perhaps) their velocities approach a limiting value in the range $(0,2L)$.  Two possible limiting values are shown in the plots as horizontal lines; clearly the solid green line is a better fit than the dotted red line.  

\begin{figure}[h]
\begin{center}
    \includegraphics[width=0.45\linewidth]{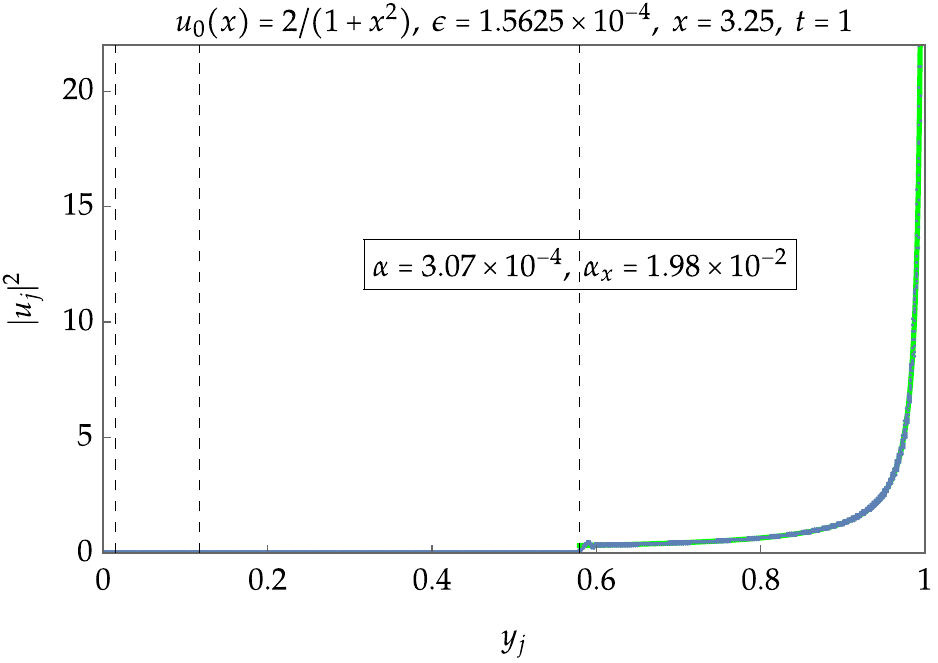}\hfill%
    \includegraphics[width=0.45\linewidth]{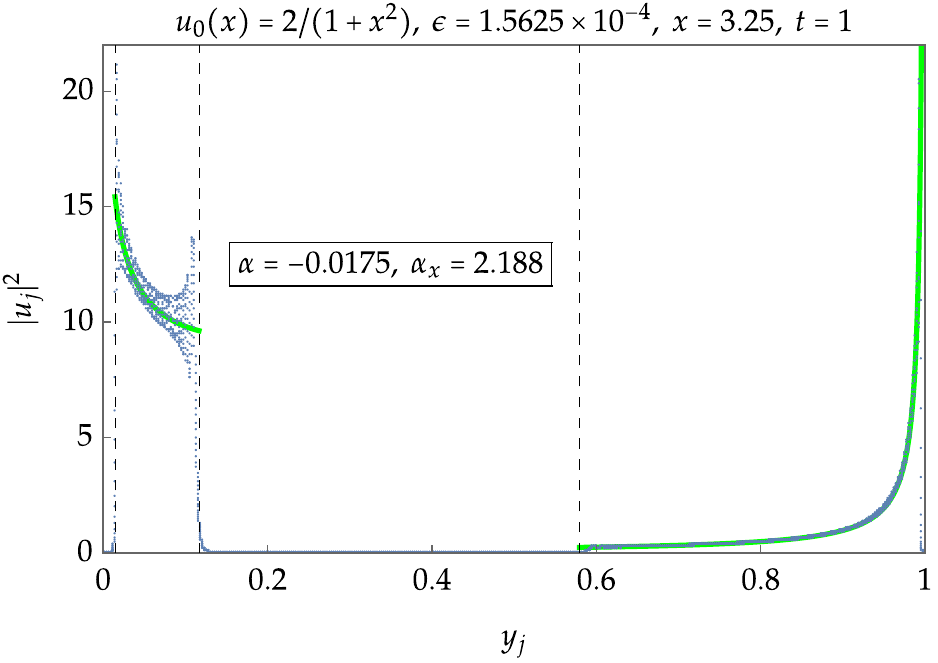}
\end{center}
\caption{Left: an (here, unnormalized) eigenvector of $\mathbf{A}(x,t)$ for a ``slow'' eigenvalue.  Right:  an unnormalized eigenvector of $\mathbf{A}(x,t)$ for a ``fast'' eigenvalue.  Values of $|u_j|^2$ are shown with blue points.  Also shown with green curves are best fits (choice of $a_0>0$) on the intervals bounded by the $y$-values corresponding to the three values of the multi-valued solution of the inviscid Burgers equation (indicated with dashed vertical lines) to the approximate squared amplitude $a_0^2\Lambda'(y)/(-\Lambda(y))$ (see \eqref{eq:amplitude} below).  For the ``slow'' eigenvalue, the best fit is for $a_0^2\approx 6.461\times 10^{-6}$.  For the ``fast'' eigenvalue, the best fit for the left support interval is $a_0^2\approx 2.983\times 10^{-4}$ and the best fit for the right support interval is $a_0^2\approx 4.873\times 10^{-6}$; the fraction of the squared $\ell^2$ norm in the right support interval is approximately $0.7543$.}
\label{fig:eigenvectors of A}
\end{figure}

To explain the predictions behind the horizontal lines in the plots of Figure~\ref{fig:alpha-x}, we may look at and compare plots of the square modulus of components of an eigenvector $\mathbf{u}_k(x,t)$ for slow and fast eigenvalues.  See Figure~\ref{fig:eigenvectors of A}.
The plots clearly show that whether a small eigenvalue is ``slow'' or ``fast'', the eigenvector is strongly localized in two subintervals of the rescaled index $y_j=\epsilon (j-\frac{1}{2})$, which generally lies in the range $0<y_j<M$ as $j\in\mathbb{Z}$ varies from $j=1$ to $j=N(\epsilon)$.  One of these intervals abuts the right edge $y_j=M$ but the other is bounded away from this edge.  The ``slow'' eigenvector appears to be supported in the abutting interval while the ``fast'' eigenvector evidently has some support in both intervals.  Moreover, there is evidence in the plots that near the right edge, $|u_{k,j}(x,t)|^2$ is proportional to $(M-y_j)^{-1}$, which is not integrable as a function of $y_j$.  Therefore, if there is any support of $\mathbf{u}_k(x,t)$ in the interval abutting this singularity, the constant of proportionality must be very small for the eigenvector to be normalized.  Now, recalling that as $j$ varies from $j=1$ to $j=N(\epsilon)$ the numbers $\lambda_j$ increase monotonically from $\lambda_1\approx -L$ to $\lambda_{N(\epsilon)}\approx 0$, and in fact one can show that for $u_0(x)=2(1+x^2)^{-1}$, $\lambda_j=\mathcal{O}(M-y_j)$ near the right edge.  Consequently, for such eigenvectors the formula \eqref{eq:alpha-x} predicts a small value of $\partial_x\alpha_{k}(x,t)$.  

This means that in order for an eigenvalue $\alpha_k(x,t)$ to have a velocity $\partial_x\alpha_k(x,t)$ that is not small in the limit $\epsilon\to 0$, it is necessary for the corresponding eigenvector to be predominantly localized in the other subinterval of $(0,M)$ that does not abut the right edge.  This could happen in two different ways:
\begin{enumerate}
    \item In the limit $\epsilon \to 0$, the eigenvector could have \emph{no support} on any interval abutting the right edge $y_j=M$.  If this is true, then the plot in the right-hand panel of Figure~\ref{fig:eigenvectors of A} is misleading in the sense that if $\epsilon$ is made smaller the evident support near the right edge should disappear rapidly.
    \item In the limit $\epsilon\to 0$, the eigenvector could have a nonzero limiting fraction of its norm in the subinterval abutting the right edge; since $|u_{k,j}(x,t)|^2$ is proportional to $(M-y_j)^{-1}$ and $M-y_{N(\epsilon)}=\mathcal{O}(\epsilon^{1/2})$ as $\epsilon \to 0$ for $u_0(x)=2(1+x^2)^{-1}$, this suggests that the proportionality constant on the abutting interval should be small of size $\mathcal{O}(\ln(\epsilon^{-1})^{-1})$.
\end{enumerate}
To see the full implications of these two alternatives requires an asymptotic theory of eigenvectors of $\mathbf{A}(x,t)$ that we will develop next.  However, we may point out at this juncture that the dotted red line in the panels of Figure~\ref{fig:alpha-x} corresponds to the first case, while the green line corresponds to a specific choice of small proportionality constant in the second case, selected to match the formal predictions of Whitham modulation theory \cite{DobrokhotovK91}.

\subsection{Microlocal analysis of eigenvectors of \texorpdfstring{$\mathbf{A}(x,t)$}{A}}
As suggested by the above numerical observations, we distinguish two families of ‘‘fast'' and ‘‘slow'' eigenvectors. An explanation for these two types of  eigenvectors may lie in semiclassical analysis, and more precisely, on Toeplitz (or Berezin-Toeplitz) quantization. Indeed, one might notice that the matrix $\mathbf{A}(x,t)$ resembles a generalized Toeplitz matrix, whose entries vary slowly along the diagonals.

To exhibit this structure more clearly, let a monotone increasing function $\Lambda:(0,M)\to(-L,0)$ (see \eqref{eq:LM} for $L$, $M$) be defined as
\begin{align}\label{def: big lambda}
    \int_{-L}^{\Lambda(y)}F(\lambda)\,\dd\lambda = y.
\end{align}
It follows that $\Lambda'(y)=1/F(\Lambda(y))$ and $\lambda_j=\Lambda(\epsilon(j-\frac{1}{2}))$, see \eqref{def: tilde lambda n}.  We can write the off-diagonal elements of $\mathbf{A}(x,t)$ as
\begin{equation}
A_{jk}(x,t)=\frac{2\Lambda(y_j)}{\Lambda'(y_j)}\frac{\ii}{k-j} + \epsilon f(y_j,\Delta y_{jk}),\quad k\neq j
\label{eq:BOD-exact}
\end{equation}
in which 
\begin{equation}
    y_j:=\epsilon\left(j-\frac{1}{2}\right),\quad
\Delta y_{jk}:=\epsilon (k-j),
\end{equation}
\begin{equation}
f(y,\Delta y):=
-2\ii\Lambda(y)\frac{\Lambda(y + \Delta y)-\Lambda(y)-\Lambda'(y)\Delta y \sqrt{\frac{\Lambda(y+\Delta y)}{\Lambda(y)}}}{ \Lambda'(y)\Delta y[\Lambda(y+\Delta y)-\Lambda(y)]}.
\end{equation}
Both the numerator and denominator of $f(y,\Delta y)$ are smooth $\epsilon$-independent functions of $\Delta y$ that vanish to second order at $\Delta y=0$.  In fact, we can let $\Delta y\to 0$ and obtain the limiting value
\begin{equation}
    f(y,0)= \ii\frac{\Lambda'(y)^2-\Lambda(y)\Lambda''(y)}{\Lambda'(y)^2}=\ii\frac{\dd}{\dd y}\frac{\Lambda(y)}{\Lambda'(y)}.
\label{eq:f(y,0)}
\end{equation}
Note that the diagonal elements of $\mathbf{A}(x,t)$ can also be expressed as the sampling of a smooth function:
\begin{equation}
A_{jj}(x,t)
	=-2\Lambda(y_j)(x+2\Lambda(y_j)t+\gamma(\Lambda(y_j))).
 \label{eq:A-diagonal-elements}
\end{equation}

Numerics shown in Figure~\ref{fig:eigenvectors of A} and in Figure~\ref{fig:eigenvectorphase} below suggest that on ranges of indices $j$ where eigenvector elements $u_{k,j}(x,t)$ of $\mathbf{A}(x,t)$ are not small, they have slowly varying amplitude and rapidly oscillating phase.  Heuristically, we therefore propose a wavepacket approximation of the eigenvectors $\mathbf{u}_k$ via a WKB-type expansion. The proof of the following proposition relies on a hypothetical but reasonable estimate; see \eqref{eq:big-claim} below.

\begin{proposition}[Wavepacket approximation at non-stationary points]\label{prop:wavepacket}
Let $J$ be a finite union of pairwise-disjoint closed subintervals of $(0,M)$, and denote by $J_\delta$ the corresponding union of closed intervals each of which is one of the intervals of $J$ extended by $\delta$ at both ends, such that the intervals of $J_\delta$ are also pairwise-disjoint and contained in $(0,M)$.  Suppose that an amplitude function $a:J_\delta\to\mathbb{R}_{>0}$ is of class $C^\infty(J_\delta)$ and strictly bounded away from zero, and that a phase function $S:J_\delta\to\mathbb{R}$ is of class $C^\infty(J_\delta)$ with derivative $S'$ strictly bounded away from $2\pi\mathbb{Z}$.  Let $\chi$ be a $C^\infty(0,M)$ cutoff function for which $\chi(y)=1$ for $y\in J$ and $\chi(y)=0$ for $y\in (0,M)\setminus J_\delta$, assume that $a(\cdot)$ and $S(\cdot)$ are functions independent of $\epsilon$, and define a wavepacket ansatz $\mathbf{u}$ with components $u_j$ given by 
\begin{equation}
u_{j}=\chi(y_j)a(y_j)\ee^{\ii S(y_j)/\epsilon},\quad y_j=\epsilon (j-1/2).
\label{eq:wavepacket-ansatz}
\end{equation}
Then $\mathbf{u}$ is an approximate eigenvector of $\mathbf{A}(x,t)$ with eigenvalue $\alpha\in\mathbb{R}$ in the sense that for each $y_j\in J$,
\begin{equation}
    [(\mathbf{A}(x,t)-\alpha\mathbb{I})\mathbf{u}]_j=(\zeta(y_j;\epsilon) + \mathcal{O}(\epsilon^2))u_j
\label{eq:residual}
\end{equation}
holds with $\zeta(y_j;\epsilon)$ being a uniformly bounded function of $y_j\in J$ that has zero mean with gridscale-wavelength (i.e., $\mathcal{O}(\epsilon)$ in $y_j$) oscillations wherever $S''(y_j)\neq 0$, provided that 
\begin{itemize}
    \item 
the phase derivative $S'(y)$ satisfies the eikonal equation
\begin{equation}\label{eq:eikonal}
-2\Lambda(y)\left(x+2\Lambda(y)t+\gamma(\Lambda(y))-\frac{U(S'(y))}{\Lambda'(y)}\right)=\alpha,\quad y\in J,
\end{equation}
where $U$ is the piecewise-linear function defined on $(-\pi,\pi)$ whose graph is shown in Figure~\ref{fig:Uplot};
\begin{figure}
\begin{center}
   \includegraphics[width= 0.4\linewidth]{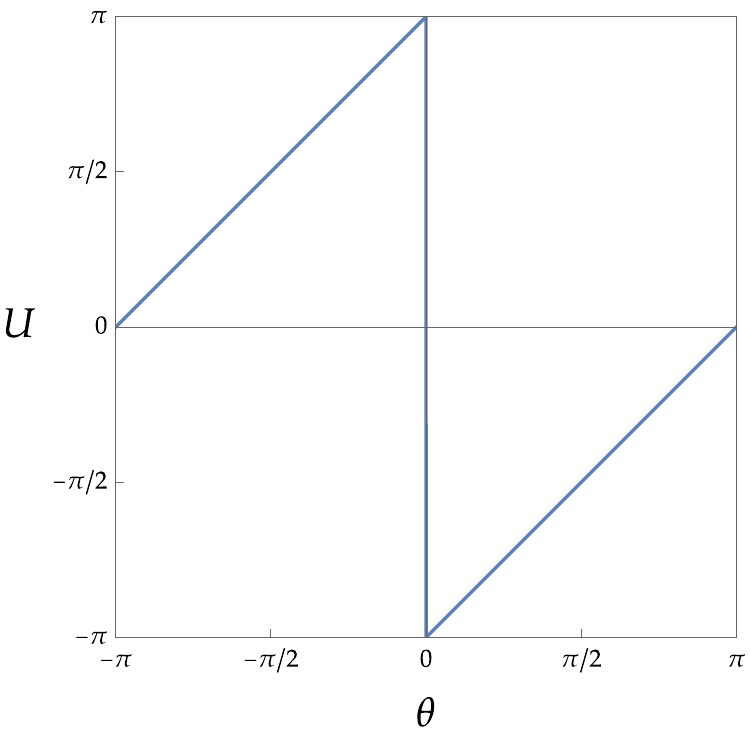} 
\end{center}
\caption{The graph of $U(\theta)$ on $(-\pi,\pi)$.}
\label{fig:Uplot}
\end{figure}
\item the amplitude $a(y)$ is subject to the equation:  
\begin{equation}
\frac{\dd}{\dd y}\frac{\Lambda(y)}{\Lambda'(y)}+\frac{2\Lambda(y)}{\Lambda'(y)}\frac{a'(y)}{a(y)}=0,\quad y\in J. \label{eq:amplitude-ODE}
\end{equation}
\end{itemize}
Note that unlike \eqref{eq:eikonal}, the amplitude equation \eqref{eq:amplitude-ODE} is independent of the eigenvalue $\alpha$.
\end{proposition}

\begin{remark}
The constraint that $U(\theta)\in [-\pi,\pi]$ gives a range of admissible values for $y$ for which the eikonal equation \eqref{eq:eikonal} can be solved for $S'(y)$, and the intervals of $J$ should consist of admissible values only.  When $\alpha=0$, the endpoints of intervals of admissible $y$ can be identified with the branches of the possibly multi-valued solution of Burgers' equation, see Corollary~\ref{cor:eikonal} below.  For the coordinates $(x,t)$ selected for the plots in Figure~\ref{fig:eigenvectors of A}, there are three branches of the solution for the indicated initial condition and the corresponding interval endpoints are indicated on the plots with dotted vertical lines.  For such admissible $y$, we can use the identity $U(U(\theta))=\theta$ to solve explicitly for $S'(y)$ when $\alpha=0$:
\begin{equation}
S'(y)=U\left(\Lambda'(y)[x+2\Lambda(y)t +\gamma(\Lambda(y))]\right).
\label{eq:Sprime}
\end{equation}
From \eqref{eq:Sprime} we can see that if $u_0$ is an analytic initial condition so that $\Lambda(y)$ is an analytic function of $y\in (0,M)$, then there are at most finitely-many points $y\in (0,M)$ for which $S''(y)=0$ and near which the function $\zeta(\cdot;\epsilon)$ fails to be rapidly oscillatory. 
For the wavepacket ansatz \eqref{eq:wavepacket-ansatz}, the meaning of $S'(y)$ is that it should be the approximate value of the site-to-site phase shift, as can be seen by Taylor-expanding the phase $S(y_j)$ about $y_{j-1}$, using $y_j=y_{j-1}+\epsilon$.  In Figure~\ref{fig:eigenvectorphase} we illustrate the remarkable accuracy of the prediction of the formula \eqref{eq:Sprime} for this phase shift.
\begin{figure}
\begin{center}
    \includegraphics[width=0.45\linewidth]{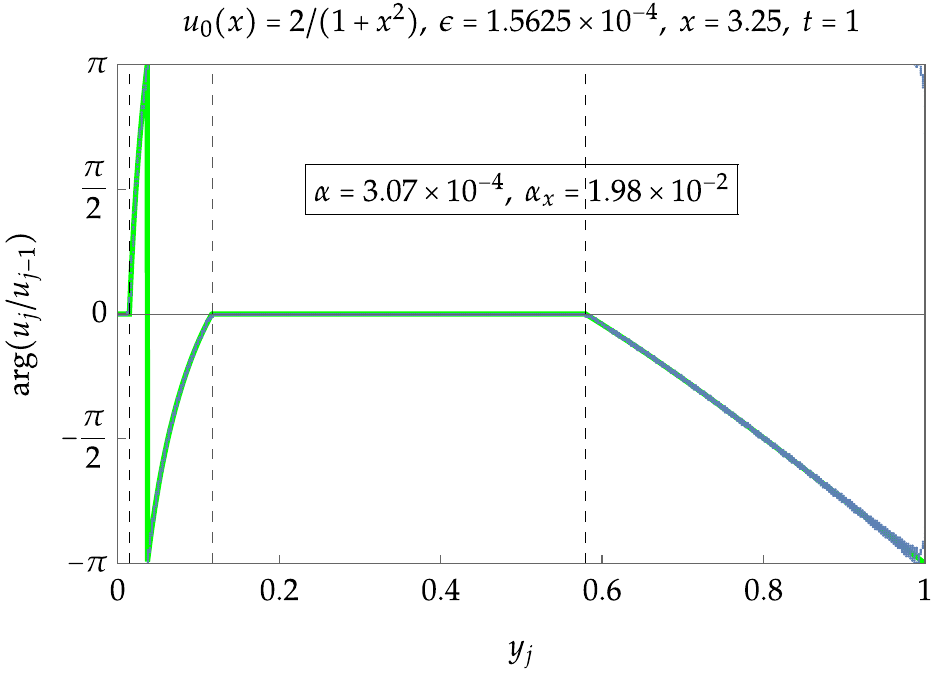}\hfill%
    \includegraphics[width=0.45\linewidth]{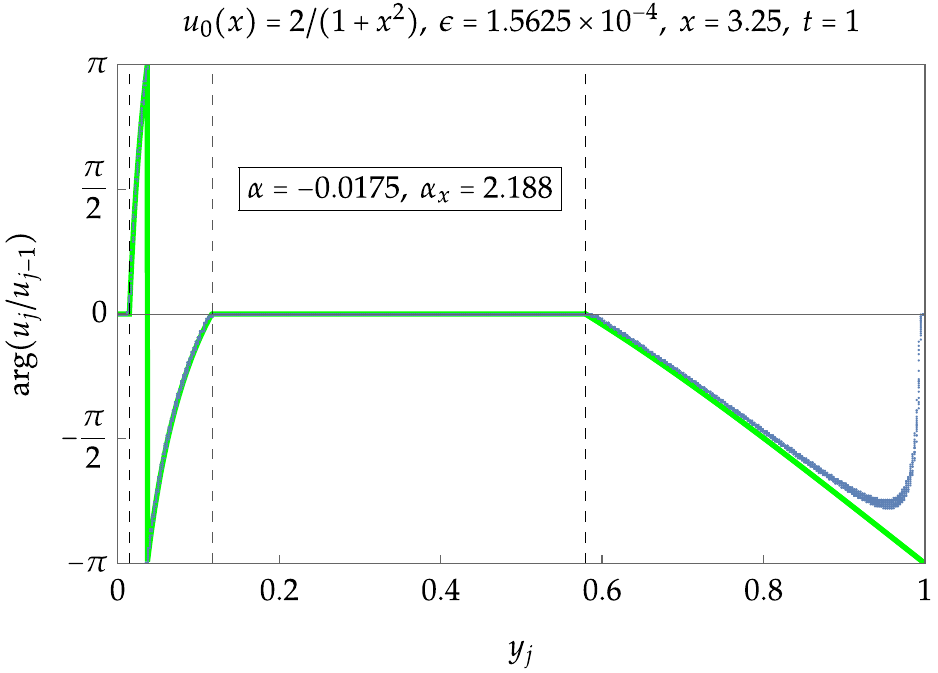}
\end{center}
\caption{Same as Figure~\ref{fig:eigenvectors of A}, except plotting the relative phase of nearest neighbor eigenvector elements in blue, and comparing with $S'(y_k)=U((x+2\Lambda(y)t+\gamma(\Lambda(y)))/F(\Lambda(y)))$ in green as determined from the eikonal equation for $\alpha=0$.  Note that for the eigenvector for a ``slow'' eigenvalue (left panel), the phase shift predicted by the eikonal equation is even accurate in the interval of admissible $y_j$ that does not abut $y=M$, where the amplitude is very small.  On the other hand, for the ``fast'' eigenvalue (right panel), some accuracy of the phase prediction is evidently lost near $y=M$, possibly because the eigenvalue $\alpha$ is not sufficiently small.}
\label{fig:eigenvectorphase}
\end{figure}
\end{remark}

\begin{remark}
Similarly, by explicit integration of the differential equation \eqref{eq:amplitude-ODE}, the amplitude $a(y)>0$ can be obtained for $y$ in the range of admissible values in the form
\begin{equation}\label{eq:amplitude}
a(y)=a_0\sqrt{\frac{\Lambda'(y)}{-\Lambda(y)}}.
\end{equation}
If there are multiple pairwise disjoint intervals of admissible values of $y\in [0,M]$, then the value of the integration constant $a_0>0$ may be different for each.  The form \eqref{eq:amplitude} can be fit by a least-squares computation to determine the value of $a_0$ in each interval of admissible $y_j$ from given eigenvector data; such fits are shown with green curves in Figure~\ref{fig:eigenvectors of A}; they are also remarkably accurate on the support subinterval of $[0,M]$ that abuts $y=M$.  The approximation appears to be slightly less accurate on the support subinterval that is separated from $y=M$ as seen in the right-hand panel of Figure~\ref{fig:eigenvectors of A}, because there are noticeable gridscale-wavelength fluctuations about the mean predicted by \eqref{eq:amplitude}.  We may expect that these terms might be captured by a refinement of the ansatz \eqref{eq:wavepacket-ansatz} to include a highly-oscillatory correction to the amplitude which could perhaps be chosen to removing the oscillatory correction $\zeta(y_j;\epsilon)$ from the residual in \eqref{eq:residual}.
\end{remark}

Proposition~\ref{prop:wavepacket} essentially shows that $\mathbf{A}(x,t)$ acts microlocally on suitable vectors.  Let us now turn to its proof. The following lemma will be useful.
\begin{lemma}[Nonstationary phase sums]
Suppose that $h:\mathbb{R}\to\mathbb{C}$ is of class $C^\infty(\mathbb{R})$ and has compact support $\mathrm{supp}(h)=[a,b]$, and suppose that $S:\mathrm{supp}(h)\to\mathbb{R}$ is of class $C^\infty(\mathbb{R})$ with $S'(y)$ bounded away from $2\pi\mathbb{Z}$ on $\mathrm{supp}(h)$.  Then recalling the notation $y=y_k=\epsilon(k-\frac{1}{2})$, 
\begin{equation}
\sum_{k\in\mathbb{Z}}
h(y_k)\ee^{\ii S(y_k)/\epsilon} = \mathcal{O}(\epsilon^\infty),\quad\epsilon\to 0.
\label{eq:cancellation}
\end{equation}
\label{lem:remainder}
\end{lemma}

\begin{proof}
By the Poisson summation formula,
\begin{equation}
\begin{split}
\sum_{k\in\mathbb{Z}}h(y_k)\ee^{\ii S(y_k)/\epsilon} &= \sum_{m\in\mathbb{Z}}\int_\mathbb{R}\ee^{-2\pi\ii km}h(y_k)\ee^{\ii S(y_k)/\epsilon}\,\dd k\\
&=\frac{1}{\epsilon}\sum_{m\in\mathbb{Z}}(-1)^m\int_\mathbb{R}h(y)\ee^{\ii (S(y)-2\pi m y)/\epsilon}\,\dd y.
\end{split}
\end{equation}
Integrating by parts $n>1$ times, we get
\begin{equation}
\int_\mathbb{R}h(y)\ee^{\ii (S(y)-2\pi my)/\epsilon}\,\dd y = \epsilon^n\int_\mathbb{R} h_n(y;m)\ee^{\ii (S(y)-2\pi my)/\epsilon}\,\dd y,
\end{equation}
wherein the function $h_n(y;m)$ is defined recursively by $h_0(y;m):=h(y)$ and 
\begin{equation}
h_n(y;m):=\ii\frac{\dd}{\dd y}\left[\frac{h_{n-1}(y;m)}{S'(y)-2\pi m}\right],\quad n\ge 1.
\end{equation}
It is easy to see that for fixed $n>1$, $h_n(y;m)$ satisfies an estimate of the form
\begin{equation}
\sup_{y\in\mathbb{R}}|h_n(y;m)|\le \frac{K_n}{\langle m\rangle^n},\quad \langle m\rangle:=\sqrt{1+m^2}.
\end{equation}
Therefore since $\mathrm{supp}(h)=[a,b]$ implies also $\mathrm{supp}(h_n)\subset [a,b]$,
\begin{equation}
\left|
\sum_{k\in\mathbb{Z}}
h(y_k)\ee^{\ii S(y_k)/\epsilon}\right|\le \epsilon^{n-1}\sum_{m\in\mathbb{Z}}
\int_a^b\frac{K_n}{\langle m\rangle^n}\,\dd y =K_n(b-a)\epsilon^{n-1}\sum_{m\in\mathbb{Z}}\langle m\rangle^{-n}= \mathcal{O}(\epsilon^{n-1}),\quad\epsilon\to 0
\end{equation}
holds for each $n>1$.  This completes the proof.
\end{proof}

\begin{proof}[Proof of Proposition~\ref{prop:wavepacket}]
Again we write $y=y_j$ for a lattice point in $J$ and set $\Delta y =\Delta y_{jk}$. We start by separating the action of $\mathbf{A}(x,t)-\alpha\mathbb{I}$ on the wavepacket vector $\mathbf{u}$ defined in \eqref{eq:wavepacket-ansatz} according to the diagonal and off-diagonal elements of $\mathbf{A}(x,t)$:
\begin{multline}
\begin{aligned}
\relax[(\mathbf{A}(x,t)-\alpha\mathbb{I})\mathbf{u}]_j
	&=\sum_{k=1}^{N(\epsilon)}(A_{jk}(x,t)-\delta_{jk}\alpha)\chi(y+\Delta y) a(y+\Delta y)\ee^{\ii S(y+\Delta y)/\epsilon} \\
	&=(A_{jj}(x,t)-\alpha)a(y)\ee^{\ii S(y)/\epsilon} \\
&\qquad\qquad\qquad{}+ \sum_{\substack{k=1\\ k\neq j}}^{N(\epsilon)}A_{jk}(x,t)\chi(y+\Delta y)a(y+\Delta y)\ee^{\ii S(y+\Delta y)/\epsilon}.
\end{aligned}
\end{multline}
Using~\eqref{eq:BOD-exact} along with $\chi(y)=1$, this shows that 
\begin{multline}
[(\mathbf{A}(x,t)\mathbf{u}-\alpha\mathbb{I})\mathbf{u}]_j 
	=(A_{jj}(x,t)-\alpha-\epsilon f(y,0))a(y)\ee^{\ii S(y)/\epsilon}\\
	{}+\frac{2\ii\Lambda(y) }{ \Lambda'(y)}\sum_{\substack{k=1\\k\neq j}}^{N(\epsilon)}
\frac{\chi(y+\Delta y)a(y+\Delta y)\ee^{\ii S(y+\Delta y)/\epsilon}}{k-j}\\
	{}+\sum_{k=1}^{N(\epsilon)}\epsilon f(y,\Delta y)\chi(y+\Delta y)a(y+\Delta y)\ee^{\ii S(y+\Delta y)/\epsilon}.
\end{multline}
We use \eqref{eq:f(y,0)} and \eqref{eq:A-diagonal-elements} to express $f(y,0)$ and the diagonal elements of $\mathbf{A}(x,t)$ respectively.  Then, extending the sum on the third line to $k\in\mathbb{Z}$ using compact support of $\chi$ and applying Lemma~\ref{lem:remainder} with $h(y_k)=f(y_j,y_k-y_j)\chi(y_k)a(y_k)$ individually on each of the intervals of $J_\delta$, we get
\begin{multline}
[(\mathbf{A}(x,t)-\alpha\mathbb{I})\mathbf{u}]_j 		
	= \left(-2\Lambda(y)\left(x+2\Lambda(y)t+\gamma(\Lambda(y))\right)-\alpha -\ii\epsilon\frac{\dd}{\dd y}\frac{\Lambda(y)}{\Lambda'(y)} \right)a(y)\ee^{\ii S(y)/\epsilon}\\
	+\frac{2\ii\Lambda(y)}{\Lambda'(y)}\sum_{\substack{k=1\\k\neq j}}^{N(\epsilon)}\frac{\chi(y+\Delta y)a(y+\Delta y)\ee^{\ii S(y+\Delta y)/\epsilon}}{k-j} + \mathcal{O}(\epsilon^\infty).
 \label{eq:Auj-last}
\end{multline}
Therefore, to very high accuracy, $\mathbf{A}(x,t)$ acts on the wavepacket $\mathbf{u}$ as multiplication by a sum of a diagonal matrix and the product of a diagonal matrix and an exact Toeplitz matrix with elements $T_{jk}=(k-j)^{-1}$ for $k\neq j$ and $T_{jj}=0$.

To deal with the Toeplitz part, we again use compact support of $\chi$ to write, for $K>0$ sufficiently large,
\begin{equation}
\sum_{\substack{k=1\\k\neq j}}^{N(\epsilon)}\frac{\chi(y+\Delta y)a(y+\Delta y)\ee^{\ii S(y+\Delta y)/\epsilon}}{k-j}
	= \sum_{\substack{k=j-K\\k\neq j}}^{j+K}\frac{\chi(y+\Delta y)a(y+\Delta y)\ee^{\ii S(y+\Delta y)/\epsilon}}{k-j}.
 \label{eq:Sum-extend}
\end{equation}
To allow $y=y_j$ to range over the full set $J$ and have exact equality in \eqref{eq:Sum-extend} for $K$ independent of $j$, we will assume that $K=N(\epsilon)$. We introduce the notation
\begin{equation}
S_2(\Delta y;y):=S(y+\Delta y)-S(y)-S'(y)\Delta y
\label{eq:S2}
\end{equation}
which vanishes to second order as $\Delta y\to 0$. We also use the Taylor expansion
\begin{equation}
    \chi(y+\Delta y)a(y+\Delta y)\ee^{\ii S_2(\Delta y;y)/\epsilon} = a(y) + a'(y)\Delta y + \frac{1}{2}R_2(\Delta y;y,\epsilon)\Delta y^2,
\label{eq:Taylor-product}
\end{equation}
where the real and imaginary parts of $R_2(\Delta y;y,\epsilon)$ are those of the second derivative with respect to $\Delta y$ of the left-hand side evaluated at two generally different points between $0$ and $\Delta y$. Hence we write the sum in \eqref{eq:Sum-extend} in the form
\begin{multline}
\sum_{\substack{k=j-K\\k\neq j}}^{j+K}\frac{\chi(y+\Delta y)a(y+\Delta y)\ee^{\ii S(y+\Delta y)/\epsilon}}{k-j} \\
{}
\begin{aligned}
&=\ee^{\ii S(y)/\epsilon}\sum_{\substack{k=j-K\\k\neq j}}^{j+K}\chi(y+\Delta y)a(y+\Delta y)\ee^{\ii S_2(\Delta y)/\epsilon}\frac{\ee^{\ii (k-j)S'(y)}}{k-j}\\
&= 
a(y)\ee^{\ii S(y)/\epsilon}\sum_{\substack{n=-K\\n\neq 0}}^{K}\left[
\frac{\ee^{\ii nS'(y)}}{n} + \epsilon \frac{a'(y)}{a(y)}\ee^{\ii nS'(y)} + \frac{\epsilon^2}{2a(y)} R_2(\epsilon n;y,\epsilon)n\ee^{\ii nS'(y)}\right],
\end{aligned}
\label{eq:Taylor-in-sum}
\end{multline}
where we used $\Delta y=\epsilon (k-j)$ and reindexed by $n=k-j$ in the last line. Since $y=y_j$ is independent of the index $n$, we just need to examine three summands, $\ee^{\ii n\theta}/n$, $\ee^{\ii n\theta}$, and $R_2(\epsilon n;\epsilon)n\ee^{\ii n\theta}$ for $\theta=S'(y)$, which by assumption is bounded away from $2\pi\mathbb{Z}$. 
According to \cite[Eqn.\@ 4.36]{MillerX11}:
\begin{equation}
\lim_{K\to\infty}\sum_{\substack{n=-K\\n\neq 0}}^K\frac{\ee^{\ii n\theta}}{n} = -\ii U(\theta),\quad \theta\neq 0\pmod{2\pi},\label{eq:Usum}
\end{equation}
where the function $U(\theta)$ is periodically extended with period $2\pi$.
More precisely, by representing the summand as an integral and exchanging the order of finite summation and integration, we can write
\begin{equation}
\begin{split}
    \sum_{\substack{n=-K\\n\neq 0}}^K\frac{\ee^{\ii n\theta}}{n} &= \ii\int_{\pm\pi}^\theta\left[\cos(K\tau)-1+\sin(K\tau)\cot(\tfrac{1}{2}\tau)\right]\,\dd \tau\\ &=\ii(\pm\pi-\theta)+\frac{\ii\sin(K\theta)}{K} + \ii\int_{\pm\pi}^\theta\sin(K\tau)\cot(\tfrac{1}{2}\tau)\,\dd\tau.
    \end{split}
\end{equation}
If we assume that the sign on $\pm\pi$ corresponds to the sign of $\theta\in [-\pi,\pi]\setminus\{0\}$, then $\ii(\pm\pi-\theta)=-\ii U(\theta)$ and the remaining integral admits repeated integration by parts.  In this way we obtain (using also $K=N(\epsilon)=\epsilon^{-1}M+\mathcal{O}(1)$)
\begin{equation}
\begin{split}
    \sum_{\substack{n=-K\\n\neq 0}}^K\frac{\ee^{\ii n\theta}}{n} &= -\ii U(\theta)-\frac{\ii}{K}\cdot \frac{\cos((K+\frac{1}{2})\theta)}{\sin(\frac{1}{2}\theta)} + \mathcal{O}(K^{-2}),\quad K\to\infty,\\
    &=-U(\theta)-\frac{\ii\epsilon}{M}\cdot\frac{\cos((N(\epsilon)+\frac{1}{2})\theta)}{\sin(\frac{1}{2}\theta)} + \mathcal{O}(\epsilon^2),\quad\epsilon\to 0,
    \end{split}
    \label{eq:sum-1}
\end{equation}
which holds uniformly for $\theta\in [-\pi,\pi]$ bounded away from zero.  Note that the singularity of the correction term at $\theta=0$, which is also the jump point for $U(\theta)$, is related to Gibbs' phenomenon\footnote{This is the main difficulty in extending the theory of Toeplitz quantization to the setting of non-smooth symbols such as $U(\theta)$.  See Section~\ref{sec:Toeplitz} below for more information.}.  A more straightforward calculation gives that the sum of $\ee^{\ii n\theta}$ is exactly a constant shift of the Dirichlet kernel:
\begin{equation}
    \sum_{\substack{n=-K\\n\neq 0}}^K\ee^{\ii n\theta} = -1+\sum_{n=-K}^K\ee^{\ii n\theta}=-1+\frac{\sin((K+\frac{1}{2})\theta)}{\sin(\frac{1}{2}\theta)}=-1+\frac{\sin((N(\epsilon)+\frac{1}{2})\theta)}{\sin(\frac{1}{2}\theta)}.
    \label{eq:sum-2}
\end{equation}

Upon evaluation for $\theta=S'(y)=S'(y_j)$ bounded away from $2\pi\mathbb{Z}$, the terms $\cos((N(\epsilon)+\frac{1}{2})\theta)/\sin(\frac{1}{2}\theta)$ and $\sin((N(\epsilon)+\frac{1}{2})\theta)/\sin(\frac{1}{2}\theta)$ will be highly oscillatory zero-mean functions of $y=y_j$ when $\epsilon$ is large, at each point $y$ with $S''(y)\neq 0$.  We introduce the notation $\widetilde{\mathcal{O}}(\epsilon^p)$ to denote such a function, whose absolute value is also $\mathcal{O}(\epsilon^p)$ in the usual sense.  With this notation, \eqref{eq:sum-1} and \eqref{eq:sum-2} respectively imply that
\begin{equation}
    \sum_{\substack{n=-K\\n\neq 0}}^K\frac{\ee^{\ii n S'(y)}}{n} = -U(S'(y))+\widetilde{\mathcal{O}}(\epsilon) + \mathcal{O}(\epsilon^2),\quad\epsilon\to 0,
    \label{eq:oscillatory-1}
\end{equation}
\begin{equation}
    \epsilon\sum_{\substack{n=-K\\n\neq 0}}^K\ee^{\ii n S'(y)} = -\epsilon+\widetilde{\mathcal{O}}(\epsilon),\quad\epsilon\to 0.
    \label{eq:oscillatory-2}
\end{equation}
We omit the details of the estimation of the contribution of sum involving the Taylor error term $R_2(\epsilon n;\epsilon)$.  However we claim that 
\begin{equation}
    \sum_{n=-K}^KR_2(\epsilon n;\epsilon)n\ee^{\ii n S'(y)} = \widetilde{\mathcal{O}}\left(\frac{1}{\epsilon^2}\right),\quad\epsilon\to 0.
    \label{eq:big-claim}
\end{equation}
Indeed, keeping one more term in the Taylor expansion \eqref{eq:Taylor-product}, the first term in $R_2(\Delta y;\epsilon)$ would just be the second derivative of $a(y+\Delta y)\ee^{\ii S_2(\Delta y;y)/\epsilon}$ at $\Delta y=0$, which one can check is of the form $\epsilon^{-1}v(y)+w(y)$ for some smooth and bounded functions $v$ and $w$.  The contribution of this term to the sum in \eqref{eq:big-claim} is then $-\ii\epsilon^{-1}v(y)-\ii w(y)$ times the derivative (from the factor of $\ii n$ in the summand) of the Dirichlet kernel $\sin((N(\epsilon)+\frac{1}{2})\theta)/\sin(\frac{1}{2}\theta)$ evaluated at $\theta=S'(y)$, so one can easily check that the contribution is a term of the form $\widetilde{\mathcal{O}}(\epsilon^{-2})$.  Similar arguments apply to the terms obtained by continuing the expansion \eqref{eq:Taylor-product} to any finite order.  In fact, the $k^\mathrm{th}$ derivative of the left-hand side of \eqref{eq:Taylor-product} evaluated at $\Delta y=0$ is $\epsilon^{-\lfloor k/2\rfloor}$ times a polynomial in $\epsilon$ with coefficients that are smooth bounded functions of $y$.  This shows that one should not replace $R_2(\Delta y;\epsilon)$ with too many explicit terms, since although their contributions will be oscillatory functions of $y$, they will also start to grow in size.

Combining \eqref{eq:Auj-last}, \eqref{eq:Sum-extend}, \eqref{eq:Taylor-in-sum}, and the estimates \eqref{eq:oscillatory-1}--\eqref{eq:big-claim} yields
\begin{multline}
 [(\mathbf{A}(x,t)-\alpha\mathbb{I})\mathbf{u}]_j
    =\left[-2\Lambda(y)\left(x+2\Lambda(y)t+\gamma(\Lambda(y))-\frac{U(S'(y))}{\Lambda'(y)}\right)-\alpha\right.\\
   \left. -\ii\epsilon\left(\frac{\dd}{\dd y}\frac{\Lambda(y)}{\Lambda'(y)}+\frac{2\Lambda(y)}{\Lambda'(y)}\frac{a'(y)}{a(y)}\right) + \widetilde{\mathcal{O}}(1)+\mathcal{O}(\epsilon^2)\right]u_j.
\end{multline}
Neglecting the highly oscillatory term and keeping only the mean yields \eqref{eq:eikonal} at the leading order in $\epsilon$ and \eqref{eq:amplitude-ODE} at second order.
\end{proof}

\begin{corollary}[Admissible values of $y$ and Burgers' equation]\label{cor:eikonal}
Let $\alpha=0$.  The admissible values of $y\in [0,M]$ for which there exists $S'(y)\in(-\pi,\pi)$ such that~\eqref{eq:eikonal} holds are determined by:
\begin{itemize}
    \item $\Lambda(y)\in [-u_0^\mathrm{B}(x,t),0]$ if $(x,t)$ is a point in the single-valued region for the solution of the inviscid Burgers equation, and
    \item $\Lambda(y)\in[-u_0^\mathrm{B}(x,t),0]\cup[-u_2^\mathrm{B}(x,t),-u_1^\mathrm{B}(x,t)]$ if $(x,t)$ is a point in the triple-valued region for the solution of the inviscid Burgers equation.
\end{itemize}
\end{corollary}

\begin{proof}
From the eikonal equation \eqref{eq:eikonal} with $\alpha=0$, the range condition $U(S'(y))\in(-\pi,\pi)$ means that $\lambda=\Lambda(y)$ and $x$ have to satisfy the condition
\begin{equation}\label{eq:eikonal_admissible}
\frac{x+2\lambda t+\gamma(\lambda)}{F(\lambda)}\in(-\pi,\pi),
\end{equation}
where we used the identity $\Lambda'(y)=1/F(\Lambda(y))$.

It turns out that the curves $x+2\lambda t+\gamma(\lambda)=\pm\pi{F(\lambda)}$ in the $(\lambda,x)$-plane essentially give the (rotated and reflected) graph of the multi-valued solution $u^\mathrm{B}(\cdot,t)$ to the Burgers equation obtained by the method of characteristics. Indeed, at time $t=0$, assume that
\begin{equation}
\frac{x+\gamma(\lambda)}{F(\lambda)}=\pi.
\end{equation}
From the definitions \eqref{eq:density of eigenvalues} and \eqref{eq:gamma-formula} of $F$ and $\gamma$, respectively, this means that
\begin{equation}
x=\pm\pi F(\lambda)-\gamma(\lambda)=x_\pm(\lambda),
\end{equation}
where $x_-(\lambda)<x_+(\lambda)$ are the turning points satisfying $u_0(x_\pm(\lambda))=-\lambda$,  see below~\eqref{eq:LM}.  
Therefore $(x,-\lambda)$ belongs to the graph of $u_0$.

Similarly, at time $t$, the method of characteristics implies that $(x_{\pm}(\lambda)-2\lambda t,-\lambda)$ belongs to the graph of $u^\mathrm{B}(\cdot,t)$. As a consequence, one sees that $(x,-\lambda)$ is in the graph of $u^\mathrm{B}(\cdot,t)$ if and only if $x=x_{\pm}(\lambda)-2\lambda t$, which is equivalent to
\[
\frac{x+2\lambda t+\gamma(\lambda)}{F(\lambda)}=\pm\pi.
\]

Finally, we observe that given $x$ such that $u^\mathrm{B}(x,t)$ has only one branch $u^\mathrm{B}_0(x,t)$, the range of admissible $(x,-\lambda)$ such that~\eqref{eq:eikonal_admissible} holds is the region between the graph of $u^\mathrm{B}_0(\cdot,t)$ and the $x$-axis. This means that the admissible values of $y$ correspond to $\lambda\in [-u^\mathrm{B}_0(x,t),0]$.
However, when $x$ is such that there are three branches $u^\mathrm{B}_0(x,t)<u^\mathrm{B}_1(x,t)<u^\mathrm{B}_2(x,t)$, the range of admissible $(x,-\lambda)$ such that~\eqref{eq:eikonal_admissible} holds is the union of the region between the graph of $u^\mathrm{B}_0(\cdot,t)$ and the $x$-axis and also the region between the graphs of $u^\mathrm{B}_1(\cdot,t)$ and $u^\mathrm{B}_2(\cdot,t)$.  Hence the admissible values of $y$ correspond to $\lambda\in [-u_2^\mathrm{B}(x,t),-u_1^\mathrm{B}(x,t)]\cup [-u_0^\mathrm{B}(x,t),0]$.
\end{proof}

\subsection{Toeplitz quantization and small eigenvalues of \texorpdfstring{$\mathbf{A}(x,t)$}{A}}
\label{sec:Toeplitz}
Proposition~\ref{prop:wavepacket} suggests that the matrix $\mathbf{A}(x,t)$ is represented by a symbol $p_x$ depending on canonical variables $(y,\theta)\in[0,M]\times[-\pi,\pi]$ and given by 
\begin{equation}\label{eq:symbol}
p_x(y,\theta):=-2\Lambda(y)\left(x+2\Lambda(y)t+\gamma(\Lambda(y))-\frac{U(\theta)}{\Lambda'(y)}\right)
\end{equation}
such that the eikonal equation~\eqref{eq:eikonal} takes the form $p_x(y,\theta)=\alpha$.  In the following, we mainly consider small eigenvalues satisfying $|\alpha|\leq \epsilon^r$, in view of \eqref{eq:large-A-eigenvalues}. In the case of the zero eigenvalue $\alpha=0$, the equation $p_x(y,\theta)=\alpha$ can be written as
\begin{equation}
p(y,\theta)=x,
\label{eq:level-curve-on-sphere}
\end{equation}
where a modified symbol is defined by
\begin{equation}\label{eq:symbol-simple}
p(y,\theta):=-2\Lambda(y)t-\gamma(\Lambda(y))+\frac{U(\theta)}{\Lambda'(y)}.
\end{equation}

Unfortunately, we are not aware of results on Toeplitz quantization on rectangles, making the study of $p_x$ and $p$ on $[0,M]\times[-\pi,\pi]$ difficult.
In order to gain insight on the eigenvector approximation, one may look instead into the results on Toeplitz quantization on the sphere by using latitude/longitude coordinates $(y,\theta)$.  The strategy of using Toeplitz quantization on the sphere to study dispersionless PDEs was rigorously implemented in the context of the dispersionless Toda system in~\cite{Bloch2003dispersionless}.  However, this approach is not rigorous in our situation because of the discontinuity of $U$ at the angle $\theta=0$ and the lack of smoothness of the symbol at the poles $y=0,M$ of the sphere.

\begin{figure}
\begin{center}
    \includegraphics[width=0.5\linewidth]{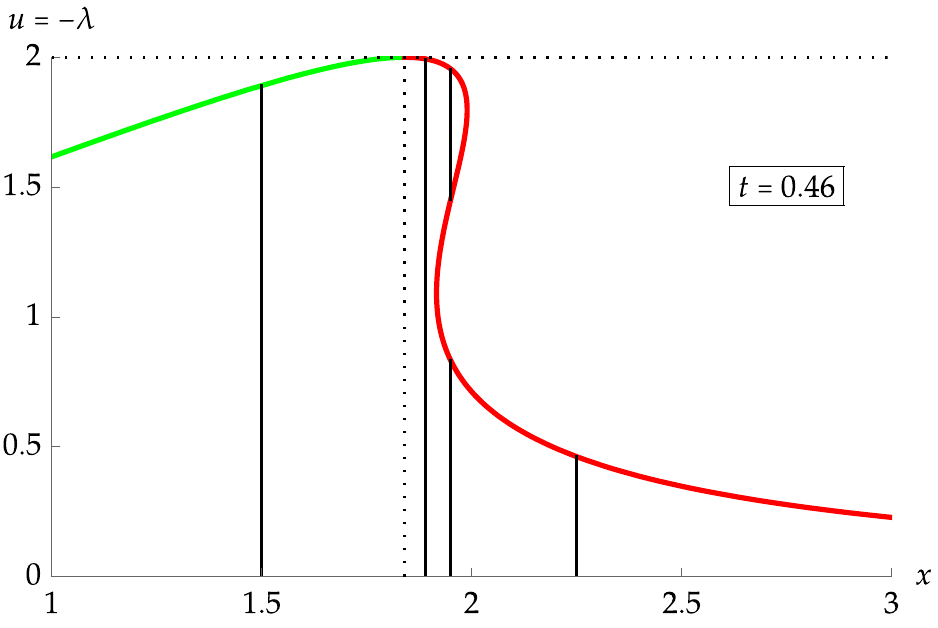}\\
    \includegraphics[width=0.3\linewidth]{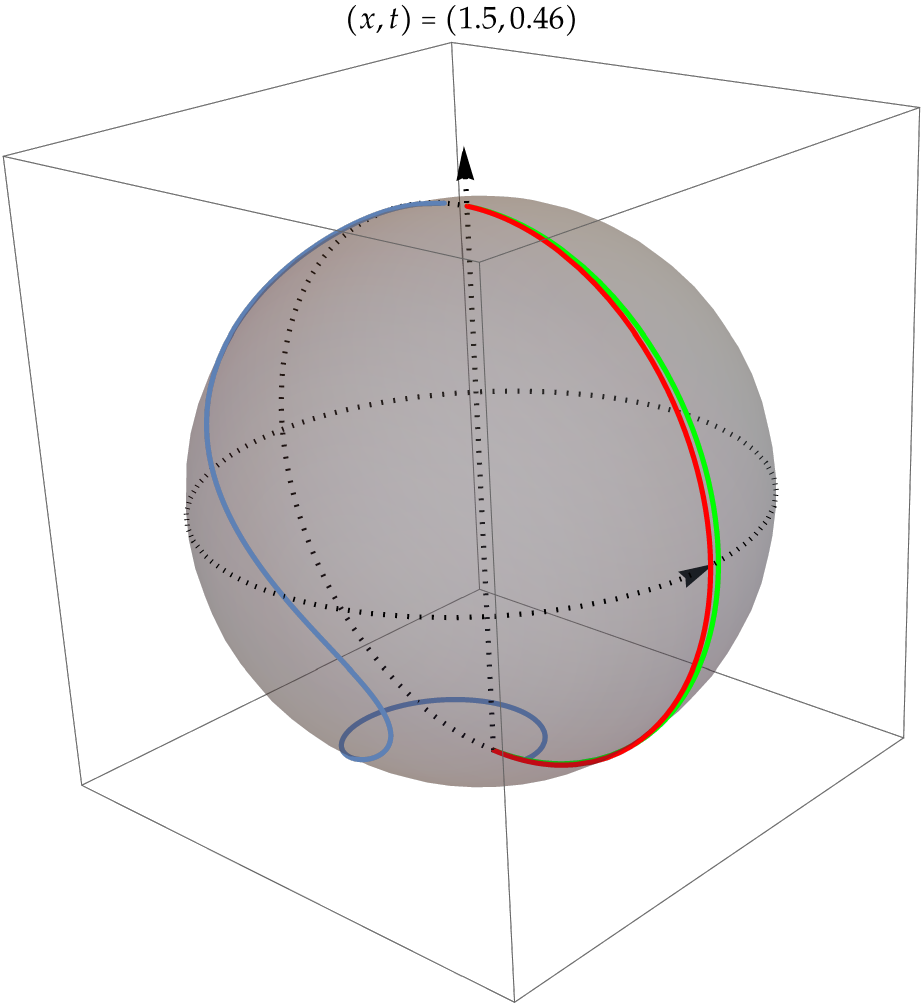}\hspace{0.2\linewidth}\includegraphics[width=0.25\linewidth]{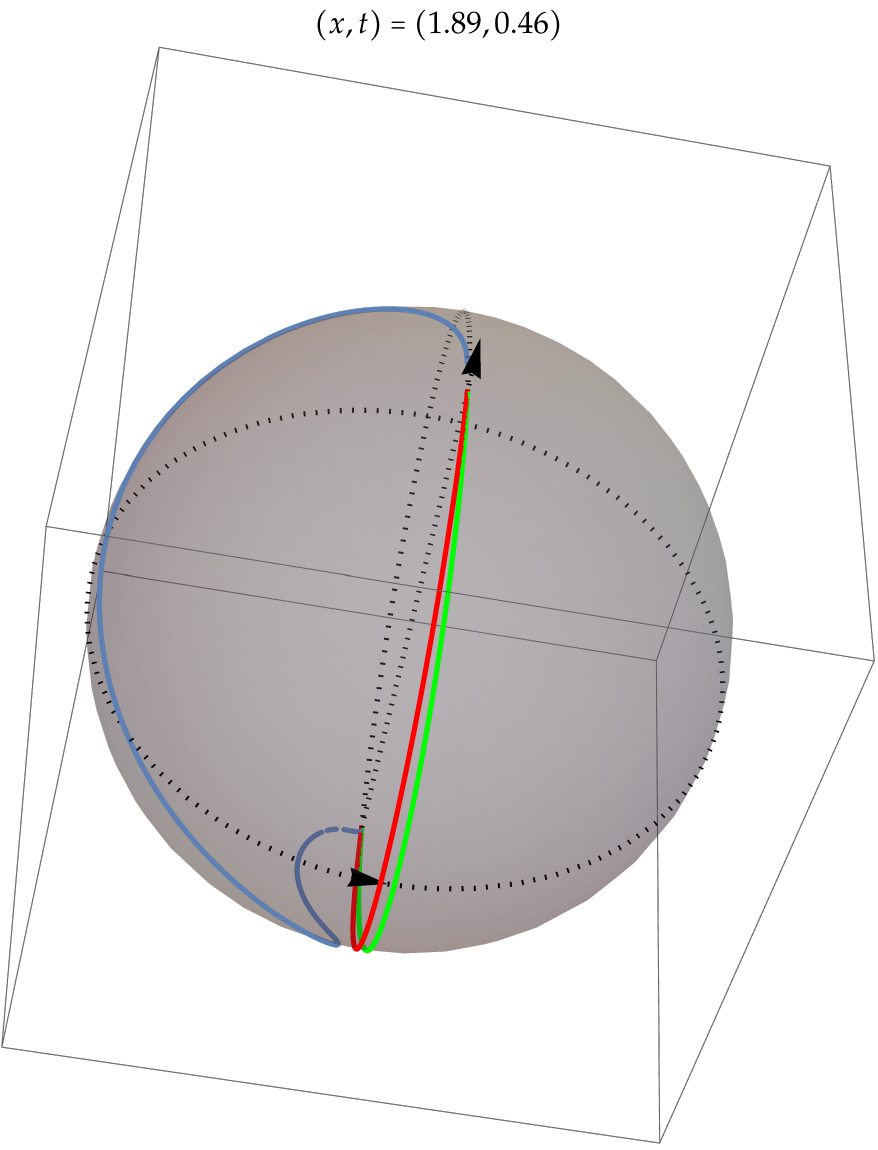}\\
    \includegraphics[width=0.3\linewidth]{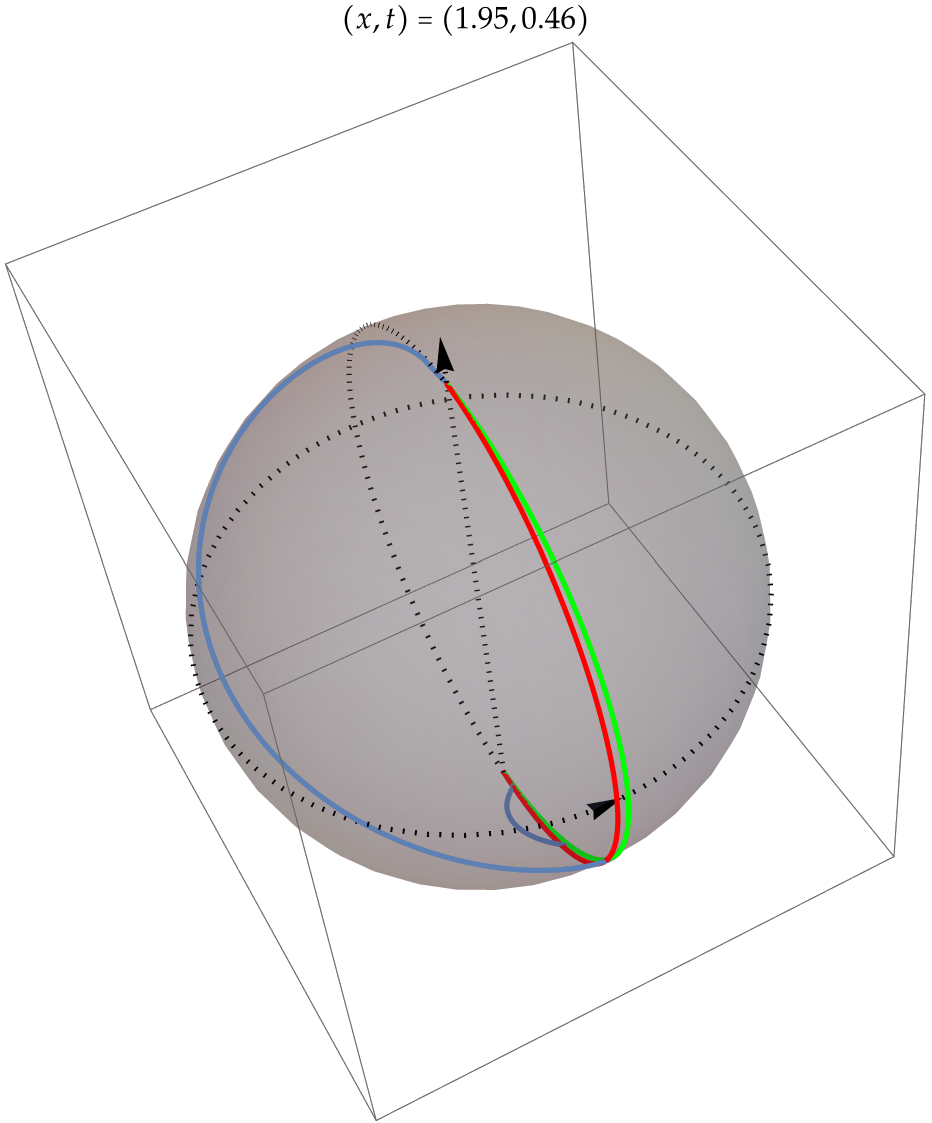}\hspace{0.2\linewidth}\includegraphics[width=0.3\linewidth]{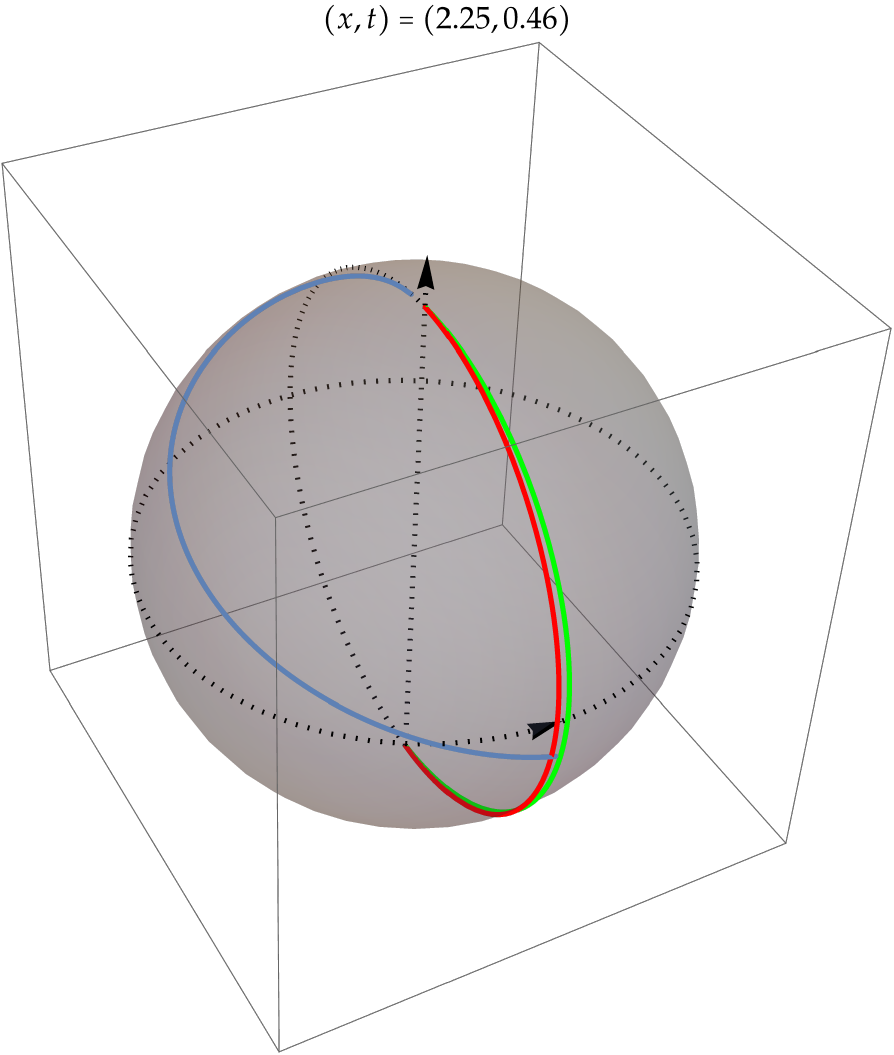}
\end{center}
\caption{Level sets $p(y,\theta)=x$ in the $(x,u)$-plane and on the sphere with vertical polar axis $y\in (0,M)$ and meridian $0<\theta<2\pi$ (green for $\theta=0$, red for $\theta=2\pi$) for $t=0.46$ and various $x$ with $u_0(x)=2/(1+x^2)$. Note that when $u=0$, $y=M$ (north pole) and $\theta=\pi$.  As $u$ increases along a vertical line in the top plot, $y$ decreases along one or two components of an orbit on the corresponding sphere plot below.}
\label{fig:levelsets1}
\end{figure}

\begin{figure}
\begin{center}
    \includegraphics[width=0.5\linewidth]{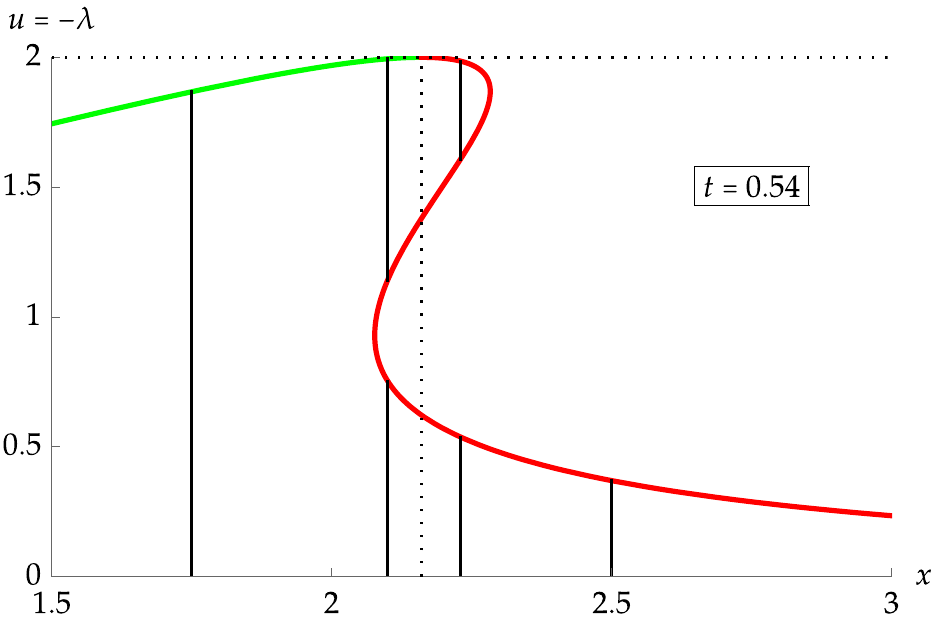}\\
    \includegraphics[width=0.3\linewidth]{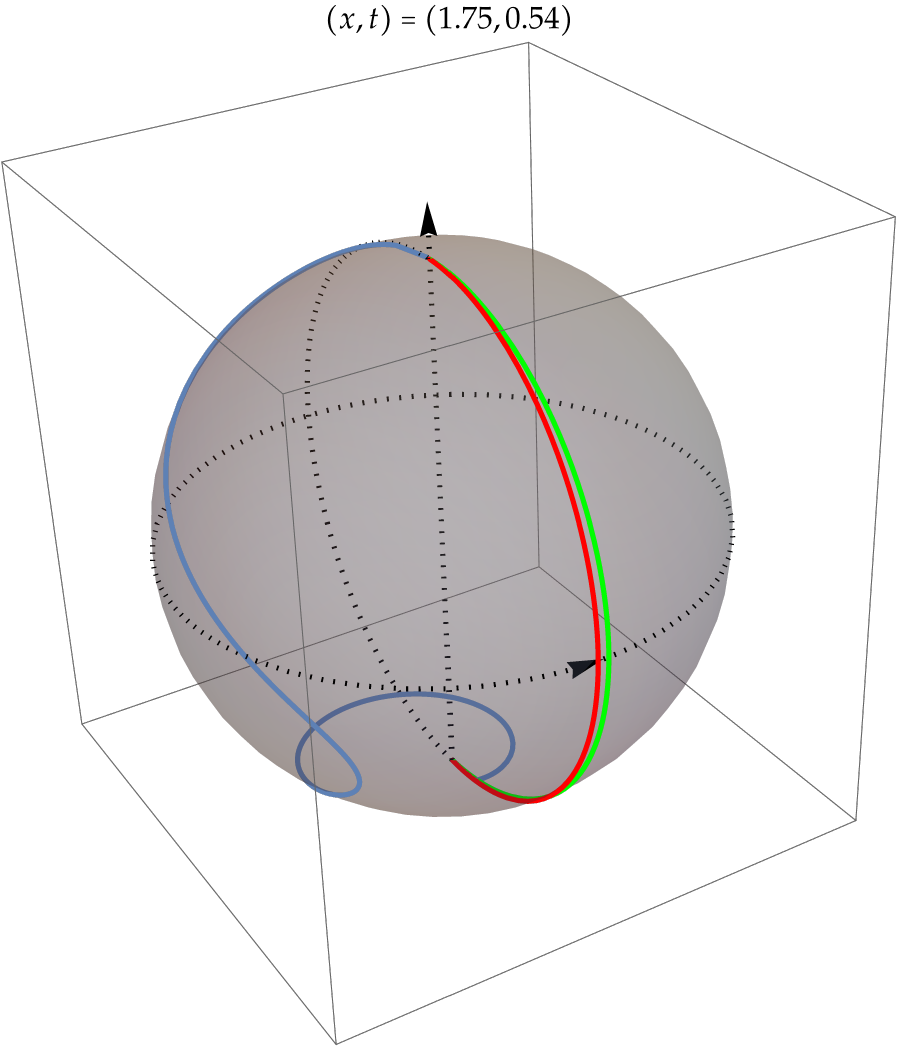}\hspace{0.2\linewidth}\includegraphics[width=0.3\linewidth]{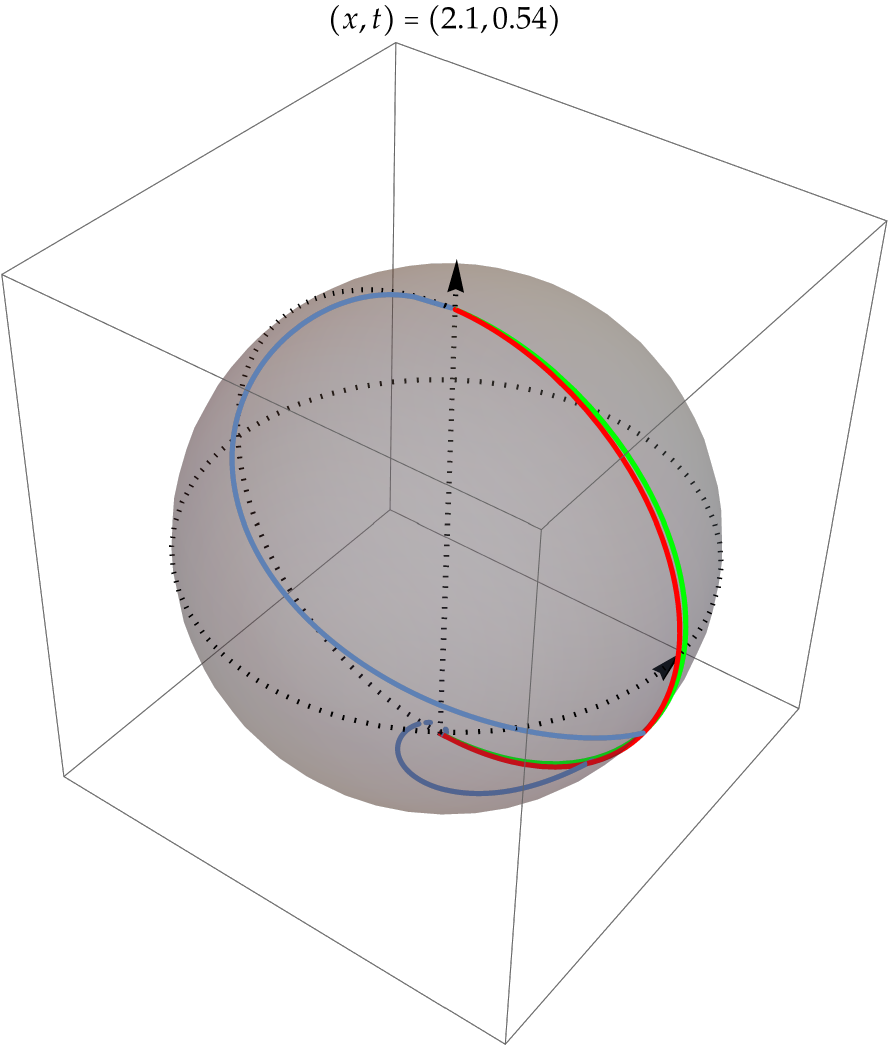}\\
    \includegraphics[width=0.3\linewidth]{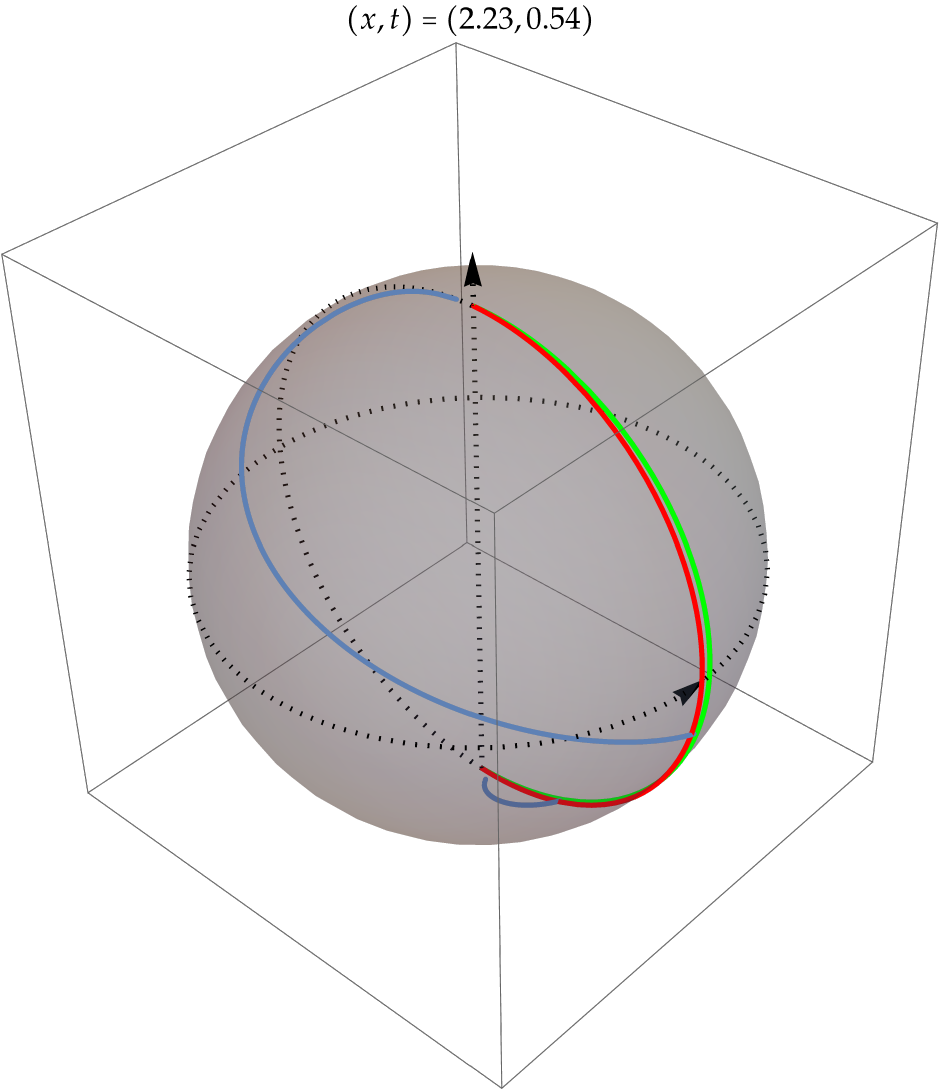}\hspace{0.2\linewidth}\includegraphics[width=0.3\linewidth]{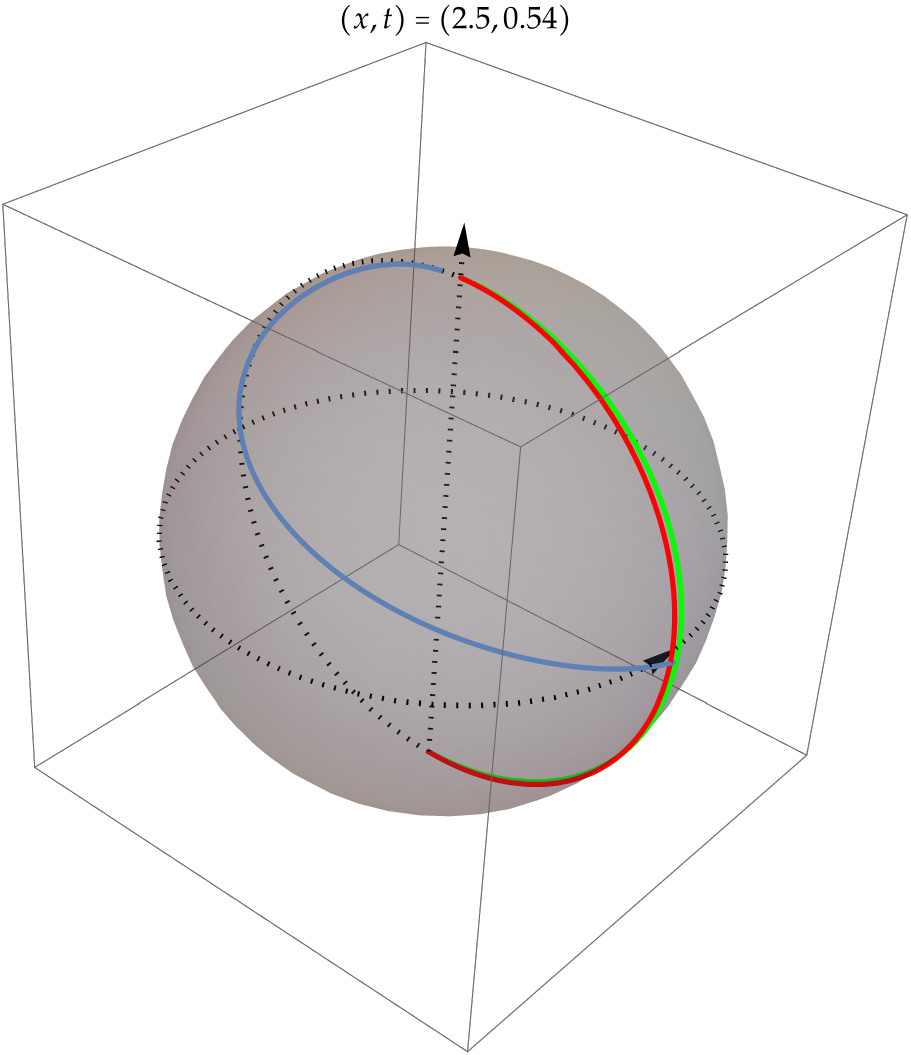}
\end{center}
\caption{Level sets $p(y,\theta)=x$ for $t=0.54$ and various $x$ with $u_0(x)=2/(1+x^2)$.}
\label{fig:levelsets2}
\end{figure}

According to Corollary~\ref{cor:eikonal}, the ranges of admissible values of $y$ for \eqref{eq:level-curve-on-sphere} to hold are determined by the branches of the Burgers solution at the point $(x,t)$.  The equation \eqref{eq:level-curve-on-sphere} for a given $(x,t)$ describes a relation between the latitude coordinate $y$ and the longitude coordinate $\theta$ that yields curves generally beginning and ending on the meridian $\theta=0$ or $\theta=2\pi$
except for one curve that always emerges from the north pole $y=M$ with longitude $\theta=\pi$.  We call different connected components of the curve given by \eqref{eq:level-curve-on-sphere} \emph{orbits} on the spherical phase space.  The connection between vertical slices through the graph of the multi-valued solution of Burgers' equation at different values of $x$ and the orbits on the sphere is illustrated in Figures~\ref{fig:levelsets1} and \ref{fig:levelsets2}.

In the context of a smooth symbol on the sphere, the paper~\cite{Charles2003quasimodes}  (see also~\cite{Charles2006symbolic}) states that the WKB expansion for the eigenvector $\mathbf{u}$ for a very small eigenvalue $\alpha$ should be valid with amplitude being a discretization of~\eqref{eq:amplitude}
on the admissible range of $y$, and such a small eigenvalue can only occur for $x$ satisfying a certain Bohr-Sommerfeld condition. More precisely, we decompose the admissible range of $y$ into one interval if $(x,t)$ is in the single-valued region, and into two disjoint intervals if $(x,t)$ is in the triple-valued region for the inviscid Burgers solution. Then for each interval, there should be an eigenvector with support exactly in that interval, and amplitude given by~\eqref{eq:amplitude}. Moreover, the values of $x$ for which $\alpha=0$ is an eigenvalue follow a Bohr-Sommerfeld condition: at first order, the parallel transport around the orbit with latitude coordinate $y$ confined to a given admissible interval is quantized, see~\cite[Theorem 0.1]{Charles2003quasimodes}.

Numerically, however, we observe that the support of the ‘‘fast'' eigenvector is not localized only in the interval $\Lambda(y)\in[-u^\mathrm{B}_2(x,t),-u^\mathrm{B}_1(x,t)]$ (the main part), as there is some smaller contribution from the other interval $\Lambda(y)\in[-u^\mathrm{B}_0(x,t),0]$ (the correction), see
Figure~\ref{fig:eigenvectors of A}. This may be due to the fact that the symbol is not smooth in our case. Moreover, we observe that this small contribution has a noticeable effect in the calculation. We will choose different values of $a_0$ in order for both the main part and the correction to have a non-negligible contribution to the Euclidean norm of the eigenvector. More precisely,  in Conjecture~\ref{conjecture:Bohr} below, we will assume that the correction has an amplitude of order $a_0\sim1/\ln(\epsilon^{-1})$ compared to the main part for which $a_0\sim 1$.  Unlike the scale $\delta(\epsilon)$ in Section 2 which was somewhat flexible (see Conjecture \ref{conjecture: nu upper} and the preceding paragraph), here the choice of precise scale inversely proportional to $\ln(\epsilon^{-1})$ is important for the proof of Theorem 3.7 below.

\subsubsection{Bohr-Sommerfeld conditions}

We guess  the Bohr-Sommerfeld conditions from the case of a smooth symbol on the sphere~\cite{Charles2003quasimodes} by using the parallel transport in the range of admissible $y$ and $\theta(\Lambda(y))=S'(y)$ given by~\eqref{eq:Sprime}, which we compute using the Levi-Civita connection on the sphere, relative to the metric inherited from the Euclidean one on $\mathbb{R}^3$. At first order, this translates as follows in the case of small eigenvalues $\alpha\approx 0$, using the symbol $p$ given in~\eqref{eq:symbol-simple}.
We make the change of variable $y=Y(\lambda)$ inverse to $\lambda=\Lambda(y)$ in \eqref{def: big lambda}, so that
\begin{equation}
Y(\lambda)=\int_{-L}^{\lambda}F(\eta)\,\dd\eta.
\label{eq:Y-of-lambda}
\end{equation}
Fixing $t>t_b$, we work in the vicinity of a point $x_0$ in the triple-valued region for the solution of Burgers' equation:  $X^-(t)<x_0<X^+(t)$. Let us define the parallel transport about an orbit corresponding to a fast zero eigenvalue by
\begin{equation}
g_0^+(x,t)=\int_{-u_2^\mathrm{B}(x,t)}^{-u_1^\mathrm{B}(x,t)}\theta'(\lambda) Y(\lambda)\, \dd\lambda,
\label{eq:g0plus}
\end{equation}
and that corresponding to a slow zero eigenvalue by
\begin{equation}
g_0^-(x,t)=\int_{-u_0^\mathrm{B}(x,t)}^{0} \theta'(\lambda) Y(\lambda)\, \dd\lambda,
\label{eq:g0minus}
\end{equation}
where $\theta(\lambda)$ is determined from the eikonal equation for $\alpha=0$ in the form $p(Y(\lambda),\theta(\lambda))=x$. 
The main idea of Toeplitz quantization is that small eigenvalues of a matrix $\mathbf{A}(x,t)$ associated with\footnote{The theory of Toeplitz quantization assumes that a sequence of matrices $\mathbf{A}(x,t)$ is constructed systematically from a given symbol, whereas in the present context we start instead with a sequence of matrices $\mathbf{A}(x,t)$ and deduce the symbol from the WKB expansion in Proposition~\ref{prop:wavepacket}.} a symbol and a corresponding orbit can only occur if, to leading order, the parallel transport around that orbit is an integer multiple of $2\pi\epsilon$, which is known as a \emph{Bohr-Sommerfeld quantization rule}.  For smooth symbols, the leading terms $g_0^\pm(x,t)$ are modified by a sum of higher-order corrections $\epsilon^ng_n^\pm(x,t)$, $n=1,2,3,\dots$ for which the formul\ae\ written in \cite{Charles2003quasimodes} are not expected to be meaningful in the non-smooth setting.

Since $g_0^\pm(x,t)$ should be quantized but also $x\approx x_0$, we define values of $x=x_k^\pm(t)\approx x_0$ corresponding to the two types of orbits by Bohr-Sommerfeld conditions of the form
\begin{equation}
g_0^{\pm}(x_k^{\pm}(t),t)+\epsilon g_1^{\pm}(x_k^\pm(t),t)
    =2\pi (k+k_0^\pm)\epsilon,\quad |k|\le C\varepsilon^{r-1},\quad k\in\mathbb{Z},
    \label{eq:BS-0}
\end{equation}
with $k_0^\pm:= \lfloor (g_0^\pm(x_0,t)+\epsilon g_1^\pm(x_0,t))/(2\pi\epsilon)\rfloor $.

Assuming moreover that $g_1^\pm$ depends smoothly in $x$, we  replace $g_1^{\pm}(x_k^\pm,t)$ by $g_1^{\pm}(x_0,t)$ up to higher-order corrections that we neglect here: with $\rho^\pm(x_0,t):=k_0^\pm-g_1^\pm(x_0,t)/(2\pi)$, the Bohr-Sommerfeld conditions become
\begin{equation}
g_0^{\pm}(x_k^{\pm}(t),t)
    =2\pi (k+\rho^\pm(x_0,t))\epsilon,\quad |k|\le C\varepsilon^{r-1},\quad k\in\mathbb{Z}.
    \label{eq:BS}
\end{equation}

The values $x=x_k^+(t)$ (resp., $x=x_k^-(t)$) are expected to approximate the values of $x\approx x_0$ for which $\alpha=0$ is a fast (resp., slow) eigenvalue, see Conjecture~\ref{conjecture:Bohr} below.

\subsubsection{Simplifying and differentiating the Bohr-Sommerfeld conditions}
We first simplify the expression
\begin{equation}
g_0^\pm(x,t)
=\int_{\lambda_{\min}}^{\lambda_{\max}}\theta'(\lambda)Y(\lambda)\dd\lambda.
\label{eq:g0-again}
\end{equation}
The integration limits are $\lambda_{\min}=-u_2^\mathrm{B}(x,t)$ and $\lambda_{\max}=-u_1^\mathrm{B}(x,t)$ in the case of fast eigenvalues (i.e. $g_0^\pm=g_0^+$), but $\lambda_{\min}=-u_0^\mathrm{B}(x,t)$ and $\lambda_{\max}=0$ in the case of slow eigenvalues (i.e. $g_0^\pm=g_0^-$).

Using that $U(\theta)=\theta-\pi$ on $(0,2\pi)$, we parametrize the relevant orbit by $y=Y(\lambda)$, $\theta=\theta(\lambda)\in(0,2\pi)$ and get that because $\Lambda'(y)=1/F(\Lambda(y))$ the equation $p(Y(\lambda),\theta(\lambda))=x$ in~\eqref{eq:symbol-simple} becomes
\begin{equation}\label{eq:theta}
x=F(\lambda)(\theta(\lambda)-\pi)-2\lambda t \implies F(\lambda)\theta(\lambda)=x+\pi F(\lambda)+2\lambda t.
\end{equation}
Substituting from \eqref{eq:Y-of-lambda} for $Y(\lambda)$ in \eqref{eq:g0-again}, integration by parts yields
\begin{equation}\label{eq:g0+-theta}
g_0^\pm(x,t)=\left[\theta(\lambda)\int_{-L}^{\lambda}F(\eta)\dd\eta\right]_{\lambda_{\min}}^{\lambda_{\max}} -\int_{\lambda_{\min}}^{\lambda_{\max}}\theta(\lambda)F(\lambda)\dd\lambda.
\end{equation}
Since $g_0^\pm$ is quantized, the positions $x_k^{\pm}(t)$ at which there is a fast ($+$) or slow ($-$) crossing of $\alpha=0$ are locally regularly spaced if the derivative of $g_0^\pm$ is bounded above and below. This derivative will determine the actual spacing. We note that $\theta(\lambda)$ evaluated on the limits of integration is locally independent of $(x,t)$, being equal to $0$ or $2\pi$, or in the case $\lambda_{\min}=0$, to $\pi$. Hence we compute
\begin{multline}
\partial_x g_0^\pm(x,t)
	=\theta(\lambda_{\max})F(\lambda_{\max})\partial_x\lambda_{\max} 
	-\theta(\lambda_{\min})F(\lambda_{\min})\partial_x \lambda_{\min}\\
	-\partial_x\lambda_{\max}	\theta(\lambda_{\max})F(\lambda_{\max})
	+\partial_x\lambda_{\min}
 \theta(\lambda_{\min})F(\lambda_{\min})-\int_{\lambda_{\min}}^{\lambda_{\max}}\partial_x(\theta(\lambda) F(\lambda))\dd\lambda.
\end{multline}
Since $\theta F$ is given by~\eqref{eq:theta}, we have that $\partial_x(\theta(\lambda) F(\lambda))=1$, hence
\begin{equation}
-\partial_x g_0^\pm(x,t)=\lambda_{\max}(x,t)-\lambda_{\min}(x,t)=
\begin{cases}
u_2^\mathrm{B}(x,t)-u_1^\mathrm{B}(x,t), &\text{ for fast eigenvalues,}\\
u_0^\mathrm{B}(x,t),&\text{ for slow eigenvalues}.
\end{cases}
\label{eq:g0prime-two-cases}
\end{equation}
Due to the assumption $|x_l^\pm-x_0|\leq C\epsilon^r$, one can replace $x$ by $x_0$ in the integrand of the Bohr-Sommerfeld condition \eqref{eq:BS} rewritten in the form
\begin{equation}
\int_{x_l^\pm}^{x_{l+1}^\pm}\partial_x g_0^\pm(x,t)\,\dd x=2\pi\varepsilon,
\end{equation}
up to an error term of order $\epsilon^r$ in the integrand.
In both cases, we deduce
\begin{equation}
(x_{l+1}^\pm-x_l^\pm)(1+o(1)) = \frac{-2\pi\varepsilon}{f^\pm(x_0,t)},
\end{equation}
where
\begin{equation}
    f^+(x,t):=u_2^\mathrm{B}(x,t)-u_1^\mathrm{B}(x,t)>0
\label{eq:f-fast}
\end{equation}
and
\begin{equation}
    f^-(x,t):=u_0^\mathrm{B}(x,t)>0.
\label{eq:f-slow}
\end{equation}
We conclude that 
\begin{equation}\label{eq:x_l+1-x_l}
x_{l+1}^\pm-x_l^\pm= \frac{-2\pi}\varepsilon{f^\pm(x_0,t)}+o(\varepsilon).
\end{equation}
Note that $\partial_x g_0^\pm(x,t)<0$, so the sequences $x_l^\pm$ are decreasing.
More precisely, with the same argument but integrating $g_0^\pm$ from $x_0$ to $x_l^\pm$, we can write for some $r^\pm(x_0,t)$:
\begin{equation}\label{eq:xl}
x_l^\pm=-\frac{2\pi l \epsilon}{ f^\pm(x_0,t)}+ r^\pm(x_0,t)\epsilon+o(\epsilon).
\end{equation}

\begin{remark}
In~\eqref{eq:g0+-theta}, one can make the change of variable $Y(\lambda)=y$ and use~\eqref{eq:Y-of-lambda}. Since $\theta(\Lambda(y))=S'(y)$, we get
\begin{equation}
\frac{1}{\epsilon}\int_{Y(\lambda_{\min})}^{Y(\lambda_{\max})} S'(y)\dd y
	=\frac{1}{\epsilon}\left[\theta(\lambda_{\max})Y(\lambda_{\max}) -\theta(\lambda_{\min})Y(\lambda_{\min})
	-g_0^{\pm}(x,t)\right].
\end{equation}
 Note that $\theta(\lambda_{\max})$ and $\theta(\lambda_{\min})$ are expected to belong to the set $\{0,\pi,2\pi\}$, and that $Y(\lambda_{\min})/\epsilon$ and $Y(\lambda_{\max})/\epsilon$ should be integers (eigenvector component indices). The Bohr-Sommerfeld conditions~\eqref{eq:BS} therefore suggest that between two consecutive eigenvalues, the accumulated phase shift of the corresponding eigenvector as measured by the left-hand side of the above identity changes by $2\pi$ at first order.
\end{remark}

\subsubsection{Classification of the small eigenvalues according to the Bohr-Sommerfeld conditions}

 Conjecture~\ref{conjecture:Bohr} below expresses the following idea. For $t>t_b$ and $x_0\in (X^-(t),X^+(t))$, we look at the small eigenvalues $|\alpha_k(x,t)|\leq \epsilon^r$ for $x\approx x_0$, where we recall that the parameter $0<r<\frac 12$ is chosen so that the large eigenvalues satisfy~\eqref{eq:large-A-eigenvalues}. If $\alpha_k$ is ‘‘slow'', then on Figure~\ref{fig:alphas}, one can follow the ‘‘slow'' line with small slope passing the point $(x,\alpha_k)$. Using the predictions of the Bohr-Sommerfeld conditions, this line crosses the axis $\alpha=0$ at one of the points~$x_l^-$: in other words  $(x_l^-,0)$ belongs to the same ‘‘slow'' line. Similarly, if $\alpha_k$ is ‘‘fast'', then on Figure~\ref{fig:alphas}, one can follow the ‘‘fast'' line with large slope passing the point $(x,\alpha_k)$: this line also passes one point of the form $(x_l^+,0)$. Then using the predictions from Proposition~\ref{prop:wavepacket} of the wavepacket approximation of the eigenvector when $\alpha=0$, the eigenvector $\mathbf{u}$ at $x=x_l^\pm$ with eigenvalue $\alpha=0$ should match the observations from Figure~\ref{fig:eigenvectors of A} and Corollary~\ref{cor:eikonal}.

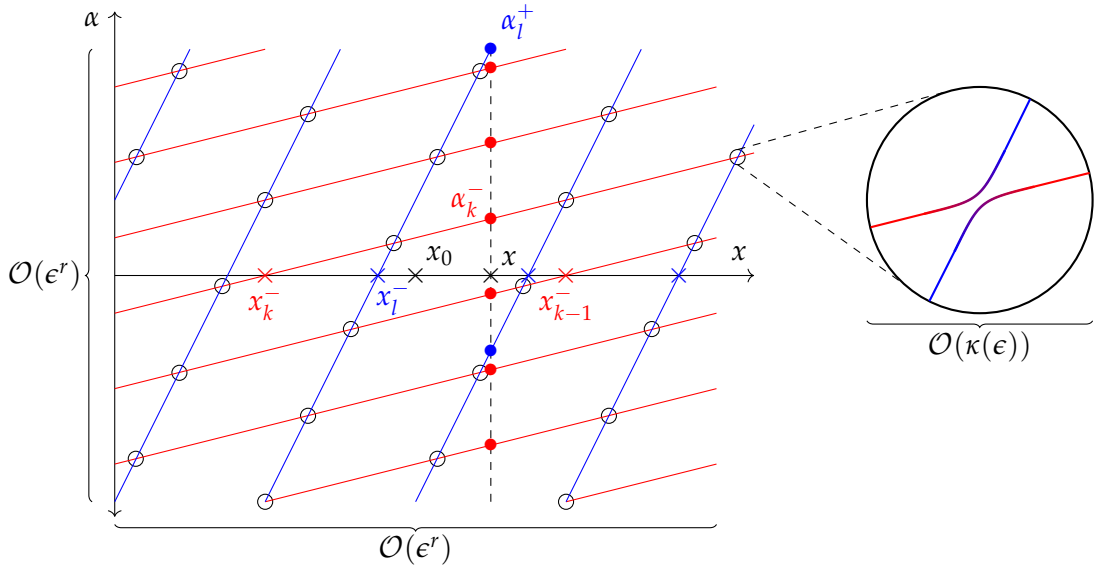
\begin{figure}
\begin{center}
\begin{tikzpicture}
% Zones à exclure 
%l=1
\draw (2,1.5) circle (1mm);
\draw (4.29,2.07) circle (1mm);
\draw (6.57,2.64) circle (1mm);
 
%l=2
\draw (0.29,2.07) circle (1mm);
\draw (2.57,2.64) circle (1mm);
\draw (4.86,3.21) circle (1mm);
 
%l=3
\draw  (0.86,3.21) circle (1mm);
% \draw  (3.14,3.79) circle (1mm);

%l=0
\draw (1.43,0.36) circle (1mm);
\draw (3.71,0.93) circle (1mm);
\draw (6,1.5) circle (1mm);
\draw (8.28,2.07) circle (1mm);

%l=-1
\draw  (0.86,-0.79) circle (1mm);
\draw  (3.14,-0.21) circle (1mm);
\draw  (5.43,0.36) circle (1mm);
\draw  (7.71,0.93) circle (1mm);

%l=-2
\draw  (0.28,-1.93) circle (1mm);
\draw  (2.57,-1.36) circle (1mm);
\draw  (4.86,-0.79) circle (1mm);
\draw  (7.14,-0.21) circle (1mm);

%l=-3
\draw (2,-2.5) circle (1mm);
\draw (4.29,-1.93) circle (1mm);
\draw (6.57,-1.36) circle (1mm);

%l=-4
\draw (6,-2.5) circle (1mm);

% Axe des abscisses
\draw[->] (0,0.5) -- (8.5,0.5);
\draw (8.3,0.8) node {$x$};
\draw (4,0.5) node {$\times$} node[above right] {$x_0$};
\draw (5,0.5) node {$\times$} node[above right] {$x$};
\draw[dashed] (5,-2.5)--(5,3.5);
% Axe des ordonnées
\draw[<->] (0,-2.7) -- (0,4);
\draw (-0.3,3.9) node {$\alpha$};

% Trait rouge
\draw[red] (0,3) -- (2,3.5);
\draw[red] (0,2) -- (6,3.5);
\draw[red] (0,1) -- (8,3);
\draw[red] (0,0) -- (8.5,2.125);
\draw[red] (0,-1) -- (8,1);
\draw[red] (0,-2) -- (8,0);
\draw[red] (2,-2.5) -- (8,-1);
\draw[red] (6,-2.5) -- (8,-2);

%Trait bleu
\draw[blue] (0,1.5) -- (1,3.5);
\draw[blue] (0,-2.5) -- (3,3.5);
\draw[blue] (2,-2.5) -- (5,3.5);
\draw[blue] (4,-2.5) -- (7,3.5);
\draw[blue] (6,-2.5) -- (8.5,2.5);

%Intersections rouges avec x
\draw[red] (5,-1.75) node {$\bullet$};
\draw[red] (5,-0.75) node {$\bullet$};
\draw[red] (5,0.25) node {$\bullet$};
\draw[red] (5,1.25) node {$\bullet$};
\draw[red] (4.7,1.5) node {$\alpha_k^-$};
\draw[red] (5,2.25) node {$\bullet$};
\draw[red] (5,3.25) node {$\bullet$};

%Intersections rouges avec abscisse
\draw[red] (2,0.5) node {$\times$} node[below] {$x_k^-$};
\draw[red] (6,0.5) node {$\times$} node[below] {$x_{k-1}^-$};

%intersections bleues avec x
\draw[blue] (5,3.5) node {$\bullet$} node[above right] {$\alpha_l^+$};
\draw[blue] (5,-0.5) node {$\bullet$};

%Intersections bleues avec abscisse
\draw[blue] (3.5,0.5) node {$\times$};
\draw[blue] (3.7,0.2) node {$x_l^-$};
\draw[blue] (5.5,0.5) node {$\times$};
\draw[blue] (7.5,0.5) node {$\times$};

% Double-flèche en-dessous avec la légende "O"
% \draw[<->,dotted] (0,-2.5) -- (8,-2.5) node[midway,below] {$\mathcal{O}(\epsilon^r)$};
\draw[decoration={brace},decorate] (8,-2.8) -- (0,-2.8) node[midway,below] {$\mathcal{O}(\epsilon^r)$};
\draw[decoration={brace},decorate] (-0.3,-2.5) -- (-0.3,3.5) node[midway,left] {$\mathcal{O}(\epsilon^r)$};

%zoom in on the right

%(8.28,2.07) circle (1mm)

\draw[dashed] (8.28,2.17) -- (11.1,2.95);
\draw[dashed] (8.28,1.97) -- (10.5,0.38);
\draw[thick] (11.5,1.5) circle (1.5cm);

\draw[red,thick] (10.05,1.14) -- (10.77,1.32);
\draw[red,thick] (12.23,1.68) -- (12.96,1.86);

\draw[blue,thick] (10.83,0.16) -- (11.17,0.83);
\draw[blue,thick] (11.84,2.17) -- (12.17,2.84);

\draw[test={0.5pt}{red}{blue}] (10.77,1.32) .. controls (11.5,1.5) .. (11.84,2.17);
\draw[test={0.5pt}{blue}{red}] (11.17,0.83) .. controls (11.5,1.5) .. (12.23,1.68);

\draw[decoration={brace},decorate] (13,-0.1) -- (10,-0.1) node[midway,below] {$\mathcal{O}(\kappa(\varepsilon))$};

\end{tikzpicture}
\end{center}
\caption{Illustration for Conjecture~\ref{conjecture:Bohr}. For $x-x_0=\mathcal{O}(\epsilon^r)$, the small eigenvalues of order $\epsilon^r$ are partitioned into two families $\alpha_k^-=\alpha_k^-(x,t)$ (``slow'' eigenvalues, in red) and $\alpha_l^+=\alpha_l^+(x,t)$ (``fast'' eigenvalues, in blue), that are modeled by parallel lines in the $(x,\alpha)$-plane of two distinct slopes. The small black circles in the figure are centered at intersections of slow and fast lines.  Eigenvalues within them cannot be reliably characterized as being slow or fast, and for them, the linear approximation is no longer expected to be accurate.  The set $X_o$ is the union of projections on the the $x$-axis of the horizontal diameters of those  circles centered on points with $\alpha=\mathcal{O}(\epsilon^s)$.}\label{fig:bohr}
\end{figure}

We illustrate the predicted position of small eigenvalues in Figure~\ref{fig:bohr} (which may be viewed as the theoretical version of the numerical plots in Figure~\ref{fig:alphas}). We model the ``slow'' and ``fast'' eigenvalues with two families of parallel lines in the $(x,\alpha)$-plane with different slopes.  True eigenvalues $\alpha_k(x,t)$ are unambiguously slow or fast when they follow a line from one or the other family; however an eigenvalue that is located near a crossing point of a slow and fast line has an ambiguous velocity $\alpha_{k,x}(x,t)$.  Note that any eigenvalues for which $\alpha_k(x,t)=\mathcal{O}(\epsilon)$ and for which the value of $\alpha_{k,x}(x,t)$ is ambiguous will make an unpredictable contribution to $u(x,t)$ via the formula \eqref{u in terms of alphas} of order $1$.  Therefore, to be able prove a convergence result without resolving the ambiguity, i.e., using linear approximations only, we should restrict the analysis to values of $x$ for which all such ambiguous eigenvalues in the small range $|\alpha_k(x,t)|\le\epsilon^r$ are sufficiently large compared to $\epsilon$.

With this motivation, we introduce a sufficiently small scale $\kappa=\kappa(\epsilon)$ and declare a small eigenvalue $\alpha(x,t)$ with $|\alpha(x,t)|\le\epsilon^r$ to be \emph{ambiguous} if it lies within a disk of radius $\kappa$ of a crossing point of the two families of straight lines in the $(x,\alpha)$-plane near $(x_0,0)$ given by the leading terms in \eqref{eq:approx-alpha-fast} below.  Then we pick an exponent $s$ with $0<r<s<1$ and define an excluded set $X_o=X_o(x_0,t;\epsilon)$ consisting of values of $x\in\mathbb{R}$ near $x_0$ such that $x\in X_o$ if there is a crossing point of the straight lines at $(x',\alpha')$ with $|x-x'|\le\kappa$ and $|\alpha'|\le\epsilon^s$.  

For fixed $x\approx x_0$ not in $X_o$, we denote by $(\alpha_k^o(x,t))_{|k|\leq C_o\epsilon^{r-1}}$ the ambiguous eigenvalues: for each such $\alpha_k^o(x,t)$, we have $|\alpha_k^o(x,t)|\gtrsim \epsilon^q$. Given that there are $\mathcal{O}(\epsilon^{r-1})$ slow lines intersecting the vertical axis passing the point $(x,0)$, we roughly estimate the number of such ambiguous crossings to be of order at most $\epsilon^{r-1}$, which explains the index range $|k|\leq C_o\epsilon^{r-1}$.  In what follows, we fix parameters $\frac 13<r<\frac 12$ and $\frac 12<q<1$ to be chosen later.

\begin{conjecture}[Bohr-Sommerfeld conditions for small eigenvalues]\label{conjecture:Bohr}
Let $t>t_b$ and fix $x_0\in (X^-(t),X^+(t))$, and let $\epsilon>0$ be sufficiently small. Let  $0<r<q<1$.

For $x$ such that $|x-x_0|\leq \epsilon^{r}$ 
and $x\not\in X_o$, small eigenvalues $\alpha$ with
\begin{equation}
|\alpha(x,t)|\leq \epsilon^r
\end{equation}
can be split into three pairwise disjoint families:
\begin{itemize}
\item a set of ``slow'' eigenvalues $(\alpha_k^-(x,t))_{|k|\leq C_-\epsilon^{r-1}}$;
\item a set of ``fast'' eigenvalues $(\alpha_k^+(x,t))_{|k|\leq C_+\epsilon^{r-1}}$;
\item a set of ``ambiguous'' eigenvalues $(\alpha_k^o(x,t))_{|k|\leq C_o\epsilon^{r-1}}$, with $|\alpha_k^o(x,t)|\gtrsim \epsilon^s$.
\end{itemize}
Moreover, for $x\not \in X_o$, $\alpha_k^\pm(x,t)$ is close to $\alpha_k^\pm(x_k^\pm(t),t)$, in the sense that there is a uniform Taylor expansion in $x-x_k^\pm(t)$:
\begin{equation}\label{eq:taylor-alpha}
\alpha_k^\pm(x,t)
	=\alpha_k^\pm(x_k^\pm(t),t)
	+\partial_x \alpha_k^\pm(x_k^\pm(t),t)(x-x_k^\pm(t))
	+o(\epsilon^{q}),
\end{equation}
\begin{equation}\label{eq:taylor-alpha_x}
\partial_x\alpha_k^\pm(x,t)
	=\partial_x \alpha_k^\pm(x_k^\pm(t),t)
	+o(\epsilon^q).
\end{equation}

The points $x_k^\pm$ satisfy the Bohr-Sommerfeld conditions in the following sense.
\begin{enumerate}
\item There holds $|x_k^\pm(t)-x_0|=\mathcal{O}(\epsilon^r)$, and  $\alpha_k^{\pm}(x_k^\pm(t),t)=\mathcal{O}(\epsilon^{1+r})$ is a small eigenvalue.
\item (Eigenvectors for slow eigenvalues.)
 The amplitude of the normalized eigenvector  for the eigenvalue $\alpha^-_k(x_k^-,t)$
 \begin{equation}
 \mathbf{u}={\bf u}^-_k(x_k^-,t)
\end{equation}
 is a discretization of~\eqref{eq:amplitude} restricted to the set $\Lambda(y)\in[-u^{\mathrm{B}}_0(x_k^-,t),0]$, up to a small remainder term.

\item (Eigenvectors for fast eigenvalues.)
We can decompose the eigenvector  for the eigenvalue $\alpha^+_k(x_k^+,t)$ as the linear combination of two normalized vectors 
\begin{equation}
{\bf u}
    ={\bf u}_k^+(x_k^+,t)
    +\frac{c}{\ln(\epsilon^{-1})} {\bf u}_l^-(x_k^+,t)
    +o\left(\frac{1}{\ln(\epsilon^{-1})}\right).
\end{equation}
The constant $c$ does not depend on $k$, it is bounded independently of $\epsilon,x_0,t$ and: 
\begin{itemize} 
\item (Fast main part.) The component ${\bf u}_k^+(x_k^+,t)$ has amplitude that is a discretization of~\eqref{eq:amplitude}  restricted to the set $\Lambda(y)\in[-u^\mathrm{B}_2(x_k^+,t),-u^\mathrm{B}_1(x_k^+,t)]$;
\item (Slow correction.) The component ${\bf u}_l^-(x_k^+,t)$ is a discretization of~\eqref{eq:amplitude} restricted the set $\Lambda(y)\in[-u^\mathrm{B}_0(x_k^+,t),0]$.
\end{itemize}
\end{enumerate}
\end{conjecture}

Note that  in the case of ``fast''  eigenvectors (as shown in the right-hand panel of Figure~\ref{fig:eigenvectors of A}), the ratio between the parameter $a_0^+$ tuned for the fast main part and the parameter $a_0^-$ tuned for the slow correction should satisfy
\begin{equation}
\frac{a_0^-}{a_0^+}=\frac{c+o(1)}{\ln(\epsilon^{-1})}.
\end{equation}

\begin{proposition}[Approximation of $\alpha_k$ and its spatial derivative]\label{prop:bohr}
Assume there exists a $p>\frac 12$ such that $u_0(x)\sim Cx^{-2p}$ as $|x|\to+\infty$, and that Conjecture~\ref{conjecture:Bohr} holds. Fix $t>t_b$ and $x_0\in (X^-(t),X^+(t))$ so that $(x_0,t)$ is within the triple-valued region for the solution of Burgers' equation with initial data $u_0$.
Let $k$ be the index of a small eigenvalue $|\alpha_k(x,t)|\leq \epsilon^r$. For $|x-x_0|\leq \epsilon^{q}$ and $x\not\in X_o$, the following estimates hold, where the remainder terms are uniform over all indices $|k|\leq C_\pm \epsilon^{r-1}$. 

\begin{enumerate}
\item  One can write 
\begin{equation}\label{eq:approx-alpha-fast}
\alpha_k^{\pm}(x,t)=\frac{2}{h_k^{\pm}}(x_0,t;\epsilon)(x f^{\pm}(x_0,t)-2\pi k\varepsilon)+R_k^{\pm}(x,t),
\end{equation}
\begin{equation}\label{eq:approx-alpha_x-fast}
\partial_x \alpha_k^{\pm}(x,t)=\frac{2{f^\pm}(x_0,t)}{h_k^\pm(x_0,t;\epsilon})+Q_k^\pm(x,t),
\end{equation}
\begin{equation}\label{eq:Rk-Qk}
R_k^\pm(x,t)= \frac{2r^\pm(x_0,t)f^\pm(x_0,t)}{h_k^\pm(x_0,t;\epsilon)}\epsilon
+ o(\epsilon^q),
\quad Q_k^\pm(x,t)=o(1).
\end{equation}
\item (Fast eigenvalues.) For some bounded $|c^+_k(x_0,t;\epsilon)|\leq C$,
\begin{equation}\label{eq:g-fast-rough}
h_k^+(x_0,t;\epsilon)=\ln\left(\frac{u_2^\mathrm{B}(x_0,t)}{u_1^\mathrm{B}(x_0,t)}\right)
	+c^+_k(x_0,t;\epsilon).
\end{equation}
\item (Slow eigenvalues.) For some bounded $|c^-_k(x_0,t;\epsilon)|\leq C$,

\begin{equation}\label{eq:g-slow}
h^-_k(x_0,t;\epsilon)=c^-_k(x_0,t;\epsilon)\ln(\epsilon^{-1}).
\end{equation}
\end{enumerate}
\end{proposition}

The rest of this subsection is devoted to the proof of Proposition~\ref{prop:bohr}.

\subsubsection{Spatial derivative of the small eigenvalues}

Choosing an eigenvector $\mathbf{u}=(u_{j})_j$ (which is not necessarily normalized) with eigenvalue $\alpha_l^\pm(x_l^\pm,t)$, the equality~\eqref{eq:alpha-x} leads to
\begin{equation}\label{eq:alphax pm}
\partial_x \alpha_l^\pm(x_l^\pm,t) =\frac{\displaystyle \sum_{j=1}^{N(\epsilon)}(-2\lambda_j)|u_{j}|^2}{\displaystyle \sum_{j=1}^{N(\epsilon)} |u_{j}|^2}=\frac{\displaystyle\frac{1}{N(\epsilon)}\sum_{j=1}^{N(\epsilon)}(-2\lambda_j)|u_j|^2}{\displaystyle\frac{1}{N(\epsilon)}\sum_{j=1}^{N(\epsilon)}|u_j|^2}.
\end{equation}
We know that the leading order of the amplitude $|u_{l,j}^\pm(x_l^\pm,t)|$ is a discretization of a constant multiple (which we take to be $1$ for the purposes of this computation as it will cancel between the numerator and denominator of $\partial_x\alpha_l^\pm$) of~\eqref{eq:amplitude}
\begin{equation}
a(y)=\sqrt{\frac{\Lambda'(y)}{-2\Lambda(y)}}
\label{eq:amplitude-again}
\end{equation}
in the range $y\in[y_{\min}(x_l^\pm,t),y_{\max}(x_l^\pm,t)]$, where we set $y_{\min}(x_l^\pm,t)=Y(\lambda_{\min}(x_l^\pm,t))$ and $y_{\max}(x_l^\pm,t)=Y(\lambda_{\max}(x_l^\pm,t))$.  Given that $x_l^\pm$ is close to $x_0$, one can replace $x_l^\pm$ by $x_0$  in the bounds of integration, up to a remainder term of order $\epsilon^r$.
This implies
\begin{equation}\label{eq:num}
\lim_{\varepsilon\to 0}
\frac{1}{N(\epsilon)}\sum_{j=1}^{N(\epsilon)}(-2\lambda_j)|u_{l,j}^{\pm}(x_l^\pm,t)|^2
	=\int_{y_{\min}(x_0,t)}^{y_{\max}(x_0,t)}(-2\Lambda(y))|a(y)|^2\dd y
	\in[0,+\infty],
\end{equation}
\begin{equation}\label{eq:denom}
\lim_{\varepsilon\to 0}
\frac{1}{N(\epsilon)}\sum_{j=1}^{N(\epsilon)}|u_{l,j}^\pm(x_l^\pm,t)|^2
	=\int_{y_{\min}(x_0,t)}^{y_{\max}(x_0,t)}|a(y)|^2\dd y\in[0,+\infty].
\end{equation}

We first show that the numerator of \eqref{eq:alphax pm} always has a finite limit. Indeed, we compute using \eqref{eq:amplitude-again} and the change of variables $\lambda=\Lambda(y)$
\begin{equation}
\begin{split}
\int_{y_{\min}(x_0,t)}^{y_{\max}(x_0,t)}(-2\Lambda(y))|a(y)|^2\dd y
	&=\int_{y_{\min}(x_0,t)}^{y_{\max}(x_0,t)}\Lambda'(y)\,\dd y \\
 &=\lambda_{\max}(x_0,t)-\lambda_{\min}(x_0,t).
 \end{split}
\end{equation}
In the case of fast eigenvalues, one has to sum up the contributions of the main part and of the correction. Since both integrals are finite, the slow correction is a negligible remainder term in the case $x_l^\pm=x_l^+$. Hence
\begin{align}
\frac{1}{N(\epsilon)}\sum_{j=1}^{N(\epsilon)}(-2\lambda_j)|u_{l,j}^{+}(x_l^+,t)|^2
	&=u^\mathrm{B}_2(x_0,t)-u^\mathrm{B}_1(x_0,t)+o(1), \\
\frac{1}{N(\epsilon)}\sum_{j=1}^{N(\epsilon)}(-2\lambda_j)|u_{l,j}^{-}(x_l^-,t)|^2
	&=u^\mathrm{B}_0(x_0,t)+o(1).
\end{align}
In both cases, we retrieve $f^\pm(x_0,t)$.

Let us now study the denominator of \eqref{eq:alphax pm}.  Again using \eqref{eq:amplitude-again} and $\lambda=\Lambda(y)$,
\begin{equation}
\int_{y_{\min}(x_0,t)}^{y_{\max}(x_0,t)}|a(y)|^2\dd y
	=
 \int_{y_{\min}(x_0,t)}^{y_{\max}(x_0,t)}\frac{\Lambda'(y)\,\dd y}{-2\Lambda(y)}
	=
 \frac{1}{2}\left(\ln(-\lambda_{\min}(x_0,t))-\ln(-\lambda_{\max}(x_0,t))\right).
\end{equation}
The integral is finite and equal to $\frac 12 g^+(x_0,t)$ in the fast case that $\lambda_{\min}(x_0,t)=-u^\mathrm{B}_2(x_0,t)$ and $\lambda_{\max}(x_0,t)=-u^\mathrm{B}_1(x_0,t)$, but infinite in the slow case because $\lambda_{\max}(x_0,t)=0$.

Let us now use the assumption that there is $p>\frac 12$ such that $u_0(x)\sim Cx^{-2p}$ as $|x|\to+\infty$.  We focus on the values of $\lambda_j$ which are close to $0$, i.e. for which $j$ is close to $N=N(\epsilon)$. Using formula~\eqref{def: tilde lambda n}, and assuming the condition \eqref{eq:result-of-epsilon-quantization}, we have $\lambda_{j}=\Lambda(y_{j})$, so
\begin{equation}
\int_{\lambda_{N+1-j}}^0F(\lambda)\dd \lambda
	=\epsilon \left(j-\tfrac{1}{2}\right). 
\end{equation}
Since $F(\lambda)\sim (-C\lambda)^{-1/(2p)}$ as $\lambda\to 0$ by assumption,  we get that when $\lambda\to 0$,
\begin{equation}
\int_{\lambda}^0F(\eta)\dd \eta\sim \frac{C^{1/(2p)}(-\lambda)^{1-1/(2p)}}{1-1/(2p)},
\end{equation}
so that for $j$ close to $1$,
\begin{equation}
\Lambda(y_{N+1-j})\sim  -\left(\frac{2p-1}{C^{1/(2p)}2p} \epsilon (j-\tfrac{1}{2})\right)^{2p/(2p-1)}.
\label{eq:smallest-lambdas}
\end{equation}
Hence $\Lambda'(y_{N+1-j})/\Lambda(y_{N+1-j})
	=[\Lambda(y_{N+1-j}) F(\Lambda(y_{N+1-j}))]^{-1}$ has the asymptotic expansion
\begin{equation}
-\frac{\Lambda'(y_{N+1-j})}{\Lambda(y_{N+1-j})}
	\sim \frac{1}{C^{-1/(2p)}(-\Lambda(y_{N+1-j}))^{1-1/(2p)}}\sim \frac{C^{1/p}2p}{2p-1}\frac{1}{\epsilon (j-\frac{1}{2})}.
\end{equation}
This approximation is valid for $\epsilon j\lesssim 1/\ln(\epsilon^{-1})$ since in this case we can check that uniformly in this range, $|\Lambda(y_{N+1-j})|=o(1)$. Hence the denominator sum has the lower bound
\begin{equation}\label{eq:lower}
\frac{1}{N(\epsilon)} \sum_{j=1}^{N(\epsilon)}|u_{l,j}^-|^2
	\gtrsim \frac{1}{N(\epsilon)}\sum_{j=1}^{\mathcal{O}(1/\epsilon\ln(\epsilon^{-1}))}\frac{\Lambda'(y_{N+1-j})}{-2\Lambda(y_{N+1-j})}
\gtrsim \frac{1}{\epsilon N(\epsilon)}\sum_{j=1}^{\mathcal{O}(1/\epsilon\ln(\epsilon^{-1}))}\frac{1}{j-\frac{1}{2}}
	\gtrsim \ln(\epsilon^{-1}).
\end{equation}
To get an upper bound, we split between the cases $\epsilon j\lesssim   1/\ln(\epsilon^{-1})$ and $\epsilon j\gtrsim   1/\ln(\epsilon^{-1})$. The small indices $j$ such that $\epsilon j\lesssim   1/\ln(\epsilon^{-1})$  are treated as in~\eqref{eq:lower}. Concerning the indices $j$ such that $\epsilon j\gtrsim   1/\ln(\epsilon^{-1})$, we use the fact that $t>t_b$ and that $\mathbf{u}_l^-$ is supported on $[-u_0^\mathrm{B}(x_l^+,t),0]$, that is, away from the left edge $\lambda=-L$.  Then, the maximum of $|\Lambda'(y_{N+1-j})/\Lambda(y_{N+1-j})|$ for $\epsilon j\gtrsim 1/\ln(\epsilon^{-1})$ occurs for the smallest $j$, and hence  $\Lambda'(y_{N+1-j})/\Lambda(y_{N+1-j})
=\mathcal{O}(\ln(\epsilon^{-1}))$ holds uniformly on this index range because $\Lambda'(y)/\Lambda(y)\sim (M-y)^{-1}$ for $y\approx M$. This yields the upper bound
\begin{equation}
\frac{1}{N(\epsilon)} \sum_{j=1}^{N(\epsilon)}|u_{l,j}^-|^2
	\lesssim \frac{1}{N(\epsilon)}\sum_{j=1}^{N(\epsilon)}
 \frac{\Lambda'(y_{N+1-j})}{-2\Lambda(y_{N+1-j})}
	\lesssim \ln(\epsilon^{-1}).
\end{equation}
Consequently, the sum over the ``slow'' part of the eigenvectors satisfies
\begin{equation}
 \frac{1}{N(\epsilon)}\sum_{j=1}^{N(\epsilon)}|u_{l,j}^-(x_l^\pm,t)|^2=c_l^-(x_0,t;\epsilon)\ln(\epsilon^{-1})
 \label{eq:386}
\end{equation}
for some bounded constant $ c_l^-(x_0,t;\epsilon)$ that may be oscillatory as $\epsilon\to 0$.
Regarding the fast eigenvectors, we sum up the contributions $u_{l,j}^+(x_l^+)$ and $u_{l,j}^-(x_l^+)$, which are supported in disjoint intervals,
\begin{equation}
\frac{1}{N(\epsilon)}\sum_{j=1}^{N(\epsilon)}|u_{l,j}(x_l^+,t)|^2
    = \frac{1}{2}\ln\left(\frac{u^{\mathrm{B}}_2(x_0,t)}{u^{\mathrm{B}}_1(x_0,t)}\right)+  \frac{c}{\ln(\epsilon^{-1})} c_l^-(x_0,t;\epsilon) \ln(\epsilon^{-1})
    +o(1).
\end{equation}
We choose $c_l^+(x_0,t;\epsilon)=c\cdot c_l^-(x_0,t;\epsilon)$ to get the Proposition.

We conclude that  spatial derivatives $\partial_x\alpha_l^\pm(x_l^\pm)$ of the fast and slow eigenvalues satisfy the following identities at $x_l^\pm$:
\begin{align}
\partial_x \alpha_l^+(x_l^+,t)
    &=\frac{2(u^{\mathrm{B}}_2(x_0,t)-u^{\mathrm{B}}_1(x_0,t))+o(1)}{\ln\left(\displaystyle\frac{u^{\mathrm{B}}_2(x_0,t)}{u^{\mathrm{B}}_1(x_0,t)}\right)+c_l^+(x_0,t;\epsilon)}+o(1), \\
\partial_x \alpha_l^- (x_l^-,t)
	&=\frac{2u^{\mathrm{B}}_0(x_0,t)+o(1)}{c_l^-(x_0,t;\epsilon)}\ln(\epsilon^{-1})(1+o(1)).
\end{align}
We deduce that~\eqref{eq:approx-alpha_x-fast} holds.  Finally, we write
\begin{equation}
\alpha_{l}^{\pm}(x,t)
	=\alpha_l^\pm(x_l^\pm,t)
+(x-x_l^\pm) \partial_x \alpha_l^\pm(x_l^\pm,t)
+o(\epsilon^{q}).
\end{equation}
Then we use~\eqref{eq:xl} and the estimate $|\alpha_l^\pm(x_l^\pm,t)|\leq C\epsilon^{1+r}$ to deduce that~\eqref{eq:approx-alpha-fast} holds.

\subsection{Small-\texorpdfstring{$\epsilon$}{epsilon} asymptotics of the sum}

Assuming now that a stronger form of Conjecture~\ref{conjecture:Bohr} holds, we establish an asymptotic expansion of the soliton ensemble solution to~\eqref{BO equation} for $(x,t)$ near $(x_0,t)$ in the oscillatory region (triple-valued region for the solution of Burgers' equation with data $u_0$).

We retrieve a result similar to Theorem~\ref{thm:oscillations} above in the oscillatory region.

\begin{theorem}\label{prop: approx after breaking via A}
Assume that Conjecture~\ref{conjecture:Bohr} holds. Fix $t>t_b$ and $x_0\in (X^-(t),X^+(t))$ so that $(x_0,t)$ is within the triple-valued region for the solution of Burgers' equation with initial data $u_0$. Let $|x-x_0|\leq \epsilon^{q}$ and $x\not\in X_o$. We assume moreover that uniformly in $k$, one can write the remainder terms from~\eqref{eq:Rk-Qk} as
\begin{equation}
R_k^\pm(x,t)=\epsilon\varphi^\pm(x_0,t)+o(\epsilon^{2-r}),
\quad
Q_k^\pm(x,t)=o(\epsilon^{1-r}),
\end{equation}
and that in~\eqref{eq:g-fast-rough} and~\eqref{eq:g-slow}, $c^\pm_k(x_0,t;\epsilon)=c^\pm(x_0,t)$ do not depend on $k$ nor on $\epsilon$.
Then
\begin{equation}
u(x,t)
	=u_0^\mathrm{B}(x_0,t) 
	+f^+(x_0,t)\frac{\sinh(\frac{1}{2}h^+(x_0,t))}{\cosh(\frac{1}{2}h^+(x_0,t))-\cos( f^+(x_0,t)\epsilon^{-1}x+\varphi^+(x_0,t))}+o(1).
 \label{eq:u-oscillatory-A}
\end{equation}
If moreover, in~\eqref{eq:g-fast-rough}, there holds
\begin{equation}\label{eq:g-fast}
h^+(x_0,t)=\ln\left(\frac{u_2^\mathrm{B}(x_0,t)-u_0^\mathrm{B}(x_0,t)}{u_1^\mathrm{B}(x_0,t)-u_0^\mathrm{B}(x_0,t)}\right),
\end{equation}
then
\begin{equation}
u(x,t)
	=u_0^\mathrm{B}(x_0,t)+\frac{(u_2^\mathrm{B}(x_0,t)-u_1^\mathrm{B}(x_0,t))(1-r(x_0,t)^2)}{1+r(x_0,t)^2-2r(x_0,t)\cos(\Theta(\epsilon^{-1}x;x_0,t))}+o(1),
 \label{eq:u-A-alt}
\end{equation}
where
\begin{equation}
\Theta(z;x_0,t)
    =(u_2^\mathrm{B}(x_0,t)-u_1^\mathrm{B}(x_0,t))z+\varphi^+(x_0,t),
\end{equation}
\begin{equation}
r(x_0,t)
	=\sqrt{\frac{u_1^\mathrm{B}(x_0,t)-u_0^\mathrm{B}(x_0,t)}{u_2^\mathrm{B}(x_0,t)-u_0^\mathrm{B}(x_0,t)}}\in (0,1).
\end{equation}
\end{theorem}

\begin{remark}The functions $\varphi^+$ and $\varphi^-$ were not present in Proposition~\ref{prop:bohr}. This is likely because we only considered the principal part of the Bohr-Sommerfeld condition in Conjecture~\ref{conjecture:Bohr}, instead of all of its semiclassical expansion. If the symbol $p$ were smooth, then the complete asymptotic expansion of the Bohr-Sommerfeld condition would be given by~\cite[Therorem 3.1]{Charles2003quasimodes}, in which the first two terms are explicitly written, and the principal part is precisely the parallel-transport integral $g_0^\pm(x,t)$ (see \eqref{eq:g0plus}--\eqref{eq:g0minus}).
\end{remark}

\begin{remark}
    Aside from details of the phase correction $\varphi^+$, the formula \eqref{eq:u-oscillatory-A} coincides with \eqref{eq: approx after breaking} from Theorem~\ref{thm:oscillations} upon proper identification of the periodic wave parameters.  Note that if one averages \eqref{eq:u-A-alt} over the fast variable $\Theta$ (requiring an integration since the cosine appears in the denominator), then the result is exactly as expected, namely the weak limit $\overline{u}(x_0,t):=u_2^\mathrm{B}(x_0,t)-u_1^\mathrm{B}(x_0,t)+u_0^\mathrm{B}(x_0,t)$ defined in \eqref{eq:ubar}.
\end{remark}

\begin{proof}
We start from the formula~\eqref{u in terms of alphas}
\begin{equation}
u(x,t)
 = \sum_{k=1}^{N(\epsilon)}\frac{2\partial_x\alpha_{k}(x,t)}{\left(\epsilon^{-1}\alpha_k(x,t)\right)^2 + 1}.
\end{equation}
According to~\eqref{eq:large-A-eigenvalues}, it is enough to focus on the small eigenvalues $\alpha$ such that $|\alpha_k|\leq \varepsilon^{r}$.  Each such eigenvalue is either fast, slow, or ambiguous, and the total number of each type in the indicated range of $\alpha$ is proportional to $\epsilon^{r-1}$.
  
First, let us estimate the contribution of the ambiguous eigenvalues. Using~\eqref{eq:bound_alpha_x} and the lower bound $|\alpha_k^o(x,t)|\gtrsim k \epsilon^s$ because $x\not\in X_o$, we have
\begin{equation}
\sum_{|k|\leq C_o\varepsilon^{r-1}} \frac{2\partial_x\alpha_{k}^o(x,t)}{\left(\epsilon^{-1}\alpha_k^o(x,t)\right)^2 + 1}
	\lesssim \sum_{|k|\leq C_o\varepsilon^{r-1}} \frac{1}{(\epsilon^{s-1})^2 + 1}
	\lesssim \epsilon^{1+r-2s}.
\end{equation} 
Given that $\frac 13<r<\frac 12$ by assumption, one can choose $\frac 12<s<1$ close enough to $\frac 12$ so that $1+r-2s>0$, and this sum goes to zero as $\epsilon\to 0$. Since the total contribution of the ambiguous eigenvalues in the range $|\alpha_k|\le\epsilon^r$ is negligible, it is harmless to ``double count'' them as both fast and slow eigenvalues and as this makes the arguments easier going forward we shall do so.

Then, let us tackle the sum over fast eigenvalues.  The ideas are remarkably similar to proof of Proposition \ref{prop: lower branch leading order}, so we omit some of the details. By assumption, we have that
\begin{multline}
\sum_{|k|\leq C_+\varepsilon^{r-1}} \frac{2\partial_x\alpha_{k}^+(x,t)}{\left(\epsilon^{-1}\alpha_k^+(x,t)\right)^2 + 1}
 \\=\sum_{|k|\le C_+\varepsilon^{r-1}}\frac{f^+(x_0,t)h^+(x_0,t)}{(-2\pi k + f^+(x_0,t)\epsilon^{-1}x+\varphi^+(x_0,t))^2 + \frac{1}{4}h^+(x_0,t)^2} + o(1).
\end{multline}
In the limit $\epsilon\to 0$, we have $\varepsilon^{r-1}\to+\infty$, so the sum on the right-hand side can be computed:
\begin{equation}\label{eq:oscillation}
\sum_{|k|\leq C_+\epsilon^{r-1}} \frac{2\partial_x\alpha_{k}^+(x,t)}{\left(\epsilon^{-1}\alpha_k^+(x,t)\right)^2 + 1}
	= \frac{f^+(x_0,t)\sinh(\frac{1}{2}h^+(x_0,t))}{\cosh(\frac{1}{2}h^+(x_0,t))-\cos( f^+(x_0,t)\epsilon^{-1}x+\varphi^+(x_0,t))}+o(1).
\end{equation}
If the expression of $h^+(x_0,t)$ is given by~\eqref{eq:g-fast}, one can further simplify
\begin{align}
\sinh\left( \frac{1}{2}h^+(x,t)\right)	
	&=\frac{u_2^\mathrm{B}(x,t)-u_1^\mathrm{B}(x,t)}{2\sqrt{(u_2^\mathrm{B}(x,t)-u_0^\mathrm{B}(x,t))(u_1^\mathrm{B}(x,t)-u_0^\mathrm{B}(x,t))}}, \\
\cosh\left( \frac{1}{2}h^+(x,t)\right)	
	&=\frac{u_2^\mathrm{B}(x,t)+u_1^\mathrm{B}(x,t)-2u^\mathrm{B}_0(x,t)}{2\sqrt{(u_2^\mathrm{B}(x,t)-u_0^\mathrm{B}(x,t))(u_1^\mathrm{B}(x,t)-u_0^\mathrm{B}(x,t))}}.
\end{align}
Hence,
\begin{equation}\label{eq:oscillation2}
\sum_{|k|\leq C_+\epsilon^{r-1}} \frac{2\partial_x\alpha_{k}^+(x,t)}{\left(\epsilon^{-1}\alpha_k^+(x,t)\right)^2 + 1}
	=\frac{(u_2^\mathrm{B}(x_0,t)-u_1^\mathrm{B}(x_0,t))(1-r(x_0,t)^2)}{1+r(x_0,t)^2-2r(x_0,t)\cos(\Theta(\epsilon^{-1}x;x_0,t))}+o(1).
\end{equation}

Now we consider the sum over slow eigenvalues.  Again, a parallel can be made to the previous section with Proposition \ref{prop: upper branch leading order}, where we obtained the leading order behavior of a sum by turning it into an integral.  Given~\eqref{eq:f-slow} of Proposition~\ref{prop:bohr}, we have that the $x$-velocities of all slow and small eigenvalues for $(x,t)$ near $(x_0,t)$ are nearly the same, i.e.
\begin{equation}
    \partial_x\alpha_k^-(x,t)\approx \frac{2 u_0^\mathrm{B}(x_0,t)}{h^-(x_0,t)\ln(\epsilon^{-1})}.
\end{equation}
However, at this point we need to approximate the quantity $h^-(x_0,t)$.  We will actually proceed more directly to approximate $\partial_x\alpha_k^-(x,t)$, which will automatically produce the logarithmic scaling in $\epsilon$.  For this purpose, let $x_-(t)$ denote the largest value of $x$ less than or equal to $x_0$ for which there is a slow eigenvalue at $\alpha=0$, let $x_+(t)$ denote the smallest value of $x$ strictly greater than $x_0$ for which the same is true, and set $\Delta x:=x_+(t)-x_-(t)$.
By setting the left-hand side of \eqref{eq:approx-alpha-fast} to zero we obtain 
\begin{equation}
    \Delta x\approx \frac{2\pi\epsilon}{f^-(x_0,t)} = \frac{2\pi\epsilon}{u_0^\mathrm{B}(x_0,t)}.
\end{equation}
Now as $x$ increases from $x_-(t)$ to $x_0$, the slow eigenvalue originally at $\alpha=0$ will increase to the value of the smallest positive slow eigenvalue $\alpha_+(x_0,t)$; likewise as $x$ decreases from $x_+(t)$ to $x_0$, the slow eigenvalue originally at $\alpha=0$ will decrease to the value of the negative slow eigenvalue $\alpha_-(x_0,t)$ of smallest absolute value.  We set $\Delta\alpha:=\alpha_+(x_0,t)-\alpha_-(x_0,t)$, and see that this is the difference between the two closest eigenvalues of opposite signs to $\alpha=0$ at $(x_0,t)$.  We can approximate $\Delta\alpha$ directly using the asymptotic 
density of all of the eigenvalues given in \eqref{eq:distribution-of-alphas}.  Indeed, the integral of $G(\alpha;x_0,t)\,\dd\alpha$ between $\alpha_-(x_0,t)$ and $\alpha_+(x_0,t)$ should be approximately 
\begin{equation}
\frac{M}{N(\epsilon)}=\epsilon(1+o(1)).
\end{equation}  Using \eqref{eq:G-small-alpha}, since the limits of integration are small, we get
\begin{multline}
\frac{1}{4\pi}\frac{2p}{2p-1}\left[\alpha_+(x_0,t)-\alpha_+(x_0,t)\ln(|\alpha_+(x_0,t)|)\right.\\
\left.{}-\alpha_-(x_0,t)+\alpha_-(x_0,t)\ln(|\alpha_-(x_0,t)|)\right]\approx \epsilon.
\end{multline}
Solving for $\Delta\alpha$ gives 
\begin{equation}
    \Delta\alpha \approx \frac{2\pi\epsilon (2p-1)}{p\ln(\epsilon^{-1})}.
\end{equation}  
Then we approximate the velocity $\partial_x\alpha(x,t)$ for a slow and small eigenvalue by 
\begin{equation}
    \partial_x\alpha(x,t)\approx \frac{\Delta\alpha}{\Delta x}\approx \frac{(2p-1)u_0^\mathrm{B}(x_0,t)}{p\ln(\epsilon^{-1})} \implies g^-(x_0,t)\approx \frac{2p}{2p-1}.
\end{equation}
Using this and the convergence in \eqref{eq:distribution-of-alphas}, the contribution of the small and slow eigenvalues (which are the majority of the small eigenvalues), is
\begin{equation}
\begin{split}
\sum_{|k|\leq C_-\epsilon^{r-1}}
\frac{2\alpha^-_{k,x}(x,t)}{(\epsilon^{-1}\alpha^-_k(x,t))^2+1} &
	\approx u_0^\mathrm{B}(x_0,t)\frac{2(2p-1)}{p\ln(\epsilon^{-1})} \sum_{|k|\leq C_-\epsilon^{r-1}}\frac{1}{(\epsilon^{-1}\alpha_k^-(x_0,t))^2+1}
\\
& = u_0^\mathrm{B}(x_0,t)\frac{2(2p-1)}{p\ln(\epsilon^{-1})}\int_{-\epsilon^{r}}^{\epsilon^r}\frac{1}{(\epsilon^{-1}\alpha)^2+1}\sum_{k=1}^{N(\epsilon)}\delta_{\alpha_k^-}(x_0,t)(\alpha)\\
&\approx u_0^\mathrm{B}(x_0,t) \frac{2(2p-1)}{p\ln(\epsilon^{-1})}\cdot\frac{1}{\epsilon}\int_{-\epsilon^{r}}^{\epsilon^r}\frac{G(\alpha;x_0,t)\,\dd\alpha}{(\epsilon^{-1}\alpha)^2+1}\\
&\approx \frac{u_0^\mathrm{B}(x_0,t)}{\pi\epsilon\ln(\epsilon^{-1})}\int_{-\epsilon^r}^{\epsilon^r}\frac{\ln(|\alpha|^{-1})\,\dd\alpha}{(\epsilon^{-1}\alpha)^2+1}.
\end{split}
\end{equation}
Scaling by $\alpha=\epsilon z$ and letting $\epsilon\to 0$, the integral above is $\pi\epsilon\ln(\epsilon^{-1})(1+o(1))$, so we conclude that 
the leading contribution of the slow and small eigenvalues is
\begin{equation}
   \lim_{\epsilon\to 0} \sum_{|k|\leq C_-\epsilon^{r-1}}
\frac{2\alpha_{k,x}^-(x,t)}{(\epsilon^{-1}\alpha_k^-(x,t))^2+1}=u_0^\mathrm{B}(x_0,t).
\end{equation}

\end{proof}

\begin{remark}
The proof above shows an interesting connection between the approaches to the strong small-$\epsilon$ asymptotic behavior of the BO soliton ensemble $u(x,t)$ based on the eigenvalues of the two matrices $\mathbf{A}(x,t)$ (the approach in this section) and $\mathbf{C}(t)$ (see Section~\ref{sec:C}).  Indeed, we can see that the contributions of the small slow/fast  eigenvalues of $\mathbf{A}(x,t)$ respectively correspond precisely to the contributions of eigenvalues of $\mathbf{C}(t)$ with real parts close to $\mu=x$ on the upper/lower branch.
\end{remark}

\begin{remark}
    Although for the proof we needed to exclude extremely small ambiguous eigenvalues, we expect that a more accurate modeling of the near-crossings of the eigenvalue curves as seen in Figure~\ref{fig:alphas} would allow such eigenvalues to be included without changing the leading asymptotics.
\end{remark}

\section*{Acknowledgments}
E. Blackstone was partially supported by the National Science Foundation under grant DMS-1812625.  L. Gassot conducted this work within the France 2030 framework program, the Centre Henri Lebesgue  ANR-11-LABX-0020-01.  P. D. Miller was partially supported by the National Science Foundation under grants DMS-1812625 and DMS-2204896, and some of this work was done with the support of a Leverhulme Trust Visiting Professorship at Bristol University, UK.  

The authors wish to thank  Alejandro Uribe, Alix Deleporte, and San V\~{u} Ng\d{o}c for useful discussions.

\bibliographystyle{siamplain}
\bibliography{references}

\end{document}